\documentclass[11pt]{amsart}

\usepackage{mathtools}
\usepackage{amsmath, amsthm, amssymb, amsfonts, enumitem}
\usepackage{dsfont}
\usepackage{color}
\usepackage{geometry}
\usepackage{todonotes}
\usepackage{soul}
\usepackage[normalem]{ulem}

\usepackage[colorlinks=true,linkcolor=blue,urlcolor=blue, citecolor=blue]{hyperref}
\usepackage{comment}

\mathtoolsset{showonlyrefs}

\def \cA{\mathcal{A}}
\def \cAcirc{\mathcal{A}^{\raise1pt\hbox{${\scriptstyle\circ}$}}}

\def \cC{\mathcal{C}}
\def \cD{\mathcal{D}}

\def \cF{\mathcal{F}}
\def \cG{\mathcal{G}}
\def \cH{\mathcal{H}}

\def \cL{\mathcal{L}}
\def \cM{\mathcal{M}}

\def \cO{\mathcal{O}}
\def \cP{\mathcal P}

\def \cT{\mathcal{T}}
\def \P{\mathsf P}
\def \hP{{\mathbb P}}
\def \hE{{\mathbb E}}

\def \E{\mathsf E}

\def \N{\mathbb{N}}
\def \bH{\mathbb{H}}
\def \R{\mathbb{R}}
\def \F{\mathbb F}
\def \G{\mathbb G}

\def \ud{\mathrm{d}}
\def \e{\mathrm{e}}

\def \tilu{u}
\def \tilv{v}
\def \Du0{U}

\newcommand\optional[1]{\tensor*[^o]{{\kern-1pt #1}}{}}
\newcommand\previsible[1]{\tensor*[^p]{{\kern-1pt #1}}{}}

\def \HYe{{T^Y_\e}}
\def \HY{T^Y}

\def \He{{T_\e}}
\def \H{T}
\def \oldh{\phi}
\def \oldtlh{\Phi}

\makeatletter
\def\punctuationSpace{\@ifnextchar.{}{\@ifnextchar,{}{\@ifnextchar;{}{ }}}}
\makeatother

\newcommand{\cadlag}{c\`adl\`ag\punctuationSpace}

\newcommand\as{\mbox{-a.s.}\punctuationSpace}

\newcommand{\eps}{\varepsilon}
\newcommand{\ind}{\mathbf{1}}

\newtheorem{theorem}{Theorem}[section]
\newtheorem{lemma}[theorem]{Lemma}
\newtheorem{corollary}[theorem]{Corollary}
\newtheorem{proposition}[theorem]{Proposition}

\newtheorem{definition}[theorem]{Definition}
\newtheorem{remark}[theorem]{Remark}
\newtheorem{assumption}[theorem]{Assumption}

\theoremstyle{definition}

\definecolor{ballblue}{rgb}{0.63, 0.37, 0.9}
\definecolor{darkgreen}{rgb}{0.2, 0.59, 0.3}

\allowdisplaybreaks

\makeatletter
\@namedef{subjclassname@2020}{\textup{2020}Mathematics Subject Classification}
\makeatother

\title[A stopping game with foresight]{Perfect Bayesian equilibrium \\ in a zero-sum stopping game with asymmetric foresight}
\thanks{{\bf Acknowledgements}: J. Palczewski was partially supported by EPSRC grant no.\ UKRI175. T.\ De Angelis was supported by – Next Generation EU – PRIN2022 (2022BEMMLZ) CUP: D53D23005780006 and PRIN-PNRR2022 (P20224TM7Z) CUP: D53D23018780001. T.\ De Angelis thanks D.\ Gerardi and I.\ Monzon for insightful discussions around perfect Bayesian equilibria.}
\author[De Angelis]{Tiziano De Angelis}
\author[Hobson]{David Hobson}
\author[Palczewski]{Jan Palczewski}
\subjclass[2020]{91A27, 91A55, 91A15, 60G07, 60G40}
\keywords{}
\address{T.\ De Angelis: School of Management and Economics, Dept.\ ESOMAS, University of Torino, Corso Unione Sovietica, 218 Bis,
10134, Torino, Italy; Collegio Carlo Alberto, Piazza Arbarello 8, 10122, Torino, Italy.}
\email{\href{mailto:tiziano.deangelis@unito.it}{tiziano.deangelis@unito.it}}
\address{D.\ Hobson: Department of Statistics, University of Warwick, Coventry, CV4 7AL, UK.}
\email{\href{mailto:d.hobson@warwick.ac.uk}{d.hobson@warwick.ac.uk}}
\address{J.\ Palczewski: Faculty of Mathematics, Wroc\l{}aw University of Science and Technology, Wybrze\.{z}e Wyspia\'{n}skiego
27,
50-370, Wroc\l{}aw, Poland.}
\email{\href{mailto:jan.palczewski@pwr.edu.pl}{jan.palczewski@pwr.edu.pl}}

\date{\today}
\numberwithin{equation}{section}

\begin{document}

\begin{abstract}
This paper provides closed-form saddle points for a novel class of continuous-time, zero-sum stopping games with one-dimensional diffusive dynamics  and asymmetric information across players. The novelty arises because although both players observe the underlying process $X$, an informed player also holds information about {\em a future state} of $X$, whereas the other player does not. We extend the notion of perfect Bayesian equilibrium (PBE) to the setting of asymmetric information in continuous time and we construct saddle points as PBEs. The methodology is also novel and PBEs are obtained by introducing the running maximum of the process $X$ as an auxiliary process but with an endogenously determined initial condition. 
That, combined with the use of randomised stopping times, yields computable formulae for the equilibrium stopping times of both players. Explicit solutions are obtained in some examples.
\end{abstract}

\maketitle

\section{Introduction}

Consider a zero-sum stopping game (Dynkin game \cite{dynkin1969}) whose underlying process $X$ is a Brownian motion killed upon leaving the interval $(0,1)$ and with payoff functions granting a second mover advantage structure for the game (i.e., at each time each player would rather let her opponent stop first). 
When both players only observe the dynamics of $X$ this game is a standard Dynkin game with one-dimensional state space and it can be addressed using a variety of methods available in the literature: see for example Ekstr\"om and Peskir \cite{ekstrom2008} for the general theory of Markovian Dynkin games, Ekstr\"om and Villeneuve \cite{ekstrom2006} for solutions of one-dimensional diffusive games, and Peskir \cite{peskir2009optimal} for the semi-harmonic characterisation of the value function for general Markov processes. The solution is generally in closed-form or explicit, thanks to the one-dimensional nature of the problem and depending on the specific assumptions on the payoff functions. Nash equilibria (saddle points) can be determined as stopping times for the Brownian motion and often as first passage times to constant thresholds.

In this paper we study a new situation by assuming that the minimiser can foresee the future and knows whether the Brownian path will exit $(0,1)$ from the upper or the lower endpoint. In mathematical terms we say that the maximiser only observes the filtration generated by $X$ whereas the minimiser also has foreknowledge of the random variable $X_{\He}\in\{0,1\}$, where $\He=\inf\{t\ge 0: X_t\notin(0,1)\}$. To the best of our knowledge, this is the first Dynkin game in which the informational advantage of one player is foresight on the commonly observed state process itself, rather than knowledge of a parameter (payoffs, drift, costs) independent of the driving noise. As a consequence, the two players observe the same trajectories of $X$ but assess its future under different laws: a Brownian motion for the uninformed player and Brownian motions conditioned on the two possible values of $X_{\He}$ known to the informed player. As in many games with asymmetric information, the informed player must randomise and the state must be augmented with the uninformed player's belief about $X_{\He}$; here, however, both ingredients take a form which has not been encountered in the literature. The minimiser who knows that $X_{\He}=1$ randomises stopping over times at which $X$ coincides with the running maximum: in the scenarios where the informed minimiser does not stop, the uninformed maximiser's belief that $X_{\He}=0$ is reinforced, inducing the maximiser to stop earlier. In equilibrium the belief is therefore a deterministic function of $X$ and its running maximum $M$, which turns an otherwise degenerate two-dimensional filtering problem into a tractable free-boundary problem in $(X,M)$. We demonstrate that the foresight is genuinely valuable. For instance, in a game with linear payoffs that is trivial under symmetric information, the foresight enjoyed by the informed player strictly lowers the value (Remark \ref{rem:3.6} and Section \ref{subsec:linear payoffs}).

We introduce a notion of perfect Bayesian equilibrium (PBE) for continuous-time Dynkin games with asymmetric information leveraging notions of Markov perfect equilibria developed for full information Dynkin games by Christensen and Lindensj\"o \cite{christensen2024general} and D\'ecamps, Gensbittel and Mariotti \cite{decamps2022mixed}. The definition we give bridges continuous-time diffusive settings with classical game-theoretical notions (cf.\ Osborne and Rubinstein \cite[Ch.\ 12.3]{Osborne1994}). We then construct an explicit PBE for our game in which the informed player stops with a singular intensity depending on the running maximum $M$ of the Brownian motion, starting from an initial value $M_0$ which must be determined endogenously as part of the optimal stopping rule. The value of $M_0$ is a device required to match the uninformed player's initial belief with admissible paths of the PBE. In particular, when the initial belief is too high, the informed player will stop with a positive probability at time $0$ in order to lower the uninformed player's belief. The uninformed player stops when the Brownian motion $X$ drops below a moving boundary $d_*(M_t)$, depending on the current value of $M$, that must also be determined endogenously in equilibrium. For this construction we prove that the value of the game is a convex combination of the so-called continuation values for two types of the informed player (one type for each realisation of the random variable $X_{\He}$). The weights of the convex combination are determined using equations that describe the trajectories of the uninformed player's belief along the equilibrium paths. We also show that, taken together, the continuation values of the two types of the informed player and the value function of the game solve a system of free boundary problems. The latter system is used for the actual computation of the equilibrium and of the optimal strategies for both players.

The literature on continuous-time stochastic games with partial and asymmetric information has been receiving increasing attention over the past decade. In particular, the study of Dynkin games with asymmetric information was initiated by Gr\"un \cite{Grun2013} and later on pursued by Gensbittel and Gr\"un \cite{GenGrun2019}, who developed techniques based on viscosity solutions of quasi-variational inequalities combining the classical obstacle constraints from optimal stopping with a convexity constraint with respect to a state variable representing players' beliefs. De Angelis, Ekstr\"om and Glover \cite{DEG2020} studied zero-sum Dynkin games on a diffusive dynamics with partially observable drift (an informed player knows the drift whereas an uninformed opponent only has a prior over two possible values of the drift). They develop a general verification theorem and solve an explicit example with one-dimensional state space. Some nonzero-sum games have also been considered. De Angelis and Ekstr\"om \cite{de2020playing} solved a Dynkin game in which players are uncertain about the existence of a competitor in the game. Those results were recently extended by Ekstr\"om and Wang \cite{ekstrom2026dynkin} who introduced asymmetric payoffs across the players. Nonzero-sum Dynkin games with private information in continuous time have also appeared in the recent literature (cf.\ Gaitsgori and Groenwald \cite{gaitsgori2025dynkin}, Kwon and Palczewski \cite{kwon2023exit}). In the latter games each player observes her own private state and may learn about her opponents' state by observing her actions. Stochastic differential games and mixed control-stopping games with asymmetric information were considered by Cardaliaguet and Rainer \cite{cardaliaguet2009} (who inspired the methods in \cite{Grun2013} and \cite{GenGrun2019}) and Ekstr\"om, Lindensj\"o and Olofsson \cite{ekstrom2022detect}, respectively.
Finally, a paper by Esmaeeli, Imkeller and Nzengang \cite{EIN2019}, considers a Dynkin game in which the informed player's filtration is an initial enlargement of the reference filtration by a finitely-valued random variable $G$, not necessarily independent of the underlying dynamics; their example $G=\ind_{\{B_T\in[a,b]\}}$ is a form of foresight. However, the game criterion is the expected payoff conditional on $\sigma(G)$, so that the uninformed player's optimisation is also carried out separately on each event $\{G=g\}$. As a consequence the value of the game can be constructed from the value of the game in which both players observe $G$, no learning takes place and pure stopping times suffice. In contrast, in our game the uninformed player assesses payoffs under her own information, the asymmetry persists in equilibrium, and both randomisation and the belief process are indispensable.

In general the Dynkin games reviewed above were studied with ad-hoc methods developed on a case-by-case basis. Works by De Angelis, Palczewski and Merkulov \cite{de2022value} and De Angelis, Palczewski and Smith \cite{de2025nash} go in the direction of developing a general framework for zero-sum Dynkin games with arbitrary information structure. The former paper established existence of saddle points for general non-Markovian games and the latter paper obtained necessary and sufficient conditions for the dynamics of players' continuation values and optimal strategies in such equilibria. However, neither of those two papers offers explicit constructions of equilibria in general diffusive setups.
In this paper we leverage results from \cite{de2025nash} and go beyond them by deriving in a Markovian framework a system of PDEs that transfers abstract martingale conditions from \cite{de2025nash} into more tangible objects amenable to further analysis. This requires us to formulate an ansatz on the form of the equilibrium and to develop detailed calculations, involving various measure changes, for the dynamics of the associated belief process. 

At the modelling level, one key difference between our work and the work reviewed above is the source of the asymmetry of information. In the literature, when a player holds additional information compared to her opponent, this usually comes in the form of a random variable which is independent of the Brownian motion driving the underlying dynamics. In \cite{Grun2013} such a random variable determines the players' payoffs, the situation is similar in \cite{GenGrun2019} although the dynamics are continuous-time Markov chains, in \cite{DEG2020} it determines the drift of the underlying SDE, in \cite{kwon2023exit} it determines each player's exit cost, etc. Instead, here the more informed player has perfect foresight about one characteristic of the underlying dynamics: she knows the value of $X_{\He}$. As a result, even though the two players observe the same realised trajectories of the process $X$ they use different models for the dynamics when assessing possible future trajectories.
The uninformed player uses Brownian dynamics whereas the two types of the informed player use two different Bessel processes, because they condition the Brownian motion to the value of the terminal state $X_{\He}$. Besides the purely mathematical interest, models of this kind may suggest applications in fields like economics and finance, where insider traders may use knowledge of future asset prices. Indeed, the idea that an informed agent knows a future characteristic of an observed process is central to the market-microstructure literature on insider trading initiated by Kyle \cite{kyle1985} and Back \cite{back1992}, see also \c{C}etin and Danilova \cite{cetin2018dynamic} for an overview.

The combination of learning mechanism and horizon is also different from the existing literature. In \cite{Grun2013} the uninformed player only learns from actions of the informed player because the underlying diffusive dynamics reveals nothing about the hidden variable. The game is set on finite horizon and the information is revealed at the terminal time, but at that point it is useless for the uninformed player. Instead, here the uninformed player learns from both the observation of the realised trajectory of $X$ and from the actions of the informed player. This is similar to the mechanism in \cite{DEG2020}. In the latter paper, it takes the uninformed player an infinite amount of time to exactly learn the drift of the SDE. Here instead the uninformed player learns the true state of the hidden variable at time $\He$ which is finite in expectation. Although the information at time $\He$ is again useless for the uninformed player, the observation of $X$ also reveals the arrival of the terminal time $\He$. Thus, the uninformed player learns simultaneously the random time horizon of the game and the terminal value of the Brownian motion. 

The main contributions of the paper are as follows. 
First, we introduce a novel class of Dynkin games with asymmetric information in which the informed player has foresight over an aspect of the future behaviour of the underlying process. Second, we extend the notion of a perfect Bayesian equilibrium to the setting of continuous-time Dynkin games with asymmetric information. Third, we show how results on general asymmetric Dynkin games can be applied to our setting to give conditions which delineate the equilibrium. Fourth, we develop an ansatz for the form of the equilibrium, in terms of the strategies of the uninformed maximiser and informed minimiser, under which the key determinants of the strategy are the underlying diffusion and its maximum process: the maximiser follows a pure threshold-based stopping rule; one type of the minimiser never stops and the other stops using a randomised rule based on the maximum process. Following the ansatz to its conclusion allows us to give concrete expressions for the equilibrium strategies in terms of the primitives of the problem. Fifth, we give a verification argument and show that we have indeed determined a perfect Bayesian equilibrium. 

The mechanism behind the randomised strategy of the informed minimiser is both novel and interesting. Even though there is a second-mover advantage she sometimes stops the game. She does this because in the scenarios where she does not stop she nonetheless changes the beliefs of the uninformed maximiser, and convinces the maximiser that it is optimal to stop, an ideal outcome for the minimiser.   

\subsection{Structure of the paper}
The paper is organised as follows. In Section \ref{sec:problem} we formulate the problem and state the main assumptions. We then provide easy sufficient conditions on the problem's primitives that imply our assumptions. In Section \ref{sec:main} we introduce the notion of PBE and we give an informal but intuitive statement about the structure of the equilibrium. In Section \ref{sec:construction} we formulate our ansatz on the form of the equilibrium and then we use it to transfer the martingale arguments from \cite{de2025nash} into a system of PDEs for the continuation values of the two types of the informed player and for the value function of the game. In Section \ref{sec:sol} we solve the system of PDEs to produce candidate values for the players and candidate optimal actions. In Section \ref{sec:verif} we show that our candidate leads to a PBE for the original game. Section \ref{sec:examples} is devoted to some explicit examples. Appendix~\ref{app0} gives some proofs of intermediate results and Appendix \ref{app} generalises our results to the case of a generic one-dimensional linear diffusion.

\section{Problem formulation}\label{sec:problem}
Let the tuple
$(\Omega,\cF,\F,(X_t)_{t\ge 0},(\P_x)_{x\in\R})$,
define a one-dimensional Brownian motion on a measurable space $(\Omega,\cF)$, equipped with a filtration $\F=(\cF_t)_{t\ge 0}$ and with a family of probability measures such that $\P_x(X_0=x)=1$. We assume that the space is complete with respect to the family $(\P_x)_{x\in\R}$ (cf.\ \cite[Sec.\ 1.5, p.27]{Blumenthal1968}). 
Although we present a model with Brownian motion stopped on first exit from $(0,1)$, standard techniques based on scaling and time change allow us to adapt our methodology to general one-dimensional, time-homogeneous diffusions on general intervals within the state space. This extension is outlined in Appendix~\ref{app} for the interested reader.

For $x\in(0,1)$ we let $X$ run until it hits either $0$ or $1$. Letting $\H_a=\inf\{t\ge 0:X_t=a\}$ for $a\in[0,1]$ we set $\He:=\H_{0}\wedge \H_{1}$ (here the subscript ``e''
stands for {\em exit}). We study a zero-sum game between two stoppers. Player 1 is a minimiser and chooses a random time $\tau \le \He$. Player 2
is a maximiser and chooses a random time $\sigma \le \He$. The exact class of admissible random times for the two players will be
specified later. At
time $\tau\wedge \sigma$ the game ends and the following payoff is paid by Player 1 to Player 2:
\begin{equation}
\label{eq:P2}
\cP(\sigma,\tau)\coloneqq
f\big(X_{\tau}\big)\ind_{\{\tau<\sigma\}}+g\big(X_{\sigma}\big)\ind_{\{\sigma<\tau\}}+h\big(X_{\tau}\big)
\ind_{\{\tau=\sigma\}},
\end{equation}
where precise assumptions on the functions $f$, $g$, $h$, are provided later, but typically $f,g,h$ are increasing functions with $g \leq h \leq f$. Both players observe the paths of the underlying
process, i.e., they have access to the filtration $\F^X\coloneqq(\cF^X_t)_{t\ge 0}$ generated by the process $X$. In addition to
that, Player 1 also knows the value of the random variable
\begin{equation}\label{eq:thetaX}
\theta^X :=  X_{\He} \, \, .
\end{equation}
That is, she knows whether the trajectory of the process exits the interval $(0,1)$ from the upper or lower
end-point. This induces asymmetry of information across the two players. Both players know the structure of the game and the
functions $f$, $g$, $h$; both players know that Player 1 is more informed. Formally, we say that Player 2's
filtration is $\F^2\coloneqq(\cF^2_t)_{t\ge 0}$ with
$\cF^2_t=\cF^X_t$, $t\ge 0$, whereas Player 1's filtration is $\F^1\coloneqq(\cF^1_t)_{t\ge 0}$ with
$\cF^1_t=\cF^X_t\vee\sigma(\theta^X)$, $t\ge 0$. It is worth noticing that although this is a so-called initial enlargement
of the
filtration, the additional information relates to the future position of the underlying process. In that sense, Player 1 foresees
the future and from time zero she has exact knowledge of the game's payoff in case neither player stops before $\He$. This type of
structure is rather different from that considered in, e.g., \cite{DEG2020} and \cite{Grun2013}. In \cite{Grun2013} the more informed player observes the value of an exogenous random variable $\theta$, independent of the commonly observed dynamics for $X$, which determines the form of the payoff functions $f_\theta(x)$, $g_\theta(x)$, $h_\theta(x)$ in the game. In \cite{DEG2020} the commonly observed dynamics for $X$ admits a drift which is known to the informed player, whereas the uninformed player only has a prior over the two possible values of $\theta\in\{0,1\}$, hence over two possible types of drift $\mu_0(x)$ and $\mu_1(x)$ (cf.\ Remark \ref{rem:erik}).

Zero-sum Dynkin games with asymmetry of information require the use of randomised stopping times for the players, as otherwise the
game is not well-posed in general (i.e., it does not admit a value). A definition of randomised stopping time is provided next.
Notice that since the game ends at time $\He$ at the latest, we are only interested in randomised stopping times bounded by
$\He$.
\begin{definition}\label{def:rst}
Given a filtration $\bH=(\cH_t)_{t\ge 0}$ with $\F^X\subseteq \bH$, we introduce a class of processes
\begin{equation*}
\begin{aligned}
\cAcirc(\bH):=\{\xi:&\,\text{$\xi$ is $\bH$-adapted with $t\mapsto \xi_t(\omega)$ c\`adl\`ag and non-decreasing,}\\
& \hspace{20mm} \text{$\xi_{0-}(\omega)=0$ and $\xi_\He(\omega) = \xi_\infty(\omega)= 1$ for all $\omega\in\Omega$}\}.
\end{aligned}
\end{equation*}
We say that $\rho$ is a \emph{$\bH$-randomised stopping time} bounded by $\He$, with {\em generating process} $\xi$ and {\em
randomisation device} $U$, if
\[
\rho=\inf\{t\ge 0 :\xi_t>U\},
\]
where $\xi\in\cAcirc(\bH)$ and $U\sim U(0,1)$ is a random variable independent of $\cH_\infty$. 

The class of randomised stopping times is denoted $\cT_R(\bH)$.
\end{definition}
For future reference we notice that $\rho\in\cT_R(\bH)$ can be expressed as
$\rho=\rho(U)$,
where $\rho(u)=\inf\{t\ge 0:\xi_t>u\}$ is an $\bH$-stopping time for any $u\in(0,1)$. It will be conveninent to extend the notation to $u=1$ by setting $\rho(1)=\inf\{t\ge 0:\xi_t = 1\}=\lim_{u \uparrow 1} \rho(u)$. The class of (pure) $\bH$-stopping times bounded by $\He$
is denoted $\cT(\bH)$ and for any $\tau\in\cT(\bH)$ we can define $\xi_t\coloneqq\ind_{\{\tau\le t\}}$ so that
\[
\inf\{t\ge 0:\xi_t>U\} =\tau.
\]
This proves the (well-known) inclusion $\cT(\bH)\subseteq\cT_R(\bH)$. Moreover, the following formulae are going to be used when
we compute expected payoffs associated to randomised stopping times: given a bounded measurable process $\Psi_t(\omega)$, a
randomised stopping time $\rho\in\cT_R(\bH)$ generated by $\xi\in\cA^\circ(\bH)$ and a stopping time $\eta\in\cT(\bH)$ we have
\begin{equation}\label{eq:useful}
\begin{aligned}
&\E\big[\ind_{\{\rho>\eta\}}\big|\cH_\infty\big]=1-\xi_\eta,\quad\E\big[\ind_{\{\rho=\eta\}}\big|\cH_\infty\big]=\xi_{\eta}-\xi_{\eta-}\\
&\text{and}\quad\E\big[\Psi_\rho\ind_{\{\rho<\eta\}}\big|\cH_\infty\big]=\E\Big[\int_{[0,\eta)}\Psi_t\ud
\xi_t\Big|\cH_\infty\Big].
\end{aligned}
\end{equation}
A proof of the above can be found, for example, in \cite[Lemmas 4.1 and 4.2]{de2022value}.

Now we can rigorously formulate the classes of admissible random times for the two players: Player 2 picks
$\sigma\in\cT_R(\F^2)=\cT_R(\F^X)$ and Player 1 picks
$\tau\in\cT_R(\F^1)=\cT_R(\F^X\vee\sigma(\theta^X))$. In the context of the
game, the randomisation device for Player 1 is independent from the randomisation device for Player 2. Later we will find that the equilibrium strategy for Player 2 is a pure strategy, but the equilibrium strategy for Player 1 is randomised.

As usual in zero-sum games we introduce an upper and a lower value of the game, denoted $\overline v$ and $\underline v$,
respectively, and defined as
\begin{equation}\label{eq:ulV}
\begin{aligned}
\underline v(x)\coloneqq\sup_{\sigma\in\cT_R(\F^2)}\inf_{\tau\in\cT_R(\F^1)}\E_x\big[\cP(\sigma,\tau)\big]\quad
\text{and}\quad
\overline v(x)\coloneqq\inf_{\tau\in\cT_R(\F^1)}\sup_{\sigma\in\cT_R(\F^2)}\E_x\big[\cP(\sigma,\tau)\big].
\end{aligned}
\end{equation}
The game admits a value when $\underline v=\overline v$ and it admits a saddle point when there is a pair
$(\sigma_*,\tau_*)\in\cT_R(\F^2)\times\cT_R(\F^1)$ such that for any other pair
$(\sigma,\tau)\in\cT_R(\F^2)\times\cT_R(\F^1)$ it
holds
\[
\E_x\big[\cP(\sigma,\tau_*)\big]\le \E_x\big[\cP(\sigma_*,\tau_*)\big]\le \E_x\big[\cP(\sigma_*,\tau)\big].
\]

In \cite{de2022value} and \cite{de2025nash} necessary and sufficient conditions for the existence of a saddle point are presented,
along with a methodology for the dynamic characterisation of optimal strategies. While those results indicate necessary and sufficient martingale properties of the equilibrium payoffs, they offer no additional insight on the actual structure of the equilibrium and on its construction. Here instead we construct an explicit equilibrium for the game, we provide computable formulae for equilibrium payoffs and we specify the dynamics of players' optimal randomised stopping times and beliefs.

Before delving into the solution method for our game, we want to set out our standing assumptions. However, in order to do that in
a meaningful way, we need to study properties of the expected payoff of the game in more detail. In particular, let us first
notice that for $\tau\in \cT_R(\F^1)$ we have the decomposition (\cite[Lemma D.2]{de2025nash})
\begin{equation}\label{eq:tau}
\tau=\ind_{\{\theta^X=0\}}\tau_0+\ind_{\{\theta^X=1\}}\tau_1,
\end{equation}
where $\tau_0,\tau_1\in\cT_R(\F^2)$. Then for any pair $(\sigma,\tau)\in\cT_R(\F^2)\times\cT_R(\F^1)$ we have
\begin{equation*}
\begin{aligned}
\E_x\big[\cP(\sigma,\tau)\big]&=\int_0^1\int_0^1\E_x\big[\ind_{\{\theta^X=0\}}\cP\big(\sigma(u),\tau_0(v)\big)+\ind_{\{\theta^X=1\}}\cP\big(\sigma(u),\tau_1(v)\big)\big]\ud u\, \ud
v\\
&=\int_0^1\int_0^1\E_x\big[\P_x\big(\theta^X=0\big|\cF^2_{\sigma(u)\wedge\tau_0(v)}\big)\cP\big(\sigma(u),\tau_0(v)\big)\big]\ud u\, \ud v\\
&\quad \quad \quad+\int_0^1\int_0^1\E_x\big[\P_x\big(\theta^X=1\big|\cF^2_{\sigma(u)\wedge\tau_1(v)}\big)\cP\big(\sigma(u),\tau_1(v)\big)\big]\ud u\, \ud v\\
&=\int_0^1\int_0^1\E_x\big[\big(1-X_{\sigma(u)\wedge\tau_0(v)}\big)\cP\big(\sigma(u),\tau_0(v)\big)\big]\ud u\,
\ud v\\
&\quad \quad \quad
+\int_0^1\int_0^1\E_x\big[X_{\sigma(u)\wedge\tau_1(v)}\cP\big(\sigma(u),\tau_1(v)\big)\big]\ud u\, \ud v,
\end{aligned}
\end{equation*}
where we used that $\P_x(\theta^X=1|\cF^2_t)=X_t$ in the third equality. Since
$t\mapsto X_{t\wedge \He}$ is a bounded non-negative martingale, we introduce the changes of measure
\begin{equation}\label{eq:P0P1}
\frac{\ud \P^0_x}{\ud\P_x}\Big|_{\cF^2_t}:=\frac{1-X_{t\wedge\He}}{1-x}\quad\text{and}\quad \frac{\ud
\P^1_x}{\ud\P_x}\Big|_{\cF^2_t}:=\frac{X_{t\wedge\He}}{x},\quad \P_x-a.s., \ x \in (0,1).
\end{equation}
By construction the measure $\P^0_x$ is the measure under which the process $X$ leaves the interval $(0,1)$ by
crossing $0$. Similarly, $\P^1_x$ is the measure under which $X$ leaves the interval $(0,1)$ by crossing $1$.
We can rewrite the expected payoff of the game as a convex combination of two expected payoffs under the measures $\P^0_x$ and
$\P^1_x$. That is
\begin{equation*}
\begin{aligned}
\E_x\big[\cP(\sigma,\tau)\big]&=\big(1-x\big)\int_0^1\int_0^1\E^0_x\big[\cP\big(\sigma(u),\tau_0(v)\big)\big]\ud u\,
\ud v
+x\int_0^1\int_0^1\E^1_x\big[\cP\big(\sigma(u),\tau_1(v)\big)\big]\ud u\, \ud v\\
&=\big(1-x\big)\E^0_x\big[\cP\big(\sigma,\tau_0\big)\big]+x\E^1_x\big[\cP\big(\sigma,\tau_1\big)\big].
\end{aligned}
\end{equation*}

By an application of It\^o's formula we have a more transparent form of the Radon-Nikodym derivatives that enables the use of
Girsanov's theorem:
\begin{equation*}
\begin{aligned}
\frac{\ud \P^0_x}{\ud\P_x}\Big|_{\cF^2_t}&=\exp\Big(\ln \frac{1-X_{t\wedge\He}}{1-x}\Big)=\exp\Big(-\int_0^{t\wedge\He}\frac{\ud X_u}{1-X_u}-\frac12\int_0^{t\wedge\He}\Big(\frac{1}{1-X_u}\Big)^2 \ud
u\Big),\\
\frac{\ud \P^1_x}{\ud\P_x}\Big|_{\cF^2_t}&=\exp\Big(\ln \frac{X_{t\wedge\He}}{x}\Big)=\exp\Big(\int_0^{t\wedge\He}\frac{\ud X_u}{X_u}-\frac12\int_0^{t\wedge\He}\Big(\frac{1}{X_u}\Big)^2 \ud
u\Big).
\end{aligned}
\end{equation*}
On the one hand, the process
$W^0_t=X_t-X_0+\int_0^t 1/(1-X_u)\ud u$,
is a standard Brownian motion under $\P^0_x$ for $t \le \He$ and the dynamics of $X$ under $\P^0_x$ reads as
\[
\ud X_t=-\frac{1}{1-X_t}\ud t+\ud W^0_t.
\]
On the other hand, the process
$W^1_t=X_t-X_0-\int_0^t1/X_u\ud u$, is a standard Brownian motion under $\P^1_x$ for $t \le \He$ and the dynamics of $X$ under $\P^1_x$ reads as
\[
\ud X_t=\frac{1}{X_t}\ud t+\ud W^1_t.
\]

Thanks to the above formulae it becomes transparent that the process $(X_{t \wedge \H_1})$ under $\P^1_x$ and the process $(1-X_{t\wedge \H_0})$ under $\P^0_x$ are $3$-dimensional Bessel processes corresponding to a Brownian motion conditioned to leave the interval $(0,1)$ from $1$.

\begin{remark}\label{rem:changeP0P1}
It is worth noticing that the measures $\P^0_x$ and $\P^1_x$ are indeed defined on $\cF$ as
\[
\frac{\ud \P^0_x}{\ud \P_x}\coloneqq\frac{\ind_{\{\theta^X=0\}}}{1-x}\quad\text{and}\quad \frac{\ud \P^1_x}{\ud \P_x}\coloneqq\frac{\ind_{\{\theta^X=1\}}}{x},
\]
whilst the formulae in \eqref{eq:P0P1} determine their projections onto $\cF^2_t$.
\end{remark}

Finally, in the statement of our technical assumptions we will use the infinitesimal generators
\[
\cL_X=\frac12\frac{\partial^2}{\partial x^2}, \quad
\cL^0_X = \frac{1}{2} \frac{\partial^2}{\partial x^2} - \frac{1}{1-x} \frac{\partial}{\partial x}\quad\text{and}\quad \cL^1_X = \frac{1}{2} \frac{\partial^2}{\partial x^2} + \frac{1}{x} \frac{\partial}{\partial x},
\]
of the process $X$ under the measures $\P_x$, $\P^0_x$ and $\P^1_x$, respectively.

\begin{assumption}\label{ass:fghscaled}
Functions $f$, $g$ and $h$ are such that:
\begin{itemize}
\item[(i)] on $[0,1]$ we have $f\ge h\ge g$ with $f$, $g$ continuous and $h$ measurable;
\item[(ii)] $h(0)=g(0)$, $h(1)=f(1)$;
\item[(iii)] on $(0,1]$ functions $f$ and $g$ are twice continuously differentiable with $f>g$, $f'>0$, $g'>0$;
\item[(iv)] $\cL_X^0 g(x)< 0$ and $\cL_X^1f(x)> 0$ for $x\in(0,1)$ and $g$ is concave on $(0,x_g]$ for some $x_g\in(0,1]$;
\item[(v)] for $0 < m \le 1$ the map 
\begin{equation}
\label{eq:phidef}
[0, m) \ni x\mapsto \oldh (x,m)\coloneqq x\frac{f(m)-g(x)}{m-x}-m f'(m),
\end{equation}
admits a root $c_*(m)\in(0,x_g \wedge m)$ and this root is unique in $(0, x_g \wedge m]$;
\item[(vi)] it holds that
\begin{equation}
\label{eq:newcondition}
1>c_*(1)\frac{g'(c_*(1))}{f'(1)} \geq \frac{g(x)-g(c_*(1))}{x-c_*(1)}\frac{1-x}{f(1) - g(x)}, \hspace{10mm} \forall x \in (c_*(1),1).
\end{equation}
\end{itemize}
\end{assumption}

\begin{remark}\label{rem:c*}
Notice that condition (v) is only about uniqueness of the root and its position relative to $x_g$. Indeed, for any $m \in (0,1]$, $\phi(0,m) = - mf'(m)<0$ and since $g(m)<f(m)$ it follows that $\lim_{x \uparrow m} \phi(x,m) = \infty$. Hence, by continuity of $\phi$ in its first argument and the intermediate value theorem, there exists a root $c=c(m)$ in $(0,m)$ such that $\phi(c(m),m)=0$. Moreover, any root of $\phi(\cdot,m)$ must be strictly contained in $(0,m)$.
\end{remark}

Some comments about the nature and implications of our assumptions are in order. We emphasise that all conditions stated above can be verified directly from the data $f(x)$ and $g(x)$ with no reference to the game. The conditions about $c_*$ are implied by the simpler ones that are illustrated in Lemma \ref{lem:jan} and \ref{lem:newsuff}. Examples with explicit solutions will be provided in Section \ref{sec:examples}.

The stopping decisions of the players are influenced by several factors. First,
the condition $f\ge h\ge g$ says that players have an incentive to let their opponent end the game (so-called second-mover
advantage). In particular, stopping simultaneously for $t<\He$ is never strictly beneficial for either player.
Second, since $f',g'>0$ at the broad level, the minimiser prefers the process to go down whereas the maximiser wants it to go up. Third, the stopping rules of the players are influenced by the super/sub-martingale properties of the payoff processes $t\mapsto (f(X_t), g(X_t))$.
The assumption that $\cL_X^0 g < 0$ means that if the uninformed maximiser knew that $\theta^X = 0$ and that the
minimiser would not stop before $\He$, they would choose to stop immediately. Similarly, the assumption that $\cL_X^1 f > 0$
means that if $\theta^X=1$ and the maximiser does not stop before $\He$, then it is optimal for the informed
minimiser to stop immediately.
Taken together, these assumptions introduce the jeopardy whereby players face a trade-off between stopping first (and foregoing the second mover advantage) and waiting (and facing adverse expected payoffs), and this is what makes the game non-trivial.
Indeed, we can be a bit more precise. When $\theta^X=0$ the informed minimiser (Player 1) has no incentive to stop before $\He$. The second mover advantage guarantees $\cP(\sigma,\tau)\ge g(X_{\sigma \wedge \tau})$ and conditionally on $\theta^X=0$ the process $t\mapsto g(X_t)$ is a supermartingale. Therefore, by choosing $\tau=\He$ the minimiser attains $g(X_\sigma)$, which is decreasing on average over time. Knowing this, when the uninformed player believes that $\theta^X=0$ they should stop as early as possible.
The situation when $\theta^X=1$ is more nuanced. If the maximiser does not stop before $\He$, then as argued above, the minimiser should stop immediately, thus forgoing the second-mover advantage. However, the minimiser can  benefit from the second-mover advantage if they can `fool' the maximiser into stopping. The minimiser can do this by following a randomised strategy and sometimes stopping. Then, in the scenarios where the minimiser does not stop, the maximiser's belief that $\theta^X=0$ is strengthened. This incentivises the maximiser to stop, bringing second-mover benefits to the minimiser. The problem facing the minimiser who knows that $\theta^X=1$ is to design a stopping rule from which the benefits from the second mover advantage more than compensates them for the fact that they may end up stopping the game at a higher value. This is what makes the game novel and interesting.

We now give a series of lemmas concerning sufficient conditions for (v) and (vi) (Lemma \ref{lem:jan} and \ref{lem:newsuff}) and regularity of the function $c_*(m)$ (Lemma \ref{lem:c*}). To streamline the exposition, we move their simple proofs to Appendix \ref{app0}.

\begin{lemma}\label{lem:jan}
Assume that $\cL^1_X g\neq 0$ on $(0, x_g]$ and that when $x_g<1$ we have $\oldh(x_g, m) > 0$ for all $m \in (x_g,1]$.
Then condition (v) in Assumption \ref{ass:fghscaled} holds.
\end{lemma}

\begin{lemma}\label{lem:newsuff}
Suppose that the first inequality of \eqref{eq:newcondition} in (vi) of Assumption \ref{ass:fghscaled} holds, i.e., $f'(1)> c_*(1) g'(c_*(1))$. Then the second inequality in \eqref{eq:newcondition} holds if
\begin{equation}
\label{eq:newsufficient}
\frac{g(x)-g(c_*(1))}{x-c_*(1)}  \leq g'(c_*(1)), \quad \forall x \in (c_*(1),1).
\end{equation}
Clearly, \eqref{eq:newsufficient} is true if $g$ is concave.
\end{lemma}

\begin{remark}
Taken together, conditions (i)--(iv) in Assumption \ref{ass:fghscaled} and the conditions in Lemma \ref{lem:jan} and \ref{lem:newsuff} offer easy to verify sufficient conditions for the equilibrium that we are going to construct in the next sections. For example,  $\cL^1_X g>0$, 
and $\cL_X g\le 0$ reflect a situation in which the uninformed player considers the first-mover's payoff process $g(X_t)$ as a $\P_x$-supermartingale (hence, late stopping is unfavourable), but at the same time knows that conditionally on $\theta^X=1$ the first-mover's payoff is a submartingale (hence, late stopping is optimal). This creates a tension between stopping as early as possible, in order to avoid the negative effect of the supermartingale dynamics, and continuing in order to learn the true value of $\theta^X$ -- if Player 2 knew that $\theta^X$ is equal to one then Player 2 should never stop.

\end{remark}

\begin{lemma}\label{lem:c*}
The function $m\mapsto c_*(m)$ is continuous on $(0,1]$ and enjoys a continuous extension to $[0,1]$ with $c_*(0) = 0$. Moreover, $c_*(1) < 1$ and $c_*$ is continuously differentiable with $c'_* > 0$ on the subset $\{m\in(0,1): m^2f'(m)>c_*(m)g'(c_*(m))\}$.
\end{lemma}

Monotonicity and continuous differentiability of $c_*$ require $m^2 f'(m) >  c_*(m)^2 g'\big(c_*(m)\big)$. From the first inequality in \eqref{eq:newcondition} of Assumption
\ref{ass:fghscaled}-(vi) and continuity of $c_*$ we know that $m^2 f'(m) > c_*(m)^2 g'\big(c_*(m)\big)$ for $m$ in a left-neighbourhood of
$1$ (notice that $c_*>c^2_*$ and $f',g'>0$). Then we can set
\begin{equation}\label{eqn:m0}
m_0 = m_0(c_*) \coloneqq\sup\{m\in[0,1):m^2 f'(m) \le c_*(m)^2 g'\big(c_*(m)\big)\},\quad\text{with $\sup\varnothing=0$},
\end{equation}
so that $m_0<1$ and $m^2 f'(m) > c_*(m)^2 g'\big(c_*(m)\big)$ for $m\in(m_0,1]$.

\section{A perfect Bayesian equilibrium}\label{sec:main}

In this section we introduce a notion of Bayesian perfection for our framework and sketch the main result of the paper -- a perfect Bayesian equilibrium with asymmetric information. The construction and the proof are developed in several steps in the subsequent sections.

The underlying dynamics of $(X_t)$ and the payoff $\cP$ have a Markovian structure. This would also be a desirable property of the equilibrium.
In full and symmetric information stopping games with a payoff based on an underlying strong Markov process, a Markov perfect equilibrium as defined in \cite{christensen2024general, decamps2022mixed} consists of a pair of Markovian randomised stopping times such that each is the best response to the other for any initial condition. A randomised stopping time with generator $(\xi_t)$ is said to be Markovian if $(1-\xi_t)$ is a multiplicative functional of the underlying Markov process (in \cite{decamps2022mixed}) or, equivalently, if $\log(1-\xi_t)$ is an additive functional of that process (in \cite{christensen2024general}).

In contrast, under asymmetric information the underlying Markovian state must be augmented to include the uninformed player's belief, thereby modifying both the definition of Markovian strategies and the equilibrium concept. To the best of our knowledge, the literature does not offer a definition of a perfect Bayesian equilibrium that applies to the game studied in this paper.
Our aim, therefore, is to give a definition which is consistent with the general principles that underpin existing definitions in related contexts. In particular, we aim to address the following points:
\begin{enumerate}[label=(\roman*)]
\item we identify a belief process $(p_t)$ representing the conditional probability of the unobserved state given the information available to the uninformed player,
\item we find a set of admissible initial points for the belief-augmented dynamics $(X_t,p_t)$,
\item we define a Markovian structure of strategies on the belief-augmented state space,
\item we state the optimality of strategies for any admissible initial point.
\end{enumerate}
The introduction of the admissible set in (ii) arises because not all pairs $(x, p)$ can be attained on the equilibrium path. This follows from an intrinsic feature of stopping games: any observable deviation by any of the players terminates the game. Hence, irrespective of the strategy pursued by Player 1, the belief process of Player 2 evolves as if Player 1 followed the equilibrium strategy, at least until the moment that the game is stopped.

Our equilibrium notion will be closely related to the concept of perfect Bayesian equilibrium (PBE) in dynamic games with incomplete information, see \cite{Osborne1994}. A formal definition of PBE for a specific non-zero sum stopping game that explicitly incorporates belief systems and Markovianity of strategies is given in \cite[Sec.~4.7]{kwon2023exit}. There, the players' strategies are pure stopping times, which greatly simplifies the theory; this is not the case here.

In the following paragraphs we provide details for all points (i)--(iv) of the perfect Bayesian equilibrium. To this end, fix a candidate saddle point $(\sigma, \tau) \in \cT_R(\F^2)\times\cT_R(\F^1)$ and recall that $\tau$ has the decomposition \eqref{eq:tau}. Denote by $(\zeta_t)$ the generating process for $\sigma$ and by $(\xi^i_t)$ the generating process for $\tau_i$, $i=0,1$.

To simplify presentation, we will restrict to candidate equilibrium strategies of the minimiser that satisfy
\begin{equation}\label{eqn:standing}
\xi^0_t \vee \xi^1_t < 1\quad \text{for $t < \He$.}
\end{equation}
This condition will be formalised in Definition \ref{def:rst-markov}.
\medskip

\subsection{The belief process}\label{sec:belief}
The asymmetry of information stems from the fact that Player 2 does not observe $\theta^X$. She, therefore, keeps track of the distribution of $\theta^X$ given the observed dynamics of $(X_t)$ and the fact that Player 1 has not terminated the game yet. As $\theta^X$ takes two values, it is sufficient to track the probability of $\theta^X = 1$. The belief-augmented state space takes the form $(0, 1)^2$ with its points denoted by $(x,p)$. For any $(x,p)$, the dynamics of the process $(X_t)$ is such that at time zero $X_0 = x$ and the probability that $X_{\He} = 1$ is equal to $p$. This determines a conditional measure $\hP^p_x$ which we formally define in the following lemma, whose simple proof is a direct consequence of Remark \ref{rem:changeP0P1} and is omitted.

\begin{lemma}\label{lem:Y-Pyp}
Under the measure $\hP^p_{x} = p \P^1_x + (1-p) \P^0_x$, we have $\hP^p_{x}(\theta^X = 1) = p$ and
\[
\frac{\ud \hP^p_{x}}{\ud \P_x} = \ind_{\{\theta^X = 1\}} \frac{p}{x} + \ind_{\{\theta^X = 0\}} \frac{1-p}{1-x}.
\]
\end{lemma}
In the notation $\hP^p_x$, the subscript $x$ has a standard meaning of indicating the initial value of the process $(X_t)$ while the superscript $p$ determines the weights in the mixture of measures $\P^1$ and $\P^0$.
The original problem formulation in \eqref{eq:ulV} uses the measure $\hP_{x}^{x}=\P_x$.

We can now formally state the belief process as
\[
p_t = \hP^p_x(\theta^X = 1\, |\, \cF^2_t, \tau > t) \coloneqq \frac{\hP^p_x(\theta^X = 1, \tau > t  \, |\, \cF^2_t)}{\hP^p_x(\tau > t\, |\, \cF^2_t)},
\]
where the final expression defines the conditional probability given the $\sigma$-algebra $\cF^2_t$ and the event $\{\tau > t\}$. Similar conditional probabilities will also be used later in the paper and the reader should interpret them in an analogous way. Notice that when  \eqref{eqn:standing} holds and $t < \He$, the event $\{\hP^p_x(\tau > t |\cF^2_t) = 0\}$ has $\hP^p_x$-measure zero, so $p_t$ is well-defined and $p_t \in (0,1)$. It will prove convenient to consider also
$\hat p_t = \hP^p_x(\theta^X = 1 | \cF^2_t)$. Thanks to the explicit form of the Radon-Nikodym derivative of $\hP^p_x$ with respect to $\P_x$ (Lemma \ref{lem:Y-Pyp}) we further have
\begin{equation}\label{eq:hatpX}
\hat p_t =\hP^p_x(\theta^X = 1|\cF^2_t) = \frac{p \frac{X_t}{x}}{p \frac{X_t}{x} + (1-p) \frac{1-X_t}{1-x}} = \frac{pX_t (1-x)}{p X_t (1-x) + (1-p)(1-X_t)x}.
\end{equation}

Using $(\hat p_t)$, the belief process $(p_t)$ can be written as
\[
p_t = \frac{\hP^p_x(\theta^X = 1, \tau > t  \, |\, \cF^2_t)}{\hP^p_x(\tau > t\, |\, \cF^2_t)} = \frac{\hP^p_x(\tau > t  \, |\, \cF^2_t,\theta^X = 1)\hP^p_x(\theta^X = 1 \, |\, \cF^2_t)}{\hP^p_x(\tau > t\, |\, \cF^2_t)}
=
\frac{(1-\xi^1_t) \hat p_t}{\hP^p_x(\tau > t\, |\, \cF^2_t)},
\]
where in the final equality we used
\[
\hP^p_x(\tau > t  \, |\, \cF^2_t,\theta^X = 1)=\hP^p_x(\tau_1 > t  \, |\, \cF^2_t)=1-\xi^1_t.
\]
We apply a similar idea to the denominator:
\begin{equation}
\begin{aligned}
\hP^p_x(\tau > t\, |\, \cF^2_t)
&=
\hP^p_x(\tau > t  \, |\, \cF^2_t,\theta^X = 1)\hP^p_x(\theta^X = 1 \, |\, \cF^2_t)\\
&\quad+
\hP^p_x(\tau > t  \, |\, \cF^2_t,\theta^X = 0)\hP^p_x(\theta^X = 0 \, |\, \cF^2_t)
=
(1-\xi^1_t) \hat p_t + (1-\xi^0_t) (1-\hat p_t)
\end{aligned}
\end{equation}
yielding
\begin{equation}\label{eqn:pt}
p_t = \frac{(1-\xi^1_t) \hat p_t}{(1-\xi^1_t) \hat p_t + (1-\xi^0_t) (1-\hat p_t)}.
\end{equation}
Hence, given $(\xi^0_t, \xi^1_t)$, there is a one-to-one correspondence between $p_t$ and $\hat p_t$ given by
\[
\hat p_t =
\frac{p_t}{p_t+\frac{1-\xi^1_t}{1-\xi^0_t}(1-p_t)}.
\]
\begin{remark}\label{rem:erik}
The dynamics of $(X_t, \hat p_t)$ can be formally derived using filtering theory. From the point of view of Player 2, the diffusion $X$ has a partially observable drift with two possible forms:
\[
\mu_0(x)=-\frac{1}{1-x}\quad\text{and}\quad\mu_1(x)=\frac{1}{x}.
\]
Following standard filtering theory, for $t<\He$ we have
\begin{equation}\label{eqn:Y-filter}
\ud X_t = \big( (1-\hat p_t) \mu_0(X_t) + \hat p_t \mu_1(X_t)\big) \ud t + \ud B_t, \quad X_0 = x,
\end{equation}
and
\begin{equation}\label{eqn:p-filter}
\ud \hat p_t = \chi(X_t) \hat p_t (1-\hat p_t) \ud B_t, \quad \hat p_0 = p,
\end{equation}
where $\chi (x) = \mu_1(x) - \mu_0(x)$. The process $(B_t)$ is the innovation process
\[
B_t = X_t-X_0 - \int_0^t \big[(1-\hat p_s) \mu_0(X_s) + \hat p_s \mu_1(X_s)\big] \ud s,
\]
which is a Brownian motion with respect to filtration $\F^2$ (generated by $(X_t)$) and the measure $\hP^p_x$ for $t<\He$. We remark that the dynamics can be extended up to $\He$ and the pair $(X_t, \hat p_t)$ forms a strong Markov process for $t\le\He$. However, the dynamics are not well-posed after $\He$. In our work this is not an issue as the the game ends at time $\He$, when both players become equally informed.

In \cite{DEG2020}, an informed player knows the drift of the process $X$ and an uninformed player does not know the drift but acquires information upon observing $X$ and the actions of the informed player. Unlike in our setting, the uninformed player in \cite{DEG2020} cannot determine the realised drift with certainty in finite time (unless the informed player ends the game) because the posterior process $(\hat p_t)$ is well-defined and takes values in $(0,1)$ for all $t\in[0,\infty)$. Moreover, the solution approach of \cite{DEG2020} is not promising for our problem because the state dynamics $(X_t,\hat p_t)$ for our problem are two-dimensional and degenerate-diffusive. The solution of the system of resulting free boundary PDEs in closed form seems impossible. Instead we follow a different route which is detailed in the next sections. 
\end{remark}

\subsection{The set of admissible initial points and Markovian strategies}

At the outset of the game, the uninformed player has no information about the distribution of $\theta^X$ apart from that arising from the dynamics of $(X_t)$ itself. Hence, given $x \in (0,1)$, we have $\hat p_0 = \P_x(\theta^X = 1) = x$. Consider the values that the process $(X_t, p_t)$ can attain given the initial point $(X_0, \hat p_0) = (x, x)$, $x \in (0,1)$, and given that $(p_t)$ is calculated using the generating process of the equilibrium randomised stopping time $\tau$ chosen by Player 1. These values are admissible initial points of the game. The set of those admissible initial points will be denoted by $A^X=A^{X,\tau}$ and called an admissible set. For convenience, we will drop the minimality and allow additional points as long as they satisfy the conditions of the following definition.
\begin{definition}\label{def:admissible}
A Borel set $A^X \subset (0,1)^2$ is an \emph{admissible set} for the process $(X_t,p_t)$ if
\begin{enumerate}
 \item $(x,x) \in A^X$, for all $x \in (0,1)$,
 \item $\hP_{x'}^{p'}\big((X_t, p_t) \in A^X \text{ for all $t \in [0, \He)$}\big) = 1$ for every $(x', p') \in A^X$.
\end{enumerate}
\end{definition}
Although in general we do not address existence of such a set, in our equilibrium we can determine $A^X$ explicitly. We also notice that the definition could be sharpened by adding the requirement that $A^X$ be the smallest set satisfying (1) and (2). This is not necessary for our use.

A necessary feature of a perfect Markov equilibrium is that the strategies enjoy an appropriate Markovian property. In \cite[Def.~2.4]{decamps2022mixed}, in the case of full and symmetric information (when both players use filtration $\F^2$) a randomised stopping time with the generating process $(\xi_t)$ is called Markovian with respect to a Markov process $(X_t)$ if $\bar \xi_t \coloneqq 1 - \xi_t$ is a multiplicative functional of the Markov process $(X_t)$. That is, for any $\F^2$-stopping time $\eta$, any $s \ge 0$ and any $x$ we have
\[
\bar \xi_{\eta + s} = \bar \xi_\eta (\bar \xi_s \circ \theta_\eta),\quad \P_x\as,
\]
where $\theta_\eta$ is the standard shift operator for the Markov process $(X_t)$. In the Bayesian case, the generating process depends on the belief process which itself depends on the strategy of Player 1, so the above definition is appropriately adapted. In preparation for that, we rewrite the above formula as $\bar \xi_{\eta + s} / \bar \xi_\eta = \bar \xi_s \circ \theta_\eta$ on $\{\bar \xi_\eta > 0\}$. It then implies a weaker condition (but sufficient for our needs):
\begin{equation}\label{eqn:mult-func-full}
\E_x[\varphi(\bar \xi_{\eta + s} / \bar \xi_\eta) | \cF^2_\eta] = \E_{X_\eta}[\varphi(\bar \xi_s)]  \quad \text{on $\{\bar \xi_\eta > 0\}$},
\end{equation}
for any (bounded) Borel function $\varphi$. Notice that $\bar\xi_{\eta + t} = 0$ for $t \ge 0$ on the event $\{\bar\xi_\eta=0\}$ is trivially satisfied since $(\bar\xi_t)$ is non-increasing and takes values in $[0,1]$. Then discarding this event in \eqref{eqn:mult-func-full} comes at no detriment.

Players' strategies could and should depend on the initial belief $p$ of Player 2 which determines the measure $\hP^p_x$ under which Player 2 assesses the game. That is, for any $p \in (0,1)$, there is a pair of randomised stopping times $(\sigma^p,\tau^p)\in\cT_R(\F^2)\times\cT_R(\F^1)$, where $\tau^p$ has the decomposition \eqref{eq:tau} of the form $\ind_{\{\theta = 1\}} \tau^{1,p} + \ind_{\{\theta = 0\}} \tau^{0,p}$ for $\tau^{i,p} \in \cT_R(\F^2)$, $i = 0,1$. The corresponding generating processes are $(\zeta^p_t)$, $(\xi^{1,p}_t)$, $(\xi^{0,p}_t)$, respectively. For the sake of compactness we will write $\sigma \coloneqq (\sigma^p)_{p \in (0,1)}$ with the family of generating processes $\zeta \coloneqq (\zeta^p)_{p \in (0,1)}$ and $\tau^i \coloneqq (\tau^{i,p})_{p \in (0,1)}$ with $\xi^i \coloneqq (\xi^{i,p})_{p \in (0,1)}$. We will refer to the pair $(\sigma, \tau)$ as \emph{belief-dependent strategy profile}. We identify a subclass of generating processes that satisfy a measurability property with respect to $p$: we will call $(\nu^p_t)_{t \ge 0}$ \emph{a belief-dependent generating process} if $\nu^p_t = \hat \nu(t, X, p)$, $t \ge 0$, and $\nu^p\in\cA^\circ(\F^2)$ for each $p\in(0,1)$, where $\hat \nu:[0, \infty) \times C([0,\infty); \R) \times (0,1) \to [0,1]$ is a measurable map.

Under $\hP^p_x$ we often use $(\sigma^p, \tau^p)$, in which case we write $(\sigma, \tau)$. We stress the distinction between $p$ in $\hP^p_x$ which determines the dynamics of the underlying process for Player 2 and $p$ in the pair $(\sigma^p, \tau^p)$, where it selects the strategy dependent on the belief of Player 2. Although the belief $p$ of Player 2 aligns with the mixing in measure $\hP^p_x$ in the game, mathematical analysis demands clear distinction between those two concepts. 

\begin{definition}\label{def:rst-markov}
A belief-dependent strategy profile $(\sigma, \tau)$ with families of belief-dependent generating processes $\zeta$, $\xi^1$, $\xi^0$, is called belief-state-\emph{Markovian} if for any bounded Borel $\varphi:\R\to\R$, any $(x,p)\in A^X$, any
$\eta\in\cT(\F^2)$ with $\eta<\He$, and any $s\ge 0$, we have, $\hP^p_{x}$\as 
\begin{equation}\label{eqn:rst-markov-xi}
\hE^p_{x}\bigg[\varphi\bigg(\frac{1-\xi^{i,p}_{\eta+s}}{1-\xi^{i,p}_{\eta}}\bigg)\,\bigg|\,
\cF^2_\eta,\ \tau^p>\eta\bigg]=\hE^{p_\eta}_{X_\eta}\big[\varphi(1-\xi^{i,p_\eta}_s)\big], \qquad i=0,1,
\end{equation}
and
\begin{equation}\label{eqn:rst-markov-zeta}
\hE^p_{x}\bigg[\varphi\bigg(\frac{1-\zeta^p_{\eta+s}}{1-\zeta^p_{\eta}}\bigg)\,\bigg|\,
\cF^2_\eta,\ \tau^p>\eta\bigg]=\hE^{p_\eta}_{X_\eta}\big[\varphi(1-\zeta^{p_\eta}_s)\big], \qquad \text{on
$\{\zeta^p_\eta < 1 \}$.}
\end{equation}
For $\eta=0$ the belief $p_0$ is the one after a possible jump of $(\xi^0,\xi^1)$ at time $0$.

Furthermore, for any $(x,p) \in A^X$, we have 
\begin{equation}\label{eqn:unfinished}
\hP^p_x(\xi^{1,p}_t \vee \xi^{0,p}_t < 1 \text{ for all }t < \He) = 1.
\end{equation}
\end{definition}

Notice that the last condition of the above definition is the formalisation of \eqref{eqn:standing} which greatly simplifies the presentation in this section.
The equalities \eqref{eqn:rst-markov-xi} and \eqref{eqn:rst-markov-zeta} use that the right-hand side is well-defined, i.e., the mapping $(x,p) \mapsto \hE_{x}^{p}[\varphi(1-\nu^p_s)]$ is appropriately measurable for $\nu = \xi^0, \xi^1, \zeta$. This, however, can be established directly. Indeed, for any $t \ge 0$ and any $A \in \cF^2_t$, the mapping $(0,1) \ni x \mapsto (\P^0_x(A),\P^1_x(A))$ is Borel measurable by \eqref{eq:P0P1}. Consequently the mapping $(0, 1)^2 \ni (x,p) \mapsto \hP^p_{x}(A)$ is Borel measurable thanks to its explicit definition in Lemma \ref{lem:Y-Pyp}. It remains to combine this with the joint measurability of belief-dependent generating processes.

The conditioning on $\tau^p > \eta$ on the left-hand side of \eqref{eqn:rst-markov-xi} and \eqref{eqn:rst-markov-zeta} stems from the fact that the process $(X_t, p_t)$ is not Markovian while the process $(X_t, \hat p_t)$ is. The belief process $(p_t)$ diverges from $(\hat p_t)$ when the informed player performs stopping.
Hence, the adjustment for the survival past time $\eta$ is required to bring back the Markov-like identity.

Notice that Definition \ref{def:rst-markov} reduces to the identity \eqref{eqn:mult-func-full} of a Markovian stopping time for games with full and symmetric information. Indeed, the process $(p_t)$ is then constant and irrelevant and there is one $\F^2$-adapted generating process for $\tau$, say $(\xi^1_t)$. 

\subsection{Optimality of strategies and an informal statement of the main theorem}
The optimality condition for a Markovian saddle point $(\sigma, \tau)$ is easily adapted from the classical theory of perfect Markov equilibrium for stopping games.
\begin{definition}\label{def:perfectBayesian}
A pair of belief-state-Markovian strategies $(\sigma, \tau)$ is a perfect Bayesian equilibrium (PBE) if for any $(x, p) \in A^X$ we have
\[
\hE^p_x[\cP(\sigma', \tau)] \le \hE^p_x[\cP(\sigma, \tau)] \le \hE^p_x[\cP(\sigma, \tau')]
\]
for any $(\sigma', \tau') \in \cT_R(\F^2) \times \cT_R(\F^1)$.
\end{definition}
Recall the discussion at the beginning of this section: the set $A^X$ and the dynamics of $(p_t)$ depend on $\tau$. In the right-most expectation the deviation from $\tau$ to $\tau'$ does not change the belief of the uninformed player as any observable deviation from $\tau$ terminates the game.

Finally, before closing the section we outline the main theorem. A complete formal statement requires the use of several functions which will be derived in the later sections. It would therefore be extremely technical at this stage to introduce all the ingredients for such a statement. Instead, we give an informal statement that summarises the structure of a PBE in our game and we postpone the fully rigorous statement to Theorem \ref{thm:mainscaled} at the end of Section \ref{sec:verif}. For the informal statement we need to complement the state process with its running maximum $(M_t)$ started from an initial value $M_0$ that must be determined in equilibrium: $M_t=M_0\vee\sup_{s\le t}X_s$; we also need a non-increasing, absolutely continuous function $\Gamma_*:(0,1]\to (0,1]$, whose derivative we denote $\Gamma'_*(x)$, also to be determined in equilibrium.

{\bf Outline of our main theorem (Theorem \ref{thm:mainscaled}):}
There is a PBE in our game in which players stop according to the following rules:
\begin{itemize}
\item Player 1, knowing that $\theta^X=0$, never stops.
\item Player 1, knowing that $\theta^X= 1$, has two actions, depending on the initial belief $p$ of Player 2:
\begin{itemize}
\item If Player 2's belief {\em is not} too high, then Player 1 stops at times of increase of the process $(M_t)$ with an intensity determined by the generating process:
\begin{equation}\label{eq:heur}
-\int_0^t \frac{1}{\Gamma_*(M_0)} \Gamma'_*(M_u) \ud M_u.
\end{equation}
Furthermore, Player 1 does not act until the process $(M_t)$ reaches some level $m^1_* \ge 0$, i.e., $\Gamma_*$ is constant for $m \leq m^1_*$.
\item If Player 2's belief {\em is} too high, then Player 1 stops at time zero with some probability so as to shift Player 2's belief to a suitably lower level and then she stops at times of increase of the process $(M_t)$ according to \eqref{eq:heur}.
\end{itemize}
\item As a result of the actions of Player 1, the belief of Player 2 after time zero evolves according to $p_t=\pi_*(X_t,M_t)$ for some function $\pi_*$.
\item Player 2 stops according to
\[
\sigma_*=\inf\big\{t\ge 0: X_t\le c_*(M_t \vee m^1_*)\big\}\wedge \He,
\]
where $c_*$ is the root defined in Assumption \ref{ass:fghscaled}-(v). 
\end{itemize}

\begin{remark}\label{rem:3.6}
We illustrate with an example how the asymmetry of information introduces genuinely new structure in the game. Suppose both players are equally uninformed and observe $\F^2$. Consider a game played with payoff functions $f(x)=h(x)=x$ and $g(x)=\alpha x$ for some $\alpha\in(0,1)$. Notice that in this game, for any pair $(\sigma,\tau)\in\cT_R(\F^2)\times\cT_R(\F^2)$ such that $\sigma,\tau\le \He$,
\[
\E_x\big[\cP(\sigma,\tau)\big]\le \E_x\big[X_{\tau\wedge\sigma}\big]=x.
\]
Player 2 (the maximiser) can attain the value on the right-hand side of the inequality by choosing $\sigma=\He$, irrespective of Player 1's choice of stopping rule. So for any $\tau\in\cT_R(\F^2)$ the payoff reads $\E_x[\cP(\He,\tau)]=x$ and Player 1 is indifferent across all choices of stopping rule. As a result, any pair $(\He,\tau)$ is an equilibrium and the value of the game is equal to $x$.

In Section \ref{subsec:linear payoffs} we study the same problem with asymmetry of information and we show that the above analysis breaks down. The more informed minimiser can attain a lower value thanks to the additional information available to her.
\end{remark}


\section{Construction of a candidate equilibrium}\label{sec:construction}

From now on we will focus on the study of our game.
In order to formalise the equilibrium outlined in the final paragraph of the previous section, we extend the state space to consider the Markovian pair $(X_t,M_t)_{t\ge 0}$ with $M_t=M_0 \vee\sup_{0\le s\le
t}X_s$ with deterministic $M_0$. Because we extended the state space from $X$ to $(X,M)$ we also introduce the family of probability measures $(\P_{x,m})_{0\le x\le m\le 1}$ such that $\P_{x,m}(X_0=x,M_0=m)=1$ and $\P_{x,x}(A)=\P_x(A)$ for any $A\in\cF$. Likewise we define the conditional measures $\P^j_{x,m}(\cdot)=\P_{x,m}(\cdot|\theta^X=j)$, $j=0,1$. Accordingly, we are going to take expectations $\E_{x,m}[\,\cdot\,]$, $\E^0_{x,m}[\,\cdot\,]$ and $\E^1_{x,m}[\,\cdot\,]$.
For simplicity we abbreviate $(\theta^X,\cL_X,\cL_X^{0,1})$ to $(\theta,\cL,\cL^{0,1})$.

It is clear that both players observe the pair $(X,M)$ which is $\F^2$-adapted. The candidate stopping rule for Player 2
is of the form
\begin{equation}\label{eq:sigmac}
\sigma_d =\inf\{t\ge 0\colon X_t\le d(M_t)\}\wedge \He,
\end{equation}
for some measurable function $d:[0,1]\to [0,1]$ which we expect to be at least increasing and continuous and satisfying $0 \leq d(m) \leq (m \wedge x_g]$ (with $d(m)<m$ on an interval of the form $(m_d,1]$ for some $0 \leq m_d <1$). Our ansatz is that the informed player who knows $\theta=0$ never stops and the informed player who knows $\theta=1$ only stops when the process reaches a new maximum. The candidate stopping rule for Player 1 is then of the form $\bar\tau_\Gamma=\ind_{\{\theta=0\}}\He+\ind_{\{\theta=1\}}\tau_\Gamma$,
where
\begin{equation}
\label{eqn:tauGamma}
\tau_\Gamma\coloneqq\inf\{t\ge 0\colon 1-\Gamma(M_t)>U\}\wedge \He,
\end{equation}
for some non-increasing continuous function $\Gamma:(0,1]\to (0,1]$.
For $u \in (0,1)$ define the $\F^2$-stopping time $\tau_\Gamma(u) \coloneqq\inf\{t\ge 0\colon 1-\Gamma(M_t)>u \}\wedge \He$ and $\tau_\Gamma(1) \coloneqq
\lim_{u \uparrow 1} \tau_\Gamma(u)$. Then, $u\mapsto\tau_\Gamma(u)(\omega)$ is measurable for every $\omega\in\Omega$ and $\tau_\Gamma=\tau_\Gamma(U)$. For the
calculations that will follow, we find it convenient to represent the increasing process that generates $\tau_\Gamma$ as $\xi_t=1-\ind_{\{M_t < 1\}}\Gamma(M_t)$, where the indicator function of the event $\{M_t < 1\}$ executes the stopping at $\He$. The pair of functions $(d(m),\Gamma(m))$ is to be determined in equilibrium and it is one of the main
objects of interest in our study. For the moment we ignore questions around Bayes perfection of the candidate equilibrium. We will come back to it in Section \ref{sec:verif}, where we will show that there is a bijection between the processes $(X_t,M_t)$ and $(X_t,p_t)$, but the former provides a more convenient parametrisation for this part of the study.
Given a pair $(\sigma_d,\bar \tau_\Gamma)$ we want to derive expressions for the associated players' values. Those expressions will lead us to a system of equations that we can use to turn the candidate equilibrium into an actual one.

The memoryless property of the Markovian structure of $(X,M)$ is hiding the fact that a game started at some earlier (negative) time and observed at time zero may have already ended as a result of the randomisation. Each player's value at time $0$ given $(X_0,M_0) = (x,m)$ needs to be calculated conditional on the game still being in play at time $0$. Due to the continuity of $\Gamma$, this is equivalent to conditioning on $\bar \tau_\Gamma > 0$. Notice that $\tau_\Gamma$ may be zero with positive probability when $\Gamma(m) < 1$, which corresponds to the game having been ended before time $0$.

For Player 1 there are two values, depending on the realisation of $\theta$. In game theory it is customary to think of Player 1 as having two
types: type $0$ ($\theta=0$) and type $1$ ($\theta=1$). Each type has her own value. Recall the first passage times $T_0$ and $T_1$ for the process $X$ at zero and one, respectively. Under $\P^0_{x,m}$ the event $\{\sigma_d=\H_1\}$ has zero
probability because the process $X$ never visits $x=1$. That, combined with the condition $h(0)=g(0)$, yields the following value for type $0$ of Player 1:
\[
u^0(x,m)=\E^0_{x,m}\big[\cP(\sigma_d,\He)\big| \bar\tau_{\Gamma} > 0\big]=
\E^0_{x,m}\big[\cP(\sigma_d,\He)\big| \H_0 > 0\big] =
\E^0_{x,m}\big[g\big(X_{\sigma_d}\big)\big],\quad 0 < x\le m < 1,
\]
where we used that $\P^0_{x,m} (\bar \tau_\Gamma = \H_0 ) = 1$ and $\P^0_{x,m} (\H_0 > 0) = 1$.

The value of type $1$ of Player 1 reads, for $0 < x\le m < 1$,
\begin{equation}\label{eqn:tilu}
\begin{aligned}
&u^1(x,m)=\E^1_{x,m}\big[\cP(\sigma_d,\bar\tau_\Gamma) \big| \bar\tau_\Gamma > 0 \big]\\
&=
\E^1_{x,m}\big[\ind_{\{\tau_\Gamma > 0\}} \big(g(X_{\sigma_d})\ind_{\{\sigma_d<\tau_\Gamma\}}+f(X_{\tau_\Gamma})\ind_{\{\tau_\Gamma<\sigma_d\}} +h(X_{\tau_\Gamma})\ind_{\{\tau_\Gamma=\sigma_d\}}\big)\big] / \P^1_{x,m} (\tau_\Gamma > 0)\\
&=\E^1_{x,m}\Big[g(X_{\sigma_d})\frac{\Gamma(M_{\sigma_d})}{\Gamma(m)}-\int_{[0,\sigma_d)}f(X_t)\frac{\ud\Gamma(M_t)}{\Gamma(m)} + \ind_{\{\sigma_d = \H_1\}} \frac{\Gamma(1)}{\Gamma(m)} \big(f(1)-g(1)\big)\Big],
\end{aligned}
\end{equation}
where for the third equality we used \eqref{eq:useful} combined with
$\E^1_{x,m}\big[\ind_{\{\sigma_d<\tau_\Gamma\}}\big|\cF^2_\infty\big]
=\Gamma(M_{\sigma_d})\ind_{\{\sigma_d<T_1\}}=\Gamma(M_{\sigma_d})-\Gamma(1)\ind_{\{\sigma_d=T_1\}}$,
$h(1) = f(1)$ and the fact that $\P^1_{x,m}(0 < \tau_\Gamma=\sigma_d < \He) = 0$ because $t \mapsto \xi_t$ is continuous for $t < \He$. 

For Player 2's value we first notice that, for $0 < x \le m < 1$,
\begin{equation}\label{eq:v0}
\begin{aligned}
v(x,m)&=\E_{x,m}\big[\cP(\sigma_d,\bar\tau_{\Gamma})\big| \bar\tau_{\Gamma} > 0\big]=\E_{x,m}\big[\cP(\sigma_d,\He)\ind_{\{\theta=0\}}+\cP(\sigma_d,\tau_\Gamma)\ind_{\{\theta=1\}}\big| \bar \tau_\Gamma > 0\big]\\
&=\E_{x,m}\big[g(X_{\sigma_d})\ind_{\{\theta=0\}} +\cP(\sigma_d,\tau_\Gamma)\ind_{\{\tau_\Gamma > 0\}} \ind_{\{\theta=1\}}\big] / \P_{x,m}(\bar\tau_\Gamma > 0),
\end{aligned}
\end{equation}
where for the second equality we used that $\{\sigma_d = \He\} \cap \{\theta = 0\} = \{\sigma_d = \H_0\} \cap \{\theta = 0\}$ together with $h(0) = g(0)$.
Now we consider the two terms under expectation separately. For the first one we notice that 
by the tower property,
\begin{equation*}
\E_{x,m}\big[g(X_{\sigma_d})\ind_{\{\theta=0\}} \big]
=
\E_{x,m}\big[g(X_{\sigma_d}) (1-X_{\sigma_d}) \big].
\end{equation*}

For the second term, noting that $\P_{x,m}(\tau_\Gamma = \sigma_d < \He)=0$ we have
\begin{equation*}
\begin{aligned}
&\E_{x,m}\big[\cP(\sigma_d,\tau_\Gamma)\ind_{\{\tau_\Gamma > 0\}}\ind_{\{\theta=1\}}\big]\\
&=\E_{x,m}\big[\big(g(X_{\sigma_d})\ind_{\{\sigma_d<\tau_\Gamma\}}+f(X_{\tau_\Gamma})\ind_{\{\tau_\Gamma<\sigma_d\}}
+h(X_{\tau_\Gamma})\ind_{\{\tau_\Gamma=\sigma_d = \He\}}\big)\ind_{\{\tau_\Gamma > 0\}} \ind_{\{\theta=1\}}\big]\\
&=\E_{x,m}\Big[\big(g(X_{\sigma_d})\Gamma(M_{\sigma_d})\ind_{\{\sigma_d<T_1\}}-\int_{[0,\sigma_d)}f(X_t)\ud
\Gamma(M_t)+h(1) \Gamma(1) \ind_{\{\sigma_d = \H_1\}} \big)\ind_{\{\theta=1\}}\Big]\\
&=\E_{x,m}\Big[\big(g(X_{\sigma_d})\Gamma(M_{\sigma_d})-\int_{[0,\sigma_d)}f(X_t)\ud
\Gamma(M_t) + \Gamma(1) \big(f(1) - g(1)\big) \ind_{\{\sigma_d = \H_1\}} \big)\ind_{\{\theta=1\}}\Big]\\
&=\E_{x,m}\Big[\big(g(X_{\sigma_d})\Gamma(M_{\sigma_d})-\int_{[0,\sigma_d)}f(X_t)\ud
\Gamma(M_t) + \Gamma(1) \big(f(1) - g(1)\big) \ind_{\{\sigma_d = \H_1\}}\big)\P_{x,m}(\theta=1|\cF^2_{\sigma_d})\Big]\\
&=\E_{x,m}\Big[X_{\sigma_d}g(X_{\sigma_d})\Gamma(M_{\sigma_d})-X_{\sigma_d}\int_{[0,\sigma_d)}f(X_t)\ud
\Gamma(M_t) + \Gamma(1) \big(f(1) - g(1)\big) \ind_{\{\sigma_d = \H_1\}}\Big],
\end{aligned}
\end{equation*}
where we used \eqref{eq:useful} for the second equality, $f(1) = h(1)$ and the set inclusions $\{\sigma_d<\tau_\Gamma\}\subset\{\tau_\Gamma>0\}$ and $\{\tau_\Gamma=\sigma_d=\He\}\subset\{\tau_\Gamma>0\}$ for the third one, the tower property for the fourth equality, and $\P_{x,m}(\theta=1|\cF^2_{\sigma_d}) = X_{\sigma_d}$ in the fifth.

As for the case of the two types of Player 1's values, Player 2's value also incorporates the probability that the game started at some earlier (negative) time
has not been terminated by Player 1 by time zero. We use that $\P_{x,m}(\bar\tau_\Gamma > 0)=\P_{x,m}(\theta=0) + \P_{x,m}(\tau_\Gamma > 0)\P_{x,m}(\theta=1) = 1 - x + \Gamma(m) x$. Combining the expressions found above with \eqref{eq:v0} yields the value function  of the uninformed player: for $0 < x\le m < 1$,
\begin{equation}\label{eq:vu0u1}
\begin{aligned}
v(x,m)&=\E_{x,m}\Big[g(X_{\sigma_d})\big(1-X_{\sigma_d}+X_{\sigma_d}\Gamma(M_{\sigma_d})\big)-X_{\sigma_d}\int_{[0,\sigma_d)}f(X_t)\ud
\Gamma(M_t)\\
&\hspace{100pt} + \Gamma(1) \big(f(1) - g(1)\big) \ind_{\{\sigma_d = \H_1\}}\Big] / \Big(1 - x + \Gamma(m) x\Big)\\
&= \frac{1-x}{1-x+x\Gamma(m)}u^0(x,m)+\frac{x\Gamma(m)}{1-x+x\Gamma(m)}u^1(x,m),
\end{aligned}
\end{equation}
where the second equality holds by the changes of probability measure in \eqref{eq:P0P1}.
Noting that under the measure $\P_x=\hP_{x}^{x}$ we have that $\hat p_t$ from \eqref{eq:hatpX} is given by $\hat p_t=X_t$, and substituting this expression into \eqref{eqn:pt} along with $\xi^0_t=0$ and $\xi^1_t = 1-\Gamma(M_t)$ we find that
\begin{equation}\label{eq:ptM}
\begin{aligned}
p_t=\P_{x}\big(\theta=1\big|\cF^2_t,\bar\tau_\Gamma>t\big)
=\frac{X_t\Gamma(M_t)}{1-X_t+X_t\Gamma(M_t)}= \pi_\Gamma(X_t,M_t),
\end{aligned}
\end{equation}
where
\begin{equation}\label{eq:pi}
\pi_\Gamma(x,m)\coloneqq\frac{x\Gamma(m)}{1-x+x\Gamma(m)},\quad 0\le x\le m < 1.
\end{equation}
To allow for arbitrary values of $M_0$ we extend \eqref{eq:ptM} to hold with $\P_x$ replaced by $\P_{x,m}$.

In order to connect with the results from \cite{de2025nash} (see Section 2.2 therein for an outline of relevant results and
Section 3 for details) we notice that if $(\sigma_d,\bar \tau_\Gamma)$ is a saddle point we
can identify the processes $u^0(X_t,M_t)$, $u^1(X_t,M_t)$ with the process $V^{*,1}(t)$ in \cite{de2025nash}, restricted to the
event $\{\theta=0\}$ and $\{\theta=1\}$, respectively. Likewise, we can identify the process $v(X_t,M_t)$ with the process
$V^{*,2}(t)$ in \cite{de2025nash}.
With this in mind, \cite{de2025nash} tells us that for a pair $(\sigma_d,\bar \tau_\Gamma)$ to be a saddle point it is necessary
and sufficient that the following conditions hold: for Player 1's two types' values it must be (cf. \cite[Propositions 3.7, 3.8 and 3.17-(i)]{de2025nash}) for $0<x\le m<1$\label{page:postulates}
\begin{equation*}
\begin{aligned}
(i)&\quad t\mapsto u^0(X_{t\wedge\sigma_d},M_{t\wedge\sigma_d})\,\,\text{is a $\P^0_{x,m}$-martingale with}\,\,
u^0(X_{t\wedge\sigma_d},M_{t\wedge\sigma_d}) < f(X_{t\wedge\sigma_d}),\,\,\text{for all $t\ge 0$},\\
(ii)&\quad t\mapsto u^1(X_{t\wedge\sigma_d},M_{t\wedge\sigma_d})\,\,\text{is a $\P^1_{x,m}$-submartingale with}\,\,
u^1(X_{t\wedge\sigma_d},M_{t\wedge\sigma_d})\le f(X_{t\wedge\sigma_d}),\,\,\text{for all $t\ge 0$},\\
(iii)&\quad t\mapsto u^1(X_{t\wedge\sigma_d\wedge\tau_\Gamma(1)},M_{t\wedge\sigma_d\wedge\tau_\Gamma(1)})\,\,\text{is
$\P^1_{x,m}$-martingale for $t\ge 0$},
\end{aligned}
\end{equation*}
where we recall that $\tau_\Gamma(1)\coloneqq\lim_{u\uparrow 1}\tau_\Gamma(u)=\inf\{t\ge 0:1-\Gamma(M_t)=1\}\wedge \He$ and,
moreover, see \cite[Proposition 3.17(ii)]{de2025nash})
\begin{equation}
\begin{aligned}
(iv)&\quad\int_{[0,\sigma_d)}\big(u^1(X_{t},M_{t}) - f(X_{t})\big)\ud \Gamma(M_t)=0,\quad\P^1_{x,m}-a.s.\\
(v)&\quad u^0(X_{\sigma_d},M_{\sigma_d})=g(X_{\sigma_d}),\,\text{$\P^0_{x,m}$-a.s.\
and}\,\,u^1(X_{\sigma_d},M_{\sigma_d})=g(X_{\sigma_d}),\,\text{$\P^1_{x,m}$-a.s.\ on $\{\sigma_d < \He\}$.}
\end{aligned}
\end{equation}
In particular, the second part of (ii), combined with (iv) imply that Player 1 only stops with intensity/probability $\ud
\Gamma(M_t)$ at times contained in the set $\{t\ge 0:\,u^1(X_{t},M_{t}) = f(X_{t})\}$. Continuing the statement about necessary
and sufficient conditions for the saddle point, for Player 2's value it must be
\begin{equation*}
\begin{aligned}
(vi)&\quad t\mapsto
\big(1-X_{t\wedge\tau_\Gamma(1)}+X_{t\wedge\tau_\Gamma(1)}\Gamma(M_{t\wedge\tau_\Gamma(1)})\big)v(X_{t\wedge\tau_\Gamma(1)},M_{t\wedge\tau_\Gamma(1)})-\int_{[0,{t\wedge\tau_\Gamma(1)})}f(X_s)X_s\ud
\Gamma(M_s),\\
&\text{is a $\P_{x,m}$-supermartingale for $t \ge 0$ and it is a martingale for $t \le \sigma_d$ (\cite[Prop.\ 3.7 and 3.8]{de2025nash})},\\
(vii)&\quad v(X_{t\wedge\tau_\Gamma(1)},M_{t\wedge\tau_\Gamma(1)}) \ge g(X_{t\wedge\tau_\Gamma(1)}),\,\,\text{for all $t\ge 0$, $\P_{x,m}$-a.s.\ (\cite[Proposition 3.17-(i)]{de2025nash})},
\end{aligned}
\end{equation*}
and, finally, we have the equation connecting the values of the two players \eqref{eq:vu0u1}: for $t \leq \He$
\begin{equation}
(viii)\quad v(X_t,M_t)=\pi_\Gamma(X_t,M_t)u^1(X_t,M_t)+\big(1-\pi_\Gamma(X_t,M_t)\big)u^0(X_t,M_t),\quad \P_{x,m}-a.s.
\end{equation}
\label{page:end-conditions}


\subsection{A system of free boundary equations}
Conditions (i)--(viii) lead us to the formulation of a system of free boundary problems, described below. Solving the system will allow us to determine the optimal strategies expressed via a pair $(d_*,\Gamma_*)$ and functions $(u_*^0,u_*^1,v_*)$. As a result we obtain a saddle point $(\sigma_*=\sigma_{d_*},\bar{\tau}_*=\bar\tau_{\Gamma_*})$ and players' values.

To begin with we work with a general pair $(d,\Gamma)$ with $d:[0,1]\to[0,1]$ and $\Gamma:(0,1]\to(0,1]$, where 
$d$ and $\Gamma$ are chosen from certain sets of regular functions. Often $d$ is of the form $d(m) = c(m \vee n) \wedge m$, where $n \in (0,1)$ and $c$ is chosen from another set of regular functions. 

\begin{definition}
\
\begin{enumerate}
\item Let $\cC$ be the set of functions $c:[0,1]\to[0,1]$ 
with the properties that (i) $c$ is continuous, increasing, and continuously differentiable on the set $\{m\in(0,1)\colon m^2 f'(m)>c^2(m)g'(c(m))\}$, (ii) $0<  c(m) < m$ on $(0,1]$, (iii) $c(m) \leq x_g$ and (iv) $c(1)g'(c(1)) < f'(1)$.
\item Let $\cD$ be the set of functions $d:[0,1]\to[0,1]$ with the properties that (i) $d$ is continuous and increasing, (ii) $0 \leq d(m) \leq m$ and (iii) there exists $m_d \in [0,1)$ such that $d(m)=m$ for $m \in [0,m_d]$ and $0 <d(m) < m$ for $m \in (m_d,1]$. We let $\cM_d =(m_d,1)$
\item Let $\cG$ be the set of functions $\Gamma:(0,1]\to (0,1]$ with the properties that (i) $\Gamma$ is non-increasing and absolutely continuous, and (ii) $\Gamma$ is strictly decreasing and continuously differentiable on the set $\{ m\in(0,1) \colon \Gamma(m)<1 \}$.  
\end{enumerate}
\end{definition}

For $c \in \cC$ we define (cf.\ \eqref{eqn:m0})
\begin{equation} \label{eq:m0cdef}
m_0(c) = \sup \{ m \in [0,1)\colon m^2 f'(m) \le c(m)^2 g'(c(m)) \},\quad \sup\varnothing=0. 
\end{equation}
Note that for $c \in \cC$ we have $f'(1)>c(1) g'(c(1)) > c(1)^2 g'(c(1))$ so that by the continuity of $c$, $m_0(c) \in [0,1)$. 
For $\Gamma \in \cG$ we define
\begin{equation}\label{eq:endpoints}
m_\Gamma^1 = \inf \{m \in (0,1)\colon \Gamma(m) < 1 \}
\end{equation}
and $\cM_\Gamma = (m^1_\Gamma,1)$. Then $\Gamma=1$ on $(0,m^1_\Gamma]$ and $\Gamma<1$ on $\cM_\Gamma$ with $\Gamma\in C^1(\cM_\Gamma)$.

The requirement $c(m) < m$ for all $m \in \cC$ is in line with properties of $c_*(m)$ in Assumption \ref{ass:fghscaled} 
and reflects our conviction that acting optimally players would never stop simultaneously, because of the second-mover advantage in the game (and the fact that Player 1 only stops on the diagonal).  

We suppose that
\[
d \in \cD,\quad \Gamma \in \cG, \quad m_d \leq m^1_\Gamma.
\]
First we consider what happens if $m_d\le X_0 \le M_0$:
the maximiser stops the game at $\sigma_d = \inf \{ u : X_u \leq d(M_u) \}$; the informed minimiser never stops the game if $\theta=0$ and stops the game at rate $- \ud \ln \Gamma(M_s) = -(\Gamma'(M_s)/\Gamma(M_s)) \ud M_s$ if $\theta=1$. If $X_0 \le M_0 \leq m_d$ then the game ends immediately.

Recall $u^0(x,m) = \E^0_{x,m}[g(X_{\sigma_d})]$
so that on $t < \sigma_d \wedge \He$ we have that $\tilu^0(X_t,M_t) = \E^0_{x,m}[g(X_{\sigma_d})|\cF^2_t]$ by standard Markovian arguments.
Then, the function $\tilu^0$ solves
\begin{equation}\label{eq:PDEu0}
\begin{cases}
\frac 12 \partial_{xx}\tilu^0(x,m)-\frac{1}{1-x}\partial_x \tilu^0(x,m)=0,&\quad d(m)<x<m,\, m \in \cM_d, \\ 
\partial_m \tilu^0(m,m)=0,&\quad m \in \cM_d, \\ 
\tilu^0(x,m)=g\big(x\big),&\quad 0 < x \leq d(m), \, m \in (0,1), 
\end{cases}
\end{equation}
where the first and second equations in \eqref{eq:PDEu0} represent the martingale property in the continuation region with the second equation the so-called normal reflection condition; and the third equation in \eqref{eq:PDEu0} follows from the stopping by the maximiser. Furthermore, from condition (i) on p.\ \pageref{page:postulates} we obtain the condition:
\begin{equation}\label{eqn:u0f}
\tilu^0(x,m) < f(x), \quad 0 < x < m.
\end{equation}

Now recall $\tilu^1(x,m) = \E^1_{x,m}[\cP(\sigma_d,\tau_\Gamma)|\bar\tau_\Gamma>0]$ and its equivalent representation in the third line of \eqref{eqn:tilu}. Then, by the strong Markov property, on $t<\He$ we have
\begin{equation}\label{eq:Mgu1}
\begin{aligned}
&\Gamma(M_{t\wedge \sigma_d})u^1(X_{t\wedge \sigma_d},M_{t\wedge \sigma_d})-\int_{[0,t\wedge \sigma_d)}f(X_s)\ud \Gamma(M_s)\\
&=\E^1_{x,m}\Big[g(X_{\sigma_d})\Gamma(M_{\sigma_d})-\int_{[0,\sigma_d)}f(X_s)\ud \Gamma(M_s)+\ind_{\{\sigma_d=T_1\}}\Gamma(1)\big(f(1)-g(1)\big)\Big|\cF^2_{t\wedge\sigma_d}\Big].
\end{aligned}
\end{equation}
The expression on the right-hand side is a martingale with initial value $\Gamma(m)u^1(x,m)$.
Using $\ud\Gamma(M_s) = \Gamma'(M_s) \ud M_s$, we apply It\^{o} formula to $\Gamma(M_s) \tilu^1(X_s, M_s)$ and obtain that
\begin{align*}
t \mapsto &\int_0^{t\wedge \sigma_d} \Gamma(M_s) \left( \tfrac12 \partial_{xx} \tilu^1(X_s,M_s) + \tfrac1x \partial_{x} \tilu^1(X_s,M_s)\right) \ud s \\
& + \int_0^{t\wedge \sigma_d} \Gamma(M_s) \partial_{x} \tilu^1(X_s,M_s) \ud W_s \\
& + \int_{[0, t\wedge \sigma_d)} \left[ \left( \Gamma(M_s) \partial_m \tilu^1 (X_s,M_s) + \Gamma'(M_s) \tilu^1(X_s, M_s) \right) \ud M_s - f(X_s) \Gamma'(M_s) \ud M_s \right], 
\end{align*}
is a zero-mean $(\F^2, \P^1_{x,m})$-martingale. Hence, various terms must be zero. Since $\Gamma > 0$, the first integral yields the first line of \eqref{eq:PDEu1} below. The third integral yields
\begin{equation}\label{eqn:tilu_condition}
\big(\tilu^1(M_s,M_s) - f(M_s) \big)\Gamma'(M_s)+ \Gamma(M_s) \partial_{m} \tilu^1(M_s,M_s) =0, \quad \text{for $s\in\mathrm{supp}(\ud M)$,}
\end{equation}
which combined with the postulated condition (iv) on p.\ \pageref{page:postulates} for $\tilu^1$ gives the second and third line in \eqref{eq:PDEu1}. The last line of the PDE problem \eqref{eq:PDEu1} arises directly from \eqref{eq:Mgu1} when $\sigma_d=0$
and it aligns with the postulated condition (v). It then follows that $\tilu^1$ satisfies
\begin{equation}\label{eq:PDEu1}
\begin{cases}
\tfrac 12 \partial_{xx} \tilu^1(x,m)+\tfrac{1}{x}\partial_x \tilu^1(x,m)=0,&\quad d(m)<x<m,\, m \in \cM_d \\ 
\partial_m \tilu^1(m,m) = 0,& \quad m \in \cM_d, \\
u^1(m,m)=f(m),&\quad m\in (m_\Gamma^1,1],\\
\tilu^1(x,m)=g\big(x\big),&\quad 0 < x \leq d(m), \, m \in (0,1). 
\end{cases}
\end{equation}
Analogously as for $\tilu^0$, item (ii) implies the condition
\begin{equation}\label{eqn:u1f}
\tilu^1(x,m) < f(x), \quad 0 < x < m.
\end{equation}

We now turn attention to the system of equations associated to $\tilv(x,m) = \E_{x,m}[\cP(\sigma_d,\bar{\tau}_\Gamma)|\bar{\tau}_\Gamma>0]$ as in the first expression on the right-hand side of \eqref{eq:vu0u1}. Notice that 
\[
X_{\sigma_d}\int_{[0,\sigma_d)}f(X_s)\ud
\Gamma(M_s)=\int_{[0,\sigma_d)}f(X_s)X_s\ud
\Gamma(M_s)+\int_0^{\sigma_d}\Big(\int_{[0,t)}f(X_s)\ud
\Gamma(M_s)\Big)\ud X_t,
\]
and the final integral is a martingale, because $X$ is a Brownian
motion under $\P_{x,m}$. So in the integral representation of $v$ in \eqref{eq:vu0u1} we can substitute the integral on the left-hand side above with the first one on the right-hand side.
Again by the strong Markov property, for $t<\sigma_d$ we have
\begin{equation*}
\begin{aligned}
&\big(1-X_{t\wedge\tau_\Gamma(1)}+X_{t\wedge\tau_\Gamma(1)}\Gamma(M_{t\wedge\tau_\Gamma(1)})\big)
\tilv(X_{t\wedge\tau_\Gamma(1)},M_{t\wedge\tau_\Gamma(1)})-\int_{[0,t\wedge\tau_\Gamma(1))}f(X_s)X_s\ud
\Gamma(M_s)\\
&=\E_{x,m}\Big[g(X_{\sigma_d})\big(1-X_{\sigma_d}+X_{\sigma_d}\Gamma(M_{\sigma_d})\big)-\int_{[0,\sigma_d)}f(X_s)X_s\ud
\Gamma(M_s)\Big|\cF^2_{t\wedge\tau_\Gamma(1)}\Big]\\
&\quad+\E_{x,m}\Big[\ind_{\{\sigma_d=T_1\}}\Gamma(1)\big(f(1)-g(1)\big)\Big|\cF^2_{t\wedge\tau_\Gamma(1)}\Big].
\end{aligned}
\end{equation*}
On the right-hand side above we have a martingale with initial value $(1-x+x\Gamma(m))v(x,m)$.
Assuming that $v$ is sufficiently smooth and applying It\^o's formula to the left-hand side 
\begin{equation}\label{eq:vIto}
\begin{aligned}
t\mapsto&\int_0^{t\wedge\tau_\Gamma(1)}\big[\big(1-X_s+X_s\Gamma(M_s)\big)\partial_x \tilv(X_s,M_s)-(1-\Gamma(M_s))\tilv(X_s,M_s)\big]\ud
X_s\\
&+\tfrac12\int_0^{t\wedge\tau_\Gamma(1)}\big[\big(1-X_s+X_s\Gamma(M_s)\big)\partial_{xx}
\tilv(X_s,M_s)-2\big(1-\Gamma(M_s)\big)\partial_x \tilv(X_s,M_s)\big]\ud s\\
&+\int_0^{t\wedge\tau_\Gamma(1)}\Big(X_s\big(v(X_s,M_s)-f(X_s)\big)
\Gamma'(M_s)+\big(1-X_s+X_s\Gamma(M_s)\big)\partial_{m} \tilv(X_s,M_s)\Big)\ud M_s,
\end{aligned}
\end{equation}
is a zero-mean $(\F^2,\P_{x,m})$-martingale.
Then, the first two equations below hold and the last equation follows by taking $\sigma_d=0$ in \eqref{eq:vu0u1}: 
\begin{equation}\label{eq:PDEv}
\begin{cases}
\tfrac12\big(1\!-\!x\!+\!x\Gamma(m)\big)\partial_{xx}\tilv(x,m)\!-\!\big(1\!-\!\Gamma(m)\big)\partial_x \tilv(x,m)=0,&\, d(m)<x <
m,\, m \in \cM_d, \\ 
m\big(\tilv(m,m)\!-\! f(m)\big)\Gamma'(m)\!+\!\big(1\!-\!m\!+\!m\Gamma(m)\big)\partial_m \tilv(m,m)=0, &m \in \cM_d, \\ 
\tilv(x,m)=g(x),&\, 0<x\le d(m),\, m \in (0,1). 
\end{cases}
\end{equation}
Postulate (vii) gives an obstacle condition
\begin{equation}\label{eq:vgL}
v(x,m)\ge g(x), \qquad 0<x\le m,\, m\in \cM_d.
\end{equation}
Once the stopping rule of Player 1 is fixed, Player 2 is solving a standard optimal stopping problem and it is therefore natural to also add the smooth-fit condition:
\begin{equation}\label{eq:smf1}
\partial_x v(d(m),m) =  g'(d(m)), \quad m\in \cM_d,
\end{equation}
which is needed to also justify the use of It\^o's formula in \eqref{eq:vIto}.
Furthermore, taking $d(m)$ as the stopping boundary, we also expect that on the continuation set $x > d(m)$ the inequality in \eqref{eq:vgL} is strict. Finally, the super-harmonic property must hold
\begin{equation}\label{eq:vsh}
\tfrac12\big(1-x+x\Gamma(m)\big)\partial_{xx}v(x,m)-\big(1-\Gamma(m)\big)\partial_x v(x,m)\le0, \quad  0<x < m,\ m \in \cM_d,
\end{equation}
because of the super-martingale condition in (vi).

In addition to the above systems of PDEs we must ensure that (viii) holds, and therefore we must impose for $0<x\le m$, $m\in(m_d,1)$,
\begin{equation}\label{eq:link}
\tilv(x,m)=\pi_\Gamma(x,m)\tilu^1(x,m)+(1-\pi_\Gamma(x,m))\tilu^0(x,m),
\end{equation}
where $\pi_\Gamma(x,m)$ is defined in \eqref{eq:pi}. The above equation is the analogue of \eqref{eq:vu0u1}, which was derived from
purely probabilistic considerations.

Luckily, the number of equations can be reduced as shown in the next lemma.  
\begin{lemma}\label{lem:vfromu}
Fix an arbitrary admissible pair $(d,\Gamma)$. 
Assume that there are functions $u^0$ and $u^1$ solving systems \eqref{eq:PDEu0} and \eqref{eq:PDEu1}, respectively.
Then, $v$ defined by \eqref{eq:link} satisfies \eqref{eq:PDEv}.
\end{lemma}
\begin{proof} 
For the first equation in \eqref{eq:PDEv} we introduce
\begin{equation}\label{eqn:tildev}
\begin{aligned}
\hat v(x,m)\coloneqq \big(1-x+x\Gamma(m)\big)v(x,m) =(1-x) u^0(x,m)+x \Gamma(m) u^1(x,m). 
\end{aligned}
\end{equation}
We then obtain, for $d(m)<x \leq m$
\[
\partial_{xx} \hat v(x,m) = 2 \Gamma(m) \partial_x \tilu^1(x,m) - 2 \partial_x u^0(x,m) + x \Gamma(m) \partial_{xx} \tilu^1(x,m) + (1-x)
\partial_{xx} \tilu^0(x,m).
\]
Substituting $\partial_{xx} \tilu^0$ and $\partial_{xx} \tilu^1$ from \eqref{eq:PDEu0} and \eqref{eq:PDEu1}, we obtain that $\partial_{xx}
\hat v(x,m) = 0$ for $x \in (d(m), m)$.
Using the latter and differentiating the identity $\hat v(x, m) = (1 - x + x \Gamma(m)) \tilv(x,m)$ twice with respect to $x$ yields
\[
0=\partial_{xx} \hat v(x,m) = 2\big(\Gamma(m) - 1\big) \partial_x \tilv(x,m) + \big(1-x+x\Gamma(m)\big) \partial_{xx} \tilv(x,m).
\]
This proves the first equation of \eqref{eq:PDEv}.

The third equation is immediate from the construction because $\tilu^0(x,m)=\tilu^1(x,m)=g(x)$ for $0<x\leq d(m)$. For the second equation, we differentiate $\hat v$ with respect to $m$ and obtain
\[
\partial_m \hat v(x,m) = x \Gamma'(m) \tilu^1(x,m) + x \Gamma(m) \partial_m \tilu^1(x,m) + (1-x) \partial_m \tilu^0(x,m).
\]
Inserting $x=m$ and using the expressions for $\partial_m \tilu^0(m,m)$ and $\partial_m \tilu^1(m,m)$ yields
\begin{equation}\label{eqn:m_tildev}
\partial_m \hat v(m,m) = m \Gamma'(m) f(m).
\end{equation}
This equality and differentiation of identity $\hat v(x, m) = (1 - x + x \Gamma(m))\tilv(x,m)$ with respect to $m$ gives the third
equality of \eqref{eq:PDEv}.
\end{proof}

Notice that the function $\hat v (x,m) = (1-x+x\Gamma(m)) v(x,m)$ defined in the above proof is affine in $x$ for $x \in
(d(m), m)$. This will prove useful in the verification argument in Section \ref{sec:verif}.


\section{Solution of the system}\label{sec:sol}

In this section we identify the optimal $d \in \cD$ and $\Gamma \in \cG$ which determine the optimal stopping rules $\sigma_d$ and $\bar\tau_\Gamma$. We do this by solving the equations introduced in the previous section. Key to our construction is the function $c_*$ from Assumption \ref{ass:fghscaled}-(v) which is obtained as the root of $\phi(\cdot,m)$. 

To aid understanding of this section, we provide below a sketch of our construction. In what follows we will \emph{derive} the quantities that the reader finds below.
\begin{itemize}[label={},leftmargin=2em, itemsep=0.5\baselineskip]
\item First, we notice that
\begin{equation}\label{eq:c*} 
c_*(m)(f(m)-g(c_*(m))) = (m-c_*(m)) m f'(m)
\end{equation}
and recall the definition of $m_0(c_*)= \sup \{ m : m^2 f'(m) \leq c_*(m)^2 g'(c_*(m))\}$; $c_*$ is the first candidate stopping rule for the maximiser and $c_* \in \cC$;
\item Second, given $c_*$ define $\Du0_{c_*}:(0,1) \to {\mathbb R}$ as the solution of a differential equation
\begin{equation} \label{eq:Ddef*}
\Du0_{c_*}'(m) = \frac{\big(1-c_*(m)\big)}{(1-m)\big(m-c_*(m)\big)} \big[\Du0_{c_*}(m)-g\big(c_*(m)\big)\big]\quad \text{subject to
$\Du0_{c_*}(1-) = \lim_{m \uparrow 1} \Du0_{c_*}(m) = g\big(c_*(1)\big)$};
\end{equation}
it will then follow that $\Du0_{c_*}(m) = \E^0_{m,m}[g(X_{\sigma_{c_*}})]$ and that $U_{c_*}(m)>g(c_*(m))$ on $(0,1)$;
\item Third, define $\gamma_{c_*}:(m_0(c_*),1) \to \R$ via
\begin{equation}
\label{eq:gamma*def} 	
\gamma_{c_*}(m)= \frac{c_*(m)}{\big(m-c_*(m)\big)}\frac{  \big(1-c_*(m)\big)\big(m-c_*(m)\big) g'\big(c_*(m)\big) -
(1-m)\big[{\Du0_{c_*}(m)}-g\big(c_*(m)\big)\big]  } {m^2 f'(m) - c_*(m)^2 g'\big(c_*(m)\big)}.
\end{equation}
and note that the denominator of $\gamma_{c_*}$ is positive on $(m_0(c_*),1)$; in addition, the numerator of \eqref{eq:gamma*def} can be shown to be positive on $(m_0(c_*),1)$ so that $\gamma_{c_*}$ restricted to $(m_0(c_*),1)$ is a map $\gamma_{c_*}: (m_0(c_*),1) \to (0,\infty)$; moreover $\gamma_{c_*}(1) <1$ and $\gamma_{c_*}$ is decreasing on $(m^1_{\gamma_{c_*}},1)$ where $m^1_{\gamma_{c_*}} = \sup \{ n \in [0,1) : \gamma_{c_*}(n) \geq 1 \} \geq m_0(c_*)$;
\item Fourth, define $\Gamma_*$ by $\Gamma_* = \gamma_{c_*} \wedge 1$ and set $m^1_* = m^1_{\gamma_{c_*}}$;
\item Fifth, define $d_*$ by $d_*(m) = c_*(m \vee m^1_{*}) \wedge m$.
\end{itemize}
Then the pair $(d_*, \Gamma_*)$ defines the optimal strategy via the stopping rules $\sigma_{d_*}$ and $\bar\tau_{\Gamma_*}$. 

As part of the construction we define $u^0_*$, $u^1_*$ and $v_*$ which satisfy \eqref{eq:PDEu0}, \eqref{eq:PDEu1} and \eqref{eq:PDEv}. Indeed, the initial choice $c=c_*$ for the maximiser's optimal stopping boundary is determined by the requirement that $u^1$ should satisfy the normal reflection condition $\partial_m u^1(m,m)=0$. The initial choice $\Gamma = \gamma_{c_*}$ for the minimiser's intensity of stopping is determined by \eqref{eq:smf1} and by the requirement of smooth-fit of $v$ at $c_*(m)$. It turns out that typically $\gamma_{c_*}$ does not take values only in $(0,1]$ and so does not directly define an element $\Gamma_* \in \cG$. Instead we must define $\Gamma_{c_*} = \gamma_{c_*} \wedge 1$. This means that the optimal stopping rules must be modified for $M_t \in (0,m^1_{\gamma_{c_*}})$. This motivates the introduction of $d_*$ as the optimal threshold for the maximiser -- this threshold is constant for $m \in (c_*(m^1_*),m^1_*)$ to reflect that for $M_t \leq m^1_*$ neither type of the minimiser stops and thus the maximiser can learn nothing from their behaviour.

As a final check, we need to verify that the resulting value functions satisfy the obstacle conditions \eqref{eqn:u0f}, \eqref{eqn:u1f} and \eqref{eq:vgL} and the super-harmonic property \eqref{eq:vsh}.

First we study the solution of the system \eqref{eq:PDEu0}. 
In determining $u^0$, $\Gamma$ plays no role. We begin by studying the case of general $c$ and  construct $u^0_c$ which has the properties required for \eqref{eq:PDEu0}, but over a different domain.

\begin{proposition}\label{prop:u0}
Fix $c \in \cC$ and define
$\Du0_c:(0,1) \to {\mathbb R}$ by
\begin{equation} \label{eq:Ddef}
\Du0_c'(m) = \frac{\big(1-c(m)\big)}{(1-m)\big(m-c(m)\big)} \big[\Du0_c(m)-g\big(c(m)\big)\big]\quad \text{subject to
$\Du0_c(1-) = \lim_{m \uparrow 1} \Du0_c(m) = g\big(c(1)\big)$}.
\end{equation} 
We extend $\Du0_c$ to $(0,1]$ by setting $\Du0_c(1)=\Du0_c (1-)$. 

For $0\le x \leq c(m)$ define $u^0_c(x,m)=g(x)$. For $c(m) < x \le m$, define
\begin{equation}\label{eq:u0def}
u^0_c(x,m)=g\big(c(m)\big) + \frac{\big(x-c(m)\big)\big(1-m\big)}{(1-x)\big(m-c(m)\big)} \big[\Du0_c(m)-g\big(c(m)\big)\big].
\end{equation}
Then $u^0_c$ is continuous and it solves $\cL^0 u^0_c(x,m) = 0$ on $c(m)<x<m$ for all $m \in (0,1)$, and $\partial_m u^0_c(m,m)=0$ for $m\in(m_0(c),1)$.  Moreover, $u^0_c(m,m)=U_c(m)$ for all $m \in (0,1)$, and $u^0_c(x,m) = \E_{x,m}[g(X_{\sigma_c})]$ for $m\in(m_0(c),1)$. 
\end{proposition}
\begin{proof}
Note that by l'H\^{o}pital's rule and the formula for $\Du0'_c$ in \eqref{eq:Ddef}, if $\Du0_c(1-)=g(c(1))$ then $\Du0'_c(1-)=\tfrac12 c'(1) g'(c(1))$. Since the only singularity for the right-hand side of \eqref{eq:Ddef} could occur in $m=1$ it follows that 
the solution to \eqref{eq:Ddef} is bounded on $(0,1)$.

Continuity of $u^0_c$ in both $x$ and $m$ is immediate from the continuity of $c$ and $U_c$ and the fact that letting $x \downarrow c(m)$ in \eqref{eq:u0def} $u^0_c(c(m)+,m) = g(c(m))$.
Further, since $(1-x)u^0_c(x,m)$ is affine in $x$ for fixed $m$ we find that $\cL^0 u^0_c(x,m) = 0$. It is also clear that $u^0_c(m,m) = \Du0_c(m)$.

Taking the derivative of \eqref{eq:u0def} with respect to $m$ at $x=m$, and using that $c$ is continuously differentiable on $(m_0(c),1)$, we find that any coefficient of $c'(m)$ disappears and then
\begin{equation*}
\begin{aligned}
\partial_m u^0_c(m,m) & =  \partial_m \left. \left(\frac{(1-m)(x-c)}{(m-c)(1-x)} (\Du0_c(m) - g(c)) \right) \right|_{x=m,c=c(m)} \\
& =  \Du0_c'(m) - \frac{1-c(m)}{(1-m)(m-c(m))} \big[\Du0_c(m) - g(c(m)) \big]
\end{aligned}
\end{equation*}
Then, since $\Du0_c$ solves \eqref{eq:Ddef} we find that $u^0_c$ satisfies the normal reflection condition $\partial_m u^0_c(m,m) = 0$.

Since $\Du0_c$ is bounded, the Feynman-Kac formula gives that $u^0_c(x,m)=\E^0_{x,m}[g(X_{\sigma_c})]$ for $c(m)\le x\le m$ and $m\in(m_0(c),1)$. 
\end{proof}

\begin{corollary}\label{cor:upbu0}
For $c \in \cC$ and $m \in (m_0(c),1)$ and $0<x \leq m$ we have $u^0_c(x,m)\le g(x) < f(x)$. 
\end{corollary}
\begin{proof} 
Note that outside $c(m) < x \leq m$  the statement $u^0_c(x,m)\le g(x) < f(x)$ is immediate since $u^0_c(x,m) = g(x)$.

For $m_0(c)<m<1$ and $c(m)<x \leq m$, we have 
$u^0_c(x,m)=\E^0_{x,m}[g(X_{\sigma_c})]$ with $\sigma_c$ as in \eqref{eq:sigmac}.
Then, using It\^o's formula we get
\[
u^0_c(x,m)=g(x)+\E^0_{x,m}\Big[\int_0^{\sigma_{d}} \cL^0 g(X_t)\ud t\Big]\le g(x) < f(x),
\]
where the first inequality is by Assumption \ref{ass:fghscaled}-(iv) and the second one by (i) in the same Assumption.
\end{proof}

We now study $u^1$. Again at a first pass we get a solution which satisfies the properties of \eqref{eq:PDEu1} but on a different domain.

\begin{lemma}\label{lem:u1c}
Fix $c \in \cC$. 
For $m\in (0,1)$ and $c(m) < x\le m$, define
\begin{equation}\label{eq:u1defa}
u^1_c(x,m)=g(c(m)) + \frac{m}{x}\frac{\big(x-c(m)\big)}{\big(m-c(m)\big)} \big[f(m)-g\big(c(m)\big)\big] .
\end{equation}
For $0 < x \leq  c(m)$ define $u^1_c(x,m) = g(x)$.

Then $u^1_c$ is continuous and for all $m \in (0,1)$, $u^1_c$ solves $\cL^1 u^1_c (x,m)= 0$ on $c(m)<x<m$.
\end{lemma}

\begin{proof}
A function $u^1$ is a solution of the first equation in \eqref{eq:PDEu1} if and only if $xu^1(x,m)$ is an affine function of $x$ and then the boundary conditions at $m$ and $c(m)$ yield \eqref{eq:u1defa}. 
\end{proof}

The next proposition motivates the choice of $c_*$ from Assumption \ref{ass:fghscaled}-(v) as a candidate stopping boundary for Player 2. Note that $c_* \in \cC$ by Lemma~\ref{lem:c*}. Note further that $c_*$ is defined on $(0,1]$, but later we will modify the optimal stopping boundary so that we only take $c=c_*$ on $(m^1_\gamma,1)$ where $m^1_\gamma$ is defined in \eqref{eq:mgamma1def} below.

\begin{lemma}\label{lem:u1*}
Suppose that $c=c_*$ where $c_*(m)$ is the root of $\oldh(\cdot,m)$. Then for $m \in (m_0(c_*),1)$ we have that $u^1_{c_*}$ solves $\partial_m u^1_{c_*}(m,m)=0$.
\end{lemma}

\begin{proof}
In Lemma~\ref{lem:c*} we proved differentiability and monotonicity of $c_*$ at least on the set $(m_0(c_*),1)$. 
Then, taking $c=c_*$, differentiating \eqref{eq:u1defa} with respect to $m$ 
and evaluating at $x=m$ yields
\[
\partial_m u^1_{c_*}(m,m) = f'(m) - \frac{c_*(m)}{m}\cdot \frac{f(m)-g(c_*(m))}{m-c_*(m)}.
\]
Since $\oldh(c_*(m),m) = 0$ we have that $\partial_m u^1_{c_*}(m,m) = 0$. Note that $c_*$ is the only choice for which this identity will hold.
\end{proof}

The next step is to identify the form of $\Gamma$ which guarantees the smooth-fit condition \eqref{eq:smf1}.
Later we will apply the result at $c=c_*$ but to begin with we work with a general  
$c \in \cC$.
Recall that 
for $c \in \cC$ we have $m_0(c) \in [0,1)$.

Define $\gamma_{c}: (m_0(c),1] \to {\mathbb R}$ by 
\begin{equation}
\label{eq:gammadefgen} 	
\gamma_c(m)= \frac{c(m)}{\big(m-c(m)\big)}\frac{  \big(1-c(m)\big)\big(m-c(m)\big) g'\big(c(m)\big) -
(1-m)\big[{\Du0_{c}(m)}-g\big(c(m)\big)\big]  } {m^2 f'(m) - c(m)^2 g'\big(c(m)\big)}.
\end{equation}
Note that the denominator on the right-hand-side of \eqref{eq:gammadefgen} is positive by the definition of $m_0(c)$.
The next result shows that the numerator is also positive on $(m_0(c),1]$ and hence that $\gamma_c$ is a map $\gamma_c : (m_0(c),1) \to (0,\infty)$.
Note that if $m_0(c)>0$ then $\lim_{u\downarrow 0}\gamma(m_0 + u)=+\infty$.
Then if we define
\begin{equation}
\label{eq:mgamma1def}
m^1_\gamma\coloneqq\sup\{m\in[0,1): \gamma(m)\ge 1\}, \quad\text{with $\sup\varnothing=0$}.
\end{equation}
we have that $m^1_{\gamma_c} \ge m_0(c)$.

\begin{lemma}\label{lem:sign}
We have $\gamma_c(m)>0$ for $m \in (m_0(c),1]$.
\end{lemma}

\begin{proof} 
For $m\in(m_0(c),1]$, the condition $\gamma_c(m)>0$ is equivalent to
\[
\big(1-c(m)\big)\big(m-c(m)\big)g'\big(c(m)\big)>\big(1-m\big)\big[{\Du0_{c}(m)}-g\big(c(m)\big)\big].
\]
In the limit $m \uparrow 1$ the strict inequality holds because the right-hand side vanishes (recall that $\Du0_{c}(1)=g(c(1))$ from \eqref{eq:Ddef}), whereas on the
left-hand
side we are left with $(1-c(1))^2g'(c(1))$ which is strictly positive because $c \in \cC$ and hence $c(1)<1$.

Let us now consider $m\in(m_0(c),1)$. 
Recall that $\Du0_{c}(m)=u^0(m,m)\le g(m)$, where the inequality is by Corollary \ref{cor:upbu0}.
Then the statement of the lemma is implied by the stronger inequality
\begin{equation}\label{eq:ineqg}
\big(1-c(m)\big)\big(m-c(m)\big)g'\big(c(m)\big)>\big(1-m\big)\big[g(m)-g\big(c(m)\big)\big].
\end{equation}
Now, for $x\in[0,m]$ let
\[
\varphi_m(x):=(1-x)(m-x)g'(x)-(1-m)\big(g(m)-g(x)\big).
\]
Clearly $\varphi_m(m)=0$, then \eqref{eq:ineqg} holds if we can prove that $\varphi'_m(x) < 0$ for $x\in(0,m)$. Direct
computations yield
\begin{equation*}
\begin{aligned}
\varphi'_m(x)&=-(m-x)g'(x)-(1-x)g'(x)+(1-x)(m-x)g''(x)+(1-m)g'(x)\\
&=-2(m-x)g'(x)+(1-x)(m-x)g''(x)\\
&=2(1-x)(m-x)\big(\tfrac12 g''(x)-\tfrac{1}{1-x}g'(x)\big)< 0,
\end{aligned}
\end{equation*}
where the inequality holds by Assumption \ref{ass:fghscaled}-(iv) which says that $\cL^0g < 0$. 
\end{proof}

Recall the definition of $\pi_\Gamma$ in \eqref{eq:pi}. 
We extend this definition from $\Gamma \in \cG$ to arbitrary positive functions $\gamma:(0,1] \to (0,\infty)$
by setting $\pi_\gamma(x,m) = x\gamma(m)/(1-x+x \gamma(m))$. If $\gamma(m)>1$ then $\pi_\gamma$ cannot be interpreted as a probability. Later we restrict attention to $\gamma$ taking values in $(0,1]$ and this interpretation is restored.

For $\gamma>0$  and for $m > m_0(c)$ and $0 < x \leq m$ define
\begin{equation}
\label{eq:vcgammadef}
 v_{c,\gamma}(x,m)=\pi_{\gamma}(x,m) u^1_{c}(x,m)+(1-\pi_{\gamma}(x,m))u^0_{c}(x,m) 
 \end{equation}
Note that the fact that $\gamma$ is positive is sufficient for $\pi_\gamma$ and therefore $v_{c,\gamma}$ to be well defined.

\begin{proposition}\label{prop:smoothfit}
Take $\gamma_c$ as in \eqref{eq:gammadefgen}. For $m\in(m_0(c),1)$, 
\[
\partial_x v_{c, \gamma_c}(c(m),m)=g'(c(m)),\quad m\in  (m_0(c),1).
\]
\end{proposition}

\begin{proof} 
By taking derivatives with respect to $x$ of $v_{c,\gamma_c}$ 
we obtain
\begin{equation*}
\begin{aligned}
\partial_x v_{c,\gamma_c}(x,m)&=\partial_x \pi_{\gamma_c}(x,m)\big(u^1_{c}(x,m)-u^0_{c}(x,m)\big)\\
&\quad+\pi_{\gamma_c}(x,m)\partial_x u^1_{c}(x,m)+\big(1-\pi_{\gamma_c}(x,m)\big)\partial_x
u^0_{c}(x,m).
\end{aligned}
\end{equation*}
Taking $x=c(m)$ the first term on the right-hand side above vanishes because $u^1_{c}(c(m),m)=u^0_{c}(c(m),m)=g(c(m))$. Direct
calculations yield
\begin{equation*}
\begin{aligned}
\partial_x u^0_{c}(x,m)&=\frac{1-m}{m-c(m)}\big[\Du0_{c}(m)-g\big(c(m)\big)\big]\frac{1-c(m)}{(1-x)^2} , \\
\partial_x u^1_{c}(x,m)&=\frac{m}{m-c(m)}\big[f(m)-g\big(c(m)\big)\big]\frac{c(m)}{x^2}.
\end{aligned}
\end{equation*}
Then we find 
\begin{equation*}
\begin{aligned}
\partial_x v_{c,\gamma_c}(c(m),m) & =  \frac{\gamma_c(m)}{1-c(m)+c(m)\gamma_c(m)} \frac{m}{m-c(m)}\big[f(m)-g(c(m))\big] \\ 
&\quad 
+\frac{1}{1-c(m)+c(m)\gamma_c(m)}\ \frac{1-m}{m-c(m)}\big[\Du0_{c}(m)-g(c(m))\big].
\end{aligned}
\end{equation*}
Substituting the expression in \eqref{eq:gammadefgen} for $\gamma_c(m)$ we obtain $\partial_x v_{c,\gamma_c}(c(m),m)=g'(c(m))$.
\end{proof}

At this stage $\gamma_c$ is a candidate choice for $\Gamma \in \cG$ as it yields the smooth fit condition for $v$. However, it is not necessarily $(0,1]$-valued and it is not necessarily decreasing and therefore we cannot use it directly. Indeed, we expect $\gamma_c$ to have a singularity at $m_0(c)$.
We can deal with the fact that $\gamma_c$ does not take only values in $(0,1]$ by setting $\Gamma_c= \gamma_c \wedge 1$ 
on $(m_0(c),1]$. Moreover, since $c(1)<1$ we have  
\begin{equation}\label{eq:Gamma1gen}
\gamma_c(1)\coloneqq\lim_{m\uparrow
1}\gamma_c(m)=\frac{c(1)\big(1-c(1)\big)g'\big(c(1)\big)}{f'(1)-c^2(1)g'\big(c(1)\big)}\in(0,1),
\end{equation}
where the upper bound $\gamma_c(1)<1$ holds 
by the definition of $\cC$ and the fact that $c(1) g'(c(1)) < f'(1)$.
Then we can define $m^1_{\gamma_c} = \sup \{ n \in [0,1) : \gamma_c(n) \geq 1 \}$. We have 
$0 \leq m_0(c) \leq m^1_{\gamma_c} < 1$ and we can extend the definition of $\Gamma_c$ to $(0,1]$ by setting
\[ 
\Gamma_c(m) = 
\begin{cases}  
1, & 0 < m \leq m^1_{\gamma_c}, \\
\gamma_c(m), & m^1_{\gamma_c} < m < 1. 
\end{cases}
\]
 
It remains the case that $\Gamma_c$ may not be non-increasing. However, in the sequel we will be interested in $\Gamma_{c_*}$ and, as the next lemma shows, $\Gamma_c$ is decreasing for the choice $c=c_*$.

\begin{lemma}\label{lem:Gamma}
For $m\in(m_0(c_*),1)$, if $\gamma_{c_*}(m)<1$, then $\gamma_{c_*}'(m)<0$.
\end{lemma}

\begin{proof}
In the proof we simplify notations by letting $c_*(m)=c$, $\gamma_{c_*}(m)=\gamma$ and $m_0(c_*)=m_0$.
Take $m\in(m_0,1)$. Notice that
\begin{equation}\label{eq:loggamma}
\log \gamma = \log \left( \frac{c}{m-c} \right) + \log A(m,c) - \log B(m,c),
\end{equation}
where
\begin{equation*}
\begin{aligned}
A(m,c) & =  (1-c)(m-c)g'(c) - (1-m)(\Du0_{c_*}(m)-g(c)), \\
B(m,c) & =  m^2 f'(m) - c^2 g'(c).
\end{aligned}
\end{equation*}
Notice that $B(m,c)>0$ because $m\in(m_0,1]$ and $A(m,c)>0$ by the proof of Lemma \ref{lem:sign}.
In order to show $\frac{\ud }{\ud m}\log \gamma<0$, we calculate the expressions for the derivative of each term on the right-hand
side of \eqref{eq:loggamma}.

For the first term,
\[
\frac{\ud }{\ud m}  \left( \frac{c}{m-c} \right) = \frac{1}{(m-c)^2}\left( m c' - c \right).
\]
For the second term,
\begin{equation*}
\begin{aligned}
\frac{\ud}{\ud m} A(m,c) & =  (1-c)g'(c) - (1-m)\Du0_{c}'(m) + \Du0_{c}(m)-g(c) + c' (m-c)\big[ (1-c)g''(c) -2 g'(c) \big] \\
          & =  (1-c)g'(c) - \frac{(1-m)}{(m-c)} \big[ \Du0_{c}(m)-g(c) \big] + 2c'(1-c)(m-c) \big[ \tfrac12 g''(c) - \tfrac{1}{1-c} g'(c)
          \big] \\
          & =  \frac{1}{(m-c)} A(m,c)
              + 2c'(1-c)(m-c) \big[ \tfrac12 g''(c) - \tfrac{1}{1-c} g'(c) \big],
\end{aligned}
\end{equation*}
where in the second line we substituted the expression for $\Du0_{c}'$ from \eqref{eq:Ddef} and in the third one we substituted the
expression for $A(m,c)$. For the third term in \eqref{eq:loggamma}, it is useful to recall \eqref{eq:c'} for $c'$, which can be rewritten as
\[
mf''(m) + 2f'(m) = c' \Big( \frac{m^2 f'(m) - c^2 g'(c)}{c(m-c)} \Big)=c' \frac{B(m,c)}{c(m-c)}.
\]
Using that formula in the second line below, we obtain
\begin{equation*}
\begin{aligned}
\frac{\ud}{\ud m} B(m,c) & =  2m f'(m) + m^2 f''(m) - 2c^2 c' \left( \tfrac12 g''(c)+\tfrac1c g'(c) \right) \\
& =  c'\left[ \frac{m}{c(m-c)}B(m,c) -  2c^2\left(\tfrac12 g''(c)+\tfrac1c g'(c) \right) \right].
\end{aligned}
\end{equation*}
Then, differentiating \eqref{eq:loggamma} and using the formulae above we get
\begin{equation*}
\begin{aligned}
\frac{\gamma'}{\gamma}&=\frac{m-c}{c}\frac{mc'-c}{(m-c)^2}+\frac{1}{m-c}+\frac{2c'(1-c)(m-c)}{A(m,c)}\big[\tfrac12
g''(c)-\tfrac{1}{1-c}g'(c)\big]\\
&\quad-c'\Big[ \frac{m}{c(m-c)} -  \frac{2c^2}{B(m,c)}\left(\tfrac12 g''(c)+\tfrac1c g'(c) \right) \Big]\\
&=2 c'\frac{m-c}{A(m,c)}\Big[\tfrac12 (1-c)g''(c)-g'(c)+\frac{c^2}{m-c}\frac{A(m,c)}{B(m,c)}\big(\tfrac12 g''(c)+\tfrac1c
g'(c)\big)\Big]\\
&=2 c'\frac{m-c}{A(m,c)}\Big[\tfrac12 (1-c)g''(c)-g'(c)+c\gamma\big(\tfrac12 g''(c)+\tfrac1c g'(c)\big)\Big],
\end{aligned}
\end{equation*}
where in the final expression we used again the formula for $\gamma$ (cf.\ \eqref{eq:loggamma}).

Now, if $\tfrac12 g''(c)+\tfrac1c g'(c)\le 0$ we obtain $\gamma'\le 0$ irrespective of whether or not $\gamma<1$, thanks to Assumption
\ref{ass:fghscaled}-(iv). If instead\footnote{Notice that this is in line with the sufficient conditions for Lemma \ref{lem:jan}.} $\tfrac12 g''(c)+\tfrac1c g'(c)> 0$, then using $\gamma<1$ we obtain
\begin{equation*}
\begin{aligned}
\frac{\gamma'}{\gamma}& < 2 c'\frac{m-c}{A(m,c)}\Big[\tfrac12 (1-c)g''(c)-g'(c)+c\big(\tfrac12 g''(c)+\tfrac1c g'(c)\big)\Big]=2
c'\frac{m-c}{A(m,c)}\tfrac12g''(c)\le 0,
\end{aligned}
\end{equation*}
where the final inequality holds because $g$ is concave on $(0,x_g)$ and $c \leq x_g$.

Thus $\gamma(m)<1$ implies $\gamma'(m)<0$ as needed.
\end{proof}
In the light of the lemma we define $\Gamma_* = \Gamma_{c_*}$ by
\begin{equation}
\label{eq:Gamma*} 
\Gamma_{*}(m) = 
\begin{cases}  
1 & 0 < m, \leq m^1_{*}, \\
\gamma_{c_*}(m), & m^1_{*} < m < 1, 
\end{cases} 
\end{equation}
where $\gamma_{c_*}$ is defined in \eqref{eq:gamma*def} and $m^1_* \coloneqq m^1_{\Gamma_*} = m^1_{\gamma_{c_*}}$.

Combinining the choices $c = c_*$ and $\Gamma=\gamma_{c} \wedge 1$ evaluated at $c=c^*$ and defining $u^0_*= u^0_{c_*}$ 
as in \eqref{eq:u0def}, $u^1_* = u^1_{c_*}$ as in \eqref{eq:u1defa} and $v_*= v_{c_*,\gamma_{c_*}}$ as in \eqref{eq:vcgammadef} we find that for $m \in [m^1_*,1]$ and $c_*(m)\le x\le m$
\begin{align} 
u^0_*(x,m) & =  g\big(c_*(m)\big) + \frac{\big(x-c_*(m)\big)\big(1-m\big)}{(1-x)\big(m-c_*(m)\big)} \big[\Du0_{c_*}(m)-g\big(c_*(m)\big)\big], \label{eq:u0*def} \\
u^1_*(x,m) & = g(c_*(m)) + \frac{m}{x}\frac{\big(x-c_*(m)\big)}{\big(m-c_*(m)\big)} \big[f(m)-g\big(c_*(m)\big)\big], \label{eq:u1*def} \\
v_*(x,m) & =  \pi_{*}(x,m) u^1_{*}(x,m)+(1-\pi_{*}(x,m))u^0_{*}(x,m),
 \label{eq:v*def}
\end{align}
where $\pi_*=\pi_{\Gamma_*}$. All functions are continuously extended to $0\le x\le c_*(m)$ by taking $u^0_*=u^1_*=g$.
Next we show that for $m \in [m^1_*,1]$ the obstacle conditions \eqref{eqn:u0f}, \eqref{eqn:u1f} and \eqref{eq:vgL} are satisfied along with the super-harmonic property \eqref{eq:vsh}. Later we extend the definitions of $u^0_*$, $u^1_*$ and $v_*$ to $m \in (0, m^1_*)$.  

It is convenient to introduce the set
\begin{equation}\label{eq:Oc*}
\cO_{c_*,\Gamma_*}\coloneqq\{(x,m):c_*(m)<x<m,\, m\in(m^1_*,1]\},
\end{equation}
and denote its closure by $\overline\cO_{c_*,\Gamma_*}$.

Notice that by construction
$u^0_*,u^1_*\in C^{2,1}(\overline \cO_{c_*,\Gamma_*})$, thanks to $c_*\in C^1((m^1_*,1])$. That is, in the open set $\cO_{c_*,\Gamma_*}$ functions $u^0_*$ and $u^1_*$ are twice continuously differentiable in $x$ and once in $m$, with both functions and their first and second derivatives admitting continuous extensions to the closure. Moreover, $u^0_*$ and $u^1_*$ solve \eqref{eq:PDEu0} and \eqref{eq:PDEu1} on $\cO_{c_*,\Gamma_*}$. Their convex combination yields $v_*$ which solves \eqref{eq:PDEv} and satisfies the smooth-fit condition \eqref{eq:smf1} for $m\in(m^1_*,1)$ (cf.\ Proposition \ref{prop:smoothfit}). 

\begin{lemma}\label{lem:apriori}
We have $u^0_* \le u^1_* \le f$ on
$\overline\cO_{c_*,\Gamma_*}$ with $u^1_*(x,m)<f(x)$ for $x\in[c_*(m),m)$, $m\in(m^1_*,1]$.
Moreover, $u^0_*(x,m)\ge g(c_*(m))$ for $(x,m)\in\overline\cO_{c_*,\gamma}$ and, if $m_0=0$, $\lim_{m\to 0}u^0_*(m,m)=g(0)$.
\end{lemma}
\begin{proof}
To prove $u^1_*(x,m)\le f(x)$ for $(x,m) \in \overline\cO_{c_*,\Gamma_*}$, let us start by noticing that direct differentiation of \eqref{eq:u1*def} yields
\[
\partial_x u^1_*(x,m) = \frac{m }{x^2} \frac{c_*(m)}{m-c_*(m)} \big[f(m)-g(c_*(m))\big].
\]
The above expression is continuous in $\overline\cO_{c_*,\Gamma_*}$ and along the diagonal it yields $\partial_x u^1_*(m,m)=f'(m)$ (cf.\ \eqref{eq:c*}). 

Let us set $\psi(x,m)=u^1_*(x,m)-f(x)$. Then
\[
\tfrac12 \partial_{xx}\psi(x,m)+\tfrac1x \partial_x \psi(x,m)=-\tfrac12 f''(x)-\tfrac1x f'(x)=-\cL^1 f(x),
\]
with boundary condition $\psi(m,m)=\partial_x\psi(m,m)=0$. Therefore,
\[
\tfrac12 \partial_{xx}\psi(m,m)=-\cL^1 f(m)<0,
\]
where the inequality holds by Assumption \ref{ass:fghscaled}-(iv).
Then there is $\eps>0$ such that
\[
\psi(x,m)<0,\quad x\in(m-\eps,m).
\]
Setting $\bar x:=\sup\{x\in[c_*(m),m)\colon\psi(x,m)\ge 0\}$ (with $\sup\varnothing=-\infty$) we have $\bar x \le m-\eps$. If $\bar x\ge
c_*(m)$
we reach a contradiction because by Feynman-Kac formula, for $x\in(\bar x,m)$
\[
0\ge \psi(x,m)=\E^1_{x,m}\Big[\int_0^{\H_{\bar x}\wedge \H_m}\cL^1 f(X_t)\ud t\Big]>0,
\]
where we recall that $\H_{a} = \inf\{t\ge 0: X_t=a\}$.
Thus $\bar x=-\infty$ and $\psi(x,m)=u^1_*(x,m)-f(x)<0$ for $x\in[c_*(m),m)$.
Now we proceed to showing $u^0_*\le u^1_*$. Recall that $f(m) > g(m) \geq u^0_*(m,m)=\Du0_{c_*}(m)$ for $m\in[m^1_*,1]$ by Corollary \ref{cor:upbu0}. 
Then, from the expression for $u^1_*$ in \eqref{eq:u1*def} we have
\begin{equation}
\begin{aligned}
u^1_*(x, m) & \geq  g(c_*(m)) + \frac{m}{x} \frac{(x - c_*(m))}{(m - c_*(m))} \big(\Du0_{c_*}(m) - g(c_*(m))\big) \\
& \geq  g(c_*(m)) + \frac{1-m}{1-x} \frac{(x - c_*(m))}{(m - c_*(m))} \big(\Du0_{c_*}(m) - g(c_*(m))\big) =  u^0_*(x,m),
\end{aligned}
\end{equation}
where the second inequality holds because $\frac{m}{x} \geq \frac{1-m}{1-x}$. This completes the proof of the first statement in the lemma.

For the remaining two statements of the lemma, let us recall that $[m^1_*,1]\ni m\mapsto c_*(m)$ is non-decreasing. Then, by the Feynman-Kac representation of $u^0_*$ (cf.\ Proposition \ref{prop:u0}) we have
\begin{equation*}
\begin{aligned}
u^0_*(x,m)&=\E^0_{x,m}\big[g(X_{\sigma_{c_*}})\big]=\E^0_{x,m}\big[g\big(c_*(M_{\sigma_{c_*}})\big)\big]\ge g\big(c_*(m)\big),
\end{aligned}
\end{equation*}
where the second equality is because $X_{\sigma_{c_*}}=c_*(M_{\sigma_{c_*}})$ a.s.\ under $\P^0_{x,m}$ with $c_*(m) \le x \le m$, and  the final inequality trivially holds because
${M_{\sigma_{c_*}}}\ge M_0=m$, $\P_{x,m}$-a.s. This proves the first part of the final statement in the lemma.

For the last remaining claim, recall that $m\mapsto c_*(m)$ is continuous on $[0,1]$ and $c_*(0) = 0$, by Lemma \ref{lem:c*}. As above,
$u^0(m,m)=\E^0_{m,m}\big[g(X_{\sigma_{c_*}})\big]$.
Recalling that $X$ has dynamics
\[
X_t =m+W^0_t- \int_0^t\frac{1}{1-X_s} \ud s,
\]
under $\P^0_{m,m}$, its scale function $s^0$ is given by $s^0(x)
= 1/(1-x)$. Then for arbitrary $m_1\in(0,1)$ and for $m \in(0, m_1)$ we have
\[
\P^0_{m,m}(\H_{m_1} < \sigma_{c_*}) < \P^0_{m,m}(\H_{m_1} < \H_0) = \frac{m}{1-m}
\frac{1-m_1}{m_1}.
\]
Using that $g(0)\le g(x) \leq g(1)$ and $g$ is increasing,
\begin{equation*}
\begin{aligned}
g(0)\le u^0_*(m,m)
&= \E^0_{m,m}\big[ g(X_{\sigma_{c_*}}) \big] = \E^0_{m,m}\big[ g(X_{\sigma_{c_*}}) \ind_{\{\H_{m_1} < \sigma_{c_*}\}}\big] + \E^0_{m,m}\big[
g(X_{\sigma_{c_*}}) \ind_{\{\H_{m_1} \ge \sigma_{c_*}\}}\big]\\
&\le g(1)\P^0_{m,m}\big( \H_{m_1} < \sigma_{c_*}\big) + g\big(c_*(m_1)\big) \P^0_{m,m}\big( \H_{m_1} \ge \sigma_{c_*}\big)\\
&\le g(1)\frac{m}{1-m} \frac{1-m_1}{m_1} + g(c_*(m_1)).
\end{aligned}
\end{equation*}
Hence we find that $g(0)\le \limsup_{m \to 0} u^0_*(m,m) \leq g(c_*(m_1))$. Since $m_1$ is arbitrary and $c_*$ and $g$ are continuous with
$g(0)=g(c_*(0))$, we conclude that $\lim_{m \to 0} u^0_*(m,m) = g(0)$.
\end{proof}

In the next proposition we verify that the conditions in \eqref{eq:vgL} are satisfied. 

\begin{proposition}\label{prop:obst}
For $m\in(m^1_*,1]$ we have $v_*(x,m)\ge g(x)$ for $0\le x\le m$.
\end{proposition}
\begin{proof}
By construction $v_*(x,m)= g(x)$ for $0\le x\le c_*(m)$.
We prove $v_*(x,m)\ge g(x)$ for $c_*(m)< x \le m$. Taking derivatives with respect to $m$ in \eqref{eq:v*def} we have
\begin{equation}\label{eq:vm}
\begin{aligned}
\partial_m v_*(x,m)  &=  x (1-x) \gamma_{c_*}'(m) \frac{(u^1_*(x,m) - u^0_*(x,m))}{(1-x+x \gamma_{c_*}(m))^2}\\
&\quad+ \pi_*(x,m)\partial_m u^1_*(x,m) +
\big(1-\pi_*(x,m)\big) \partial_m u^0_*(x,m).
\end{aligned}
\end{equation}
By Lemma \ref{lem:Gamma}, $\gamma_{c_*}'(m)<0$ for $m\in(m^1_*,1]$, and by Lemma \ref{lem:apriori} $u^0_* \le u^1_*$. Thus the first term on the right-hand side above is non-positive. 

For the second term, once again we use our shorthand notation $c=c_*(m)$ and $\gamma= \gamma_{c_*}(m)$.
Taking derivatives with respect to $m$ in \eqref{eq:u1*def} yields
\begin{equation}
\begin{aligned}
\partial_m u^1_*(x,m) & =  c'  \left[ g'(c) \left(1 - \frac{m(x-c)}{(m-c)x} \right) + [f(m)-g(c)] \frac{\partial}{\partial c} \left(
\frac{m(x-c)}{(m-c)x} \right) \right] \\
&\quad+ \frac{m(x-c)}{(m-c)x} f'(m) + [f(m)-g(c)] \frac{\partial}{\partial m} \left( \frac{m(x-c)}{(m-c)x}  \right).
\end{aligned}
\end{equation}
The term on the second line above vanishes using \eqref{eq:c*} and
\[
\frac{\partial}{\partial m} \left( \frac{m(x-c)}{(m-c)x}  \right)=-\frac{c}{(m-c)^2}\frac{x-c}{x}.
\]
For the term on the first line we notice that
\begin{equation}
\label{eq:partialc}
\frac{\partial}{\partial c} \left( \frac{x-c}{m-c} \right)=-\frac{m-x}{(m-c)^2},
\end{equation}
and use again \eqref{eq:c*} to arrive at
\begin{equation}\label{eq:u1m}
\partial_m u^1_*(x,m) = c' \frac{c(m-x)}{(m-c)x} \left( g'(c) - \frac{m^2}{c^2} f'(m) \right).
\end{equation}

Now we look at the third term on the right-hand side of \eqref{eq:vm}. From \eqref{eq:u0*def} we have
\begin{equation}
\begin{aligned}
\partial_m u^0_*(x,m) & =  c'  \left[ g'(c) \left(1 - \frac{(1-m)(x-c)}{(m-c)(1-x)} \right) + [\Du0_{c}(m)-g(c)]
\frac{(1-m)}{(1-x)}\frac{\partial}{\partial c} \left( \frac{x-c}{m-c} \right) \right] \\
&\quad+ \frac{(1-m)(x-c)}{(m-c)(1-x)} \Du0'_c(m) + [{\Du0_{c}(m)}-g(c)] \frac{(x-c)}{(1-x)} \frac{\partial}{\partial m} \left(
\frac{1-m}{m-c}
\right).
\end{aligned}
\end{equation}
The term on the second line vanishes because
\[
\frac{\partial}{\partial m} \left( \frac{1-m}{m-c}  \right)=-\frac{1-c}{(m-c)^2},
\]
and using the expression \eqref{eq:Ddef} for $\Du0'_c(m)$. Then, again using \eqref{eq:c*}, the definition of $\gamma_{c_*}(m)$ from
\eqref{eq:gamma*def} and \eqref{eq:partialc}
we can further simplify the term on the first line and obtain
\begin{equation}
\begin{aligned}
\partial_m u^0_*(x,m) & = c' \left[ g'(c) \frac{(m-x)(1-c)}{(m-c)(1-x)} - [\Du0_{c}(m)-g(c)] \frac{(1-m)(m-x)}{(1-x)(m-c)^2}  \right]
\\
& =  c' \frac{(m-x)}{c(m-c)(1-x)} \left[ \gamma(m) (m^2 f'(m) - c^2 g'(c) ) \right] \\
& =  - \frac{x}{1-x}\gamma(m) \partial_m u^1_*(x,m),
\end{aligned}
\end{equation}
where the final equality holds by \eqref{eq:u1m}.

Writing again $c_*$ and $\gamma_{c_*}$ in full,
the above calculations show that for $m \in (m^1_*,1)$
\[
\pi_*(x,m)\partial_m u^1_*(x,m) + \big(1-\pi_*(x,m)\big) \partial_m u^0_*(x,m)=0,
\]
and therefore from \eqref{eq:vm} we deduce $\partial_m v_*(x,m)\le 0$ for $c_*(m)<x \le m$. Since $c'_*>0$ on $(m^1_*,1)$
(Lemma \ref{lem:c*}), then it admits a continuous inverse and
\[
v_*(x,m) = v_*(x,c_*^{-1}(x)) - \int_m^{c_*^{-1}(x)} \partial_m v_*(x,y) \ud y \geq v_*(x,c_*^{-1}(x)) = g\big(c_*(c_*^{-1}(x))\big) =
g(x),
\]
concluding the proof of $v_*(x,m)\ge g(x)$ for $c_*(m)<x \leq c_*(1)$.
It remains to show the result for $c_*(1)<x \leq m<1$.

For $c_*(1)<x \leq m$ we have
\[
v_*(x,m) = v_*(x, 1) - \int_m^{1} \partial_m v_*(x,y) \ud y \geq v_*(x,1)
\]
so to show that $v_*(x,m) \geq g(x)$ for all $m \geq x$ for a given $c_*(1) \leq x < 1$ it is necessary and sufficient to show that
$v_*(x,1) \geq g(x)$.

Fix $c_*(1)<x\le m <1$. By construction (use \eqref{eq:v*def}) we have
\[
v(x,1) = g(c_*(1)) \left[ 1 - \frac{\gamma_{c_*}(1)}{(1-c_*(1))} \frac{x-c_*(1)}{1-x+x\gamma_{c_*}(1)} \right]  + f(1) \frac{\gamma_{c_*}(1)}{(1-c_*(1))} \frac{x-c_*(1)}{1-x+x\gamma_{c_*}(1)} \]
and it is necessary and sufficient to show that
\[
0 \leq v(x,1) - g(x) = \big[f(1)- g(c_*(1))\big] \frac{\gamma_{c_*}(1)}{(1-c_*(1))} \frac{x-c_*(1)}{1-x+x\gamma_{c_*}(1)} - \big[g(x)- g(c_*(1))\big],
\]
or equivalently
\[
\frac{f(1)- g(c_*(1))}{g(x)- g(c_*(1))} \geq  \frac{(1-c_*(1))}{\gamma_{c_*}(1)} \frac{1-x+x\gamma_{c_*}(1)}{x-c_*(1)}.
\]
Subtracting $1$ from both sides we find this is equivalent to
\begin{equation}
\label{eq:newversion}
\frac{\gamma_{c_*}(1)}{1 - (1- \gamma_{c_*}(1))c_*(1)} \frac{x - c_*(1)}{1-x} \geq \frac{g(x)- g(c_*(1))}{f(1)- g(x)}
\end{equation}
Since
\[
\gamma_{c_*}(1) = \frac{c_*(1)(1-c_*(1))g'(c_*(1))}{f'(1) - c_*(1)^2g'(c_*(1))},
\]
then we have
\[
1 - (1- \gamma_{c_*}(1))c_*(1) = \frac{(1-c_*(1))f'(1)}{f'(1) - c_*(1)^2 g'(c_*(1))}.
\]
Combining the two yields
\[
\frac{\gamma_{c_*}(1)}{1 - (1- \gamma_{c_*}(1))c_*(1)} = \frac{c_*(1)g'(c_*(1))}{f'(1)},
\]
and substituting this expression into
\eqref{eq:newversion} we obtain
\[ 
\frac{x-c_*(1)}{1-x} \frac{c_*(1) g'(c_*(1))}{f'(1)} \geq \frac{g(x)-g(c_*(1))}{f(1) - g(x)} \hspace{10mm} \forall x \in (c_*(1),1).
\]
The latter holds because of the second inequality in \eqref{eq:newcondition}.
\end{proof}

\begin{remark}
A careful read of the proof above reveals that when $c_*(1)<1$, the second inequality in \eqref{eq:newcondition} is indeed necessary and sufficient for $v(x,m)\ge g(x)$ for $x\in(c_*(1),1]$. Thus it cannot be weakened in our construction of the equilibrium.
\end{remark}

\begin{lemma}
The super-harmonic condition \eqref{eq:vsh} holds for $v=v_*$.
\end{lemma}
\begin{proof}
Thanks to Lemma \ref{lem:vfromu} 
we know that the condition holds for
$c_*(m)<x<m$, $m\in(m^1_*,1]$. So we only need to show that it holds for $0<x<c_*(m)$. The latter is trivial, because
$v_*(x,m)=g(x)$ in that set and
\[
\tfrac12\big(1-x+x\gamma_{c_*}(m)\big) g''(x)- \big(1-\gamma_{c_*}(m)\big)g'(x)\le0,
\]
because $g''\le 0$ on $(0,x_g]$, due to Assumption \ref{ass:fghscaled}-(iv), and $c_*(m)<x_g$ by Assumption \ref{ass:fghscaled}-(v).
\end{proof}

So far we have constructed a tuple $(c_*,\Gamma_*,u^0_*,u^1_*,v_*)$ that solves the systems \eqref{eq:PDEu0}-\eqref{eqn:u0f},
\eqref{eq:PDEu1}-\eqref{eqn:u1f} and \eqref{eq:PDEv}-\eqref{eq:smf1} with \eqref{eq:link} for $m\in[m^1_*,1]$. In
the remainder of the section we construct the solution on $[0,m^1_*)$. The analysis above suggests that for $m\in[0,m^1_*)$ neither of the two types of the informed player would care to stop (recall $\Gamma_*=1$ on $[0,m^1_*)$). Based on this insight, we construct equilibrium values for the game that are consistent with the form of $\Gamma_*$ and $c_*$ from \eqref{eq:Gamma*} and Assumption \ref{ass:fghscaled}-(v), respectively. That is, starting from $X_0 \leq M_0 \leq m^1_*$, 
both types of the informed player do not stop until $\H_{m^1_*}\wedge \He=\inf\{t\ge 0: X_t\notin (0,m^1_*)\}$; then, on
the event $\{\theta=1\}$ Player 1 stops according to $\Gamma_*(m)$ and on the event $\{\theta=0\}$ she never stops; the uninformed
player stops at $c^1_*:=c_*(m^1_*)$.

Instead of $c_*$ being the candidate optimal threshold we find that the optimal threshold is given by $d_*$ where 
\[
d_*(m)\coloneqq c_*(m\vee m^1_*) \wedge m.
\]
Notice that by doing so, the function $d_*$ is the unique root of $x\mapsto
\phi(x,m\vee m^1_*)$ in the interval $(0,m \wedge x_g)$, for all $m\in (c^1_*,1)$. Moreover, $d_*(m) = m$ for $m\in(0,c^1_*]$.

Given $d_*(m)$, let the definitions of functions $u^0_*$, $u^1_*$ from \eqref{eq:u0*def}--\eqref{eq:u1*def} be extended as follows:
for $j=0,1$,
\begin{equation}\label{eq:u0u1*}
u^j_*(x,m)=
\begin{cases}
u^j_{c_*}(x,m\vee m^1_*), & d_*(m)\le x\le m,\,m\in[c^1_*,1],\\
g(x), & 0< x \le m\le c^1_*\quad\text{or}\quad 0< x \le d_*(m),\, m\in[c^1_*,1].
\end{cases}
\end{equation}
As a result $v_*(x,m)$ is also extended to $0\le x\le m$, $m\in(0,c^1_*)$ according to \eqref{eq:v*def}.
Notice that
\[
v_*(x,m)=x u^1_*(x,m)+(1-x)u^0_*(x,m)=x u^1_*(x,m^1_*)+(1-x)u^0_*(x,m^1_*),\quad m\in (0,m^1_*],
\]
so that $v_*$ is constant in the $m$ variable for $m\in(0,m^1_*]$.
With a similar notation as in \eqref{eq:Oc*} we introduce the set
\begin{equation}\label{eq:Od}
\cO_*\coloneqq\{(x,m):d_*(m)<x<m,\, m\in(c^1_*,1)\},
\end{equation}
with closure denoted $\overline\cO_*$.
The next proposition collects properties of functions $u^0_*$, $u^1_*$ and $v_*$ which will be needed in the verification argument in Section \ref{sec:verif}.

\begin{proposition}\label{prop:smooth}
Functions $u^0_*$, $u^1_*$ and $v_*$ along with their derivatives $\partial_x u^0_*$, $\partial_{xx} u^0_*$, $\partial_x u^1_*$,
$\partial_{xx} u^1_*$, $\partial_x v_*$, $\partial_{xx} v_*$ are continuous on
$\overline{\cO}_{*}$.
The derivatives $\partial_m u^0_*$, $\partial_m u^1_*$ and $\partial_m v_*$ are continuous separately on
$\overline{\cO}_{*}\cap\{c^1_*\le m\le m^1_*\}$ and $\overline{\cO}_{*}\cap\{m^1_*\le m\le 1\}$.

Finally, separately on the sets $\cO_{*}\cap\{c^1_*< m<m^1_*\}$ and $\cO_{*}\cap\{m^1_*< m<1\}$, the function $u^0_*$ solves
\eqref{eq:PDEu0}-\eqref{eqn:u0f}, the function $u^1_*$ solves \eqref{eq:PDEu1}-\eqref{eqn:u1f} and the function $v_*$ solves
\eqref{eq:PDEv}-\eqref{eq:smf1}.
\end{proposition}
\begin{proof}
Continuity of the functions and of their $x$-derivatives is obvious by their definition and analogous properties of $u^0_{c_*}$, $u^1_{c_*}$
and $v_{c_*,\gamma_{c_*}}$. Moreover, it is clear that $\partial_m u^0_*=\partial_m u^1_*=\partial_m v_*=0$ on $\cO_{*}\cap\{c^1_*\le m\le m^1_*\}$.
Those derivatives may be discontinuous across $m=m^1_*$, but
\[
\partial_m u^0_*(m,m)=\partial_m u^1_*(m,m)=0,\quad\text{for all $m\in[c^1_*,1]$}.
\]

We have already shown that $u^0_{c_*}$, $u^1_{c_*}$
and $v_{c_*,\gamma_{c_*}}$ solve the relevant PDEs in $\cO_{*}\cap\{m^1_*<m<1\}$. Since
$u^j_*(x,m)=u^j_{c_*}(x,m^1_*)$ in $\cO_{*}\cap\{c^1_*\le m\le m^1_*\}$, for $j=0,1$, it then becomes clear that $u^1_*$ solves
\eqref{eq:PDEu1}-\eqref{eqn:u1f} also in $\cO_{*}\cap\{c^1_*< m< m^1_*\}$, upon observing that the second condition in
\eqref{eq:PDEu1} holds by the above paragraph and the third condition only needs to be evaluated at $(x,m)=(m^1_*,m^1_*)$
when $m\le m^1_*$. By analogous arguments we also see that $u^0_*$ and $v_*$ solve \eqref{eq:PDEu0}-\eqref{eqn:u0f} and
\eqref{eq:PDEv}-\eqref{eq:smf1}, respectively.
\end{proof}


\section{An equilibrium for the game}\label{sec:verif}

In this section we construct a perfect Bayesian equilibrium for our game (Theorem \ref{thm:mainscaled}).
We state two intermediate
results concerning the interpretation of functions $u^0_*$ and $u^1_*$, which will be useful for the proof of the main theorem.

Throughout the section we let $\rho_{a}\coloneqq \inf\{t\ge 0: M_t\ge a\}$ for $a\ge 0$, $\sigma_*\coloneqq\inf\{t\ge 0\colon X_t\le d_*(M_t)\}\wedge \He$ and $\tau_{\Gamma_*}(u)\coloneqq\inf\{t\ge 0: 1-\Gamma_*(M_t)>u\}\wedge \He$. Notice that $\rho_{m^1_*}=\tau_{\Gamma_*}(0)$ so we may use the two interchangeably.
\begin{lemma}\label{lem:semim}
Let $0<x\le m<1$.
For $j=0,1$, $\P^j_{x,m}$-a.s.\ we have
\begin{equation}\label{eq:semim}
u^j_*\big(X_{t\wedge \sigma_*},M_{t\wedge \sigma_*}\big)=u^j_*(x,m)+\int_0^{t\wedge
\sigma_*}\partial_xu^j_*\big(X_s,M_s\big)\ud W^j_s,\quad\text{for any $t\ge 0$}.
\end{equation}
\end{lemma}
\begin{proof}
The claim is trivial for $x\in(0,d_*(m)]$, because $\P^j_{x,m}(\sigma_*=0)=1$. Let us now consider $d_*(m)<x\le m<1$.
For any $t\ge 0$ we can write
\begin{equation}\label{eq:ito0}
\begin{aligned}
u^j_*(X_{t\wedge\sigma_*},M_{t\wedge\sigma_*})&=u^j_*(X_{t\wedge\sigma_*},M_{t\wedge\sigma_*})-u^j_*(X_{t\wedge\sigma_*\wedge\rho_{m^1_*}
},M_{t\wedge\sigma_*\wedge\rho_{m^1_*}})\\
&\quad +u^j_*(X_{t\wedge\sigma_*\wedge\rho_{m^1_*} },M_{t\wedge\sigma_*\wedge\rho_{m^1_*}}).
\end{aligned}
\end{equation}
We can apply It\^o's formula separately to the expressions in the two lines on the right-hand side above, because $u^j_*$ is
sufficiently regular on $\cO_{*}\cap\{c^1_*\le m\le m^1_*\}$ and on $\cO_{*}\cap\{m^1_*\le m\le 1\}$ separately. Under the
measure $\P^0_{x,m}$ we have
\begin{equation*}
\begin{aligned}
&u^0_*(X_{t\wedge\sigma_*},M_{t\wedge\sigma_*})-u^0_*(X_{t\wedge\sigma_*\wedge\rho_{m^1_*}
},M_{t\wedge\sigma_*\wedge\rho_{m^1_*}})\\
&=\int_{t\wedge\sigma_*\wedge\rho_{m^1_*}
}^{t\wedge\sigma_*}\Big(\tfrac12\partial_{xx}u^0_*\big(X_s,M_s\big)-\frac{1}{1-X_s}\partial_xu^0_*\big(X_s,M_s\big)\Big)\ud s\\
&\quad+\int_{t\wedge\sigma_*\wedge\rho_{m^1_*} }^{t\wedge\sigma_*}\partial_m u^0_*\big(X_s,M_s\big)\ud
M_s+\int_{t\wedge\sigma_*\wedge\rho_{m^1_*} }^{t\wedge\sigma_*}\partial_x u^0_*\big(X_s,M_s\big)\ud W^0_s\\
&=\int_{t\wedge\sigma_*\wedge\rho_{m^1_*} }^{t\wedge\sigma_*}\partial_x u^0_*\big(X_s,M_s\big)\ud W^0_s,
\end{aligned}
\end{equation*}
where the second equality holds thanks to \eqref{eq:PDEu0} and using
\[
\int_{t\wedge\sigma_*\wedge\rho_{m^1_*} }^{t\wedge\sigma_*}\partial_m u^0_*\big(X_s,M_s\big)\ud
M_s=\int_{t\wedge\sigma_*\wedge\rho_{m^1_*} }^{t\wedge\sigma_*}\partial_m u^0_*\big(M_s,M_s\big)\ud M_s=0
\]
By analogous arguments we also obtain
\begin{equation}
u^0_*(X_{t\wedge\sigma_*\wedge\rho_{m^1_*}
},M_{t\wedge\sigma_*\wedge\rho_{m^1_*}})=u^0_*(x,m)+\int_0^{t\wedge\sigma_*\wedge\rho_{m^1_*}}\partial_x u^0_*\big(X_s,M_s\big)\ud
W^0_s.
\end{equation}
Plugging those two expressions back into \eqref{eq:ito0} we obtain \eqref{eq:semim} for $j=0$. The same arguments with obvious
changes lead to \eqref{eq:semim} for $j=1$.
\end{proof}

Now we can use \eqref{eq:semim} to prove the next proposition. Notice that under $\P_{x,m}$ the stopping time $\H_m$ is the first
time the process $X$ hits $M_0=m$. Thus, $\H_m$ is the first time of increase of the process $M$. We are also going to use the
stopping time $\rho_{c^1_*}\coloneqq\inf\{t\ge 0: X_t\ge c^1_*\}$, where we recall $c^1_*=c_*(m^1_*)$. In particular,
it turns out that it is suboptimal for the $\theta=1$ type of Player 1 to act before $\H_m\vee \rho_{c^1_*} \vee \tau_{\Gamma_*}(0)$.
Indeed, $\P_{x,m}(\sigma_* = 0) = 1$ if $x \le c^1_*$, so stopping before $\rho_{c^1_*}$ would result in the sub-optimal middle payoff. Stopping before $\tau_{\Gamma_*}(0) = \rho_{m^1_*}$ would result in a payoff $f$ which may be strictly larger than $u^1_*$ as \eqref{eq:PDEu1} guarantees $u^1_*(m,m) = f(m)$ only for $m \ge m^1_*$. 

\begin{proposition}\label{prop:u0u1}
The functions $u^0_*$ and $u^1_*$ are the game's continuation values for the two types of Player 1. That is,
\begin{equation}\label{eq:ito3}
u^0_*(x,m)=
\E^0_{x,m}\Big[\cP\big(\sigma_*,\He\big)\Big]=\inf_{\tau_0\in\cT_R(\F^2)}\E^0_{x,m}\Big[\cP\big(\sigma_*,\tau_0\big)\Big],
\end{equation}
and
\begin{equation}\label{eq:ito8}
u^1_*(x,m)=\E^1_{x,m}\Big[\cP\big(\sigma_*,\tau_\Gamma\big)\Big]=\inf_{\tau_1\in\cT_R(\F^2)}\E^1_{x,m}\Big[\cP\big(\sigma_*,\tau_1\big)\Big],
\end{equation}
where $\tau_\Gamma\coloneqq\inf\{t\ge 0:1-\Gamma_t>U\}\wedge \He$ for any non-increasing $\F^2$-adapted \cadlag process $(\Gamma_t)$ with values in $(0,1]$ satisfying
$\Gamma_{\rho_{c^1_*}\vee \tau_{\Gamma_*}(0)} = 1$ and
\begin{equation}\label{eq:Gamma-t}
\int_0^1\ind_{\{X_{\tau_\Gamma(u)}\neq M_{\tau_\Gamma(u)}\}}\ud u =-\int_{[0,\He]} \ind_{\{X_t \ne M_t\}} \ud \Gamma_t=0, \quad \P^1_{x,m} \as
\end{equation}
\end{proposition}
\begin{proof}
The first equality in \eqref{eq:Gamma-t} follows from \cite[Ch.\ 0, Prop.\ 4.9]{revuzyor}. The second equality is by the definition of the process $(\Gamma_t)$.

Fix $0< x\le m < 1$ and take arbitrary $\tau\in\cT_R(\F^1)$ with decomposition
$\tau=\ind_{\{\theta=0\}}\tau_0+\ind_{\{\theta=1\}}\tau_1$ for some $\tau_0,\tau_1\in\cT_R(\F^2)$.
Recall the notation
$\tau_j(u)=\inf\{t\ge 0: \xi^j_t>u\}\wedge \He$, where $\xi^j\in\cA^\circ(\F^2)$ is the generating process of $\tau_j$. For any
$u\in(0,1)$, thanks to Lemma \ref{lem:semim} we have
\begin{equation}
u^j_*\big(X_{\tau_j(u)\wedge \sigma_*},M_{\tau_j(u)\wedge \sigma_*}\big)=u^j_*(x,m)+\int_0^{\tau_j(u)\wedge
\sigma_*}\partial_xu^j_*\big(X_s,M_s\big)\ud W^j_s.
\end{equation}
Taking expectations with respect to $\P^j_{x,m}$ and using Doob's optional sampling yields
\begin{equation}\label{eq:semim0}
u^j_*(x,m)=\E^j_{x,m}\Big[u^j_*\big(X_{\tau_j(u)\wedge \sigma_*},M_{\tau_j(u)\wedge \sigma_*}\big)\Big].
\end{equation}

In studying the right-hand side of the above equation we first consider $j=0$. Notice that on the event $\{\sigma_*\le
\tau_0(u)\}$ we have $X_{\sigma_*}\le d_*(M_{\sigma_*})$
under $\P^0$, because $\P^0_{x,m}(\sigma_*< \He)=1$ due to the fact that $d_*(m)>0$ for all $m\in(0,1]$ and $X$ does not hit
one. Recalling that $u^0_*(x,m)=g(x)$ for $x\le d_*(m)$ and using \eqref{eqn:u0f} we deduce
\begin{equation}\label{eq:ito1}
\begin{aligned}
u^0_*\big(X_{\tau_0(u)\wedge\sigma_*},M_{\tau_0(u)\wedge\sigma_*}\big)&=\ind_{\{\sigma_*\le \tau_0(u)\}}
u^0_*\big(X_{\sigma_*},M_{\sigma_*}\big)+\ind_{\{\sigma_*>\tau_0(u)\}}u^0_*\big(X_{\tau_0(u)},M_{\tau_0(u)}\big)\\
&\le \ind_{\{\sigma_*\le \tau_0(u)\}}g\big(X_{\sigma_*}\big)+ \ind_{\{\sigma_*>\tau_0(u)\}}f\big(X_{\tau_0(u)}\big).
\end{aligned}
\end{equation}	
Plugging the latter back into the expectation \eqref{eq:semim0} yields
\begin{equation*}
\begin{aligned}
u^0_*(x,m)&\le
\E^0_{x,m}\Big[\ind_{\{\sigma_*< \tau_0(u)\}}g\big(X_{\sigma_*}\big)+\ind_{\{\sigma_*= \tau_0(u)\}}g\big(X_{\sigma_*}\big)+
\ind_{\{\sigma_*>\tau_0(u)\}}f\big(X_{\tau_0(u)}\big)\Big]\\
&\le\E^0_{x,m}\Big[\ind_{\{\sigma_*< \tau_0(u)\}}g\big(X_{\sigma_*}\big)+\ind_{\{\sigma_*=
\tau_0(u)\}}h\big(X_{\sigma_*}\big)+ \ind_{\{\sigma_*>\tau_0(u)\}}f\big(X_{\tau_0(u)}\big)\Big],
\end{aligned}
\end{equation*}
where the second inequality holds because $g\le h$. Integrating over $u\in(0,1)$ we get
\begin{equation*}
\begin{aligned}
u^0_*(x,m)&\le\int_0^1\E^0_{x,m}\Big[\ind_{\{\sigma_*< \tau_0(u)\}}g\big(X_{\sigma_*}\big)+\ind_{\{\sigma_*=
\tau_0(u)\}}h\big(X_{\sigma_*}\big)+ \ind_{\{\sigma_*>\tau_0(u)\}}f\big(X_{\tau_0(u)}\big)\Big]\ud u\\
&=\E^0_{x,m}\Big[\ind_{\{\sigma_*< \tau_0\}}g\big(X_{\sigma_*}\big)+\ind_{\{\sigma_*= \tau_0\}}h\big(X_{\sigma_*}\big)+
\ind_{\{\sigma_*>\tau_0\}}f\big(X_{\tau_0}\big)\Big],
\end{aligned}
\end{equation*}
where in the equality we used $\tau_0=\tau_0(U)$ with $U\sim U(0,1)$ independent of $\F^2$. By the arbitrariness of
$\tau_0\in\cT_R(\F^2)$ we deduce
\begin{equation}\label{eq:ito2}
u^0_*(x,m)\le \inf_{\tau_0\in\cT_R(\F^2)}\E^0_{x,m}\Big[\ind_{\{\sigma_*< \tau_0\}}g\big(X_{\sigma_*}\big)+\ind_{\{\sigma_*=
\tau_0\}}h\big(X_{\sigma_*}\big)+ \ind_{\{\sigma_*>\tau_0\}}f\big(X_{\tau_0}\big)\Big].
\end{equation}
Choosing $\tau_0=\He$, the inequality in \eqref{eq:ito1} is replaced by
\begin{equation}
\begin{aligned}
u^0_*\big(X_{\He\wedge\sigma_*},M_{\He\wedge\sigma_*}\big)&=u^0_*\big(X_{\sigma_*},M_{\sigma_*}\big)=
g\big(X_{\sigma_*}\big),\quad \P^0_{x,m}-a.s.,
\end{aligned}
\end{equation}
where in the second equality we used that $\sigma_*<\He$, $\P^0_{x,m}$-a.s., because $d_*(m)>0$ for all $m\in(0,1]$ and $X$
does not hit one. For the same reasons
\begin{equation*}
g\big(X_{\sigma_*}\big)=g\big(X_{\sigma_*}\big)\ind_{\{\sigma_*<\He\}}+\ind_{\{\sigma_*= \He\}}h\big(X_{\sigma_*}\big)+
\ind_{\{\sigma_*>\He\}}f\big(X_{\He}\big),\quad \P^0_{x,m}-a.s.
\end{equation*}
Thus, repeating the remaining steps above we obtain
\begin{equation}\label{eq:ito3a}
u^0_*(x,m)= \E^0_{x,m}\Big[g\big(X_{\sigma_*}\big)\ind_{\{\sigma_*<\He\}}+\ind_{\{\sigma_*= \He\}}h\big(X_{\sigma_*}\big)+
\ind_{\{\sigma_*>\He\}}f\big(X_{\He}\big)\Big].
\end{equation}
Combining the above with \eqref{eq:ito2} we get \eqref{eq:ito3}.

Now we repeat analogous arguments for $u^1_*$. We start once again from the right-hand side of \eqref{eq:semim0}. Recall that $u^1_*(x,m)=g(x)$ for $x\le d_*(m)$. Since $X_{\sigma_*}\le d_*(M_{\sigma_*})$ on the event $\{\sigma_*<\He\}$,
we deduce
\begin{equation}\label{eq:ito5}
\begin{aligned}
&u^1_*\big(X_{\tau_1(u)\wedge\sigma_*},M_{\tau_1(u)\wedge\sigma_*}\big)\\
&=\ind_{\{\sigma_*\le \tau_1(u)\}}
u^1_*\big(X_{\sigma_*},M_{\sigma_*}\big)+\ind_{\{\sigma_*>\tau_1(u)\}}u^1_*\big(X_{\tau_1(u)},M_{\tau_1(u)}\big)\\
&=  \ind_{\{\sigma_* < \tau_1(u)\}\cap\{\sigma_*<\He\}}g\big(X_{\sigma_*}\big)
+ \ind_{\{\sigma_* = \tau_1(u)\}\cap\{\sigma_*<\He\}}g\big(X_{\sigma_*}\big)\\
&\quad+\ind_{\{\sigma_*= \tau_1(u)\}\cap\{\sigma_*=\He\}} u^1_*\big(X_{\sigma_*},M_{\sigma_*}\big)
+\ind_{\{\sigma_*>\tau_1(u)\}}u^1_*\big(X_{\tau_1(u)},M_{\tau_1(u)}\big).
\end{aligned}
\end{equation}
Since $g \leq h$ we know that $\ind_{\{\sigma_* = \tau_1(u)\}\cap\{\sigma_*<\He\}}g\big(X_{\sigma_*}\big) \leq \ind_{\{\sigma_* = \tau_1(u)\}\cap\{\sigma_*<\He\}}h\big(X_{\sigma_*}\big)$.
By the second equation in \eqref{eqn:u1f}, we know that
\[
\ind_{\{\sigma_*>\tau_1(u)\}}u^1_*\big(X_{\tau_1(u)},M_{\tau_1(u)}\big)\le\ind_{\{\sigma_*>\tau_1(u)\}}f(X_{\tau_1(u)}).
\]
Moreover, under $\P^1_{x,m}$ we have $\He=\H_1$
and therefore, using $u^1_*(1,1)=f(1)=h(1)$ (cf.\ \eqref{eq:PDEu1}) we obtain
\begin{equation}\label{eq:ito6}
\begin{aligned}
\ind_{\{\sigma_*= \tau_1(u)\}\cap\{\sigma_*=\He\}} u^1_*\big(X_{\sigma_*},M_{\sigma_*}\big) &=\ind_{\{\sigma_*=
\tau_1(u)\}\cap\{\sigma_*=\He\}} u^1_*(1,1) \\
&=\ind_{\{\sigma_*= \tau_1(u)\}\cap\{\sigma_*=\He\}} h\big(X_{\sigma_*}\big).
\end{aligned}
\end{equation}
Putting the terms together and recombining the two terms involving the function $h$ we arrive at our desired expression
\begin{equation}
u^1_*(x,m)\le \E^1_{x,m}\Big[\ind_{\{\sigma_*<
\tau_1(u)\}}g\big(X_{\sigma_*}\big)+\ind_{\{\sigma_*>\tau_1(u)\}}f(X_{\tau_1(u)})+\ind_{\{\sigma_*= \tau_1(u)\}}
h\big(X_{\sigma_*}\big)\Big].
\end{equation}
Integrating over $u\in(0,1)$ we finally get
\begin{equation}\label{eq:verif0}
u^1_*(x,m)\le \inf_{\tau_1\in\cT_R(\F^2)}\E^1_{x,m}\Big[\ind_{\{\sigma_*<
\tau_1\}}g\big(X_{\sigma_*}\big)+\ind_{\{\sigma_*>\tau_1\}}f(X_{\tau_1})+\ind_{\{\sigma_*= \tau_1\}}
h\big(X_{\sigma_*}\big)\Big].
\end{equation}

To saturate the inequality let us take $\tau_1(u)=\tau_\Gamma(u)=\inf\{t\ge 0: 1-\Gamma_t>u\}\wedge \He$ with the process $(\Gamma_t)$ as in the statement of the proposition. Note that since $\Gamma_{\rho_{c^1_*}}=1$ we must have $M_{\tau_1(u)} > c^1_*$, $\P^1_{x,m} $-a.s.  Also, as $\tau_1(u) \ge \tau_{\Gamma_*}(0)$, we have $M_{\tau_1(u)} \ge m^1_*$ and $u^1_*\big(M_{\tau_1(u)},M_{\tau_1(u)}\big)=f(M_{\tau_1(u)})$. By \eqref{eq:Gamma-t}, $X_{\tau_1(u)} = M_{\tau_1(u)}$ on $\{\tau_1(u) < \He\}$ for a.e.\ $u \in (0,1)$. Then $\{\sigma_* = \tau_1(u)\}\cap\{\sigma_*<\He\} = \varnothing$, because $d_*(m)<m$ for $m\in(c^1_*,1]$.
Hence, using \eqref{eq:semim0}, \eqref{eq:ito5}, \eqref{eq:ito6} and the third equation in \eqref{eq:PDEu1} we have for almost every $u \in (0,1)$,
\begin{equation*}
\begin{aligned}
u^1_*(x, m) &= \E^1_{x,m} \Big[u^1_*\big(X_{\tau_1(u)\wedge\sigma_*},M_{\tau_1(u)\wedge\sigma_*}\big)\Big]\\
&= \E^1_{x,m} \Big[\ind_{\{\sigma_* < \tau_1(u)\}\cap\{\sigma_*<\He\}}g\big(X_{\sigma_*}\big)  +\ind_{\{\sigma_*= \tau_1(u)\}\cap\{\sigma_*=\He\}} h\big(X_{\sigma_*}\big) +\ind_{\{\sigma_*>\tau_1(u)\}}f\big(X_{\tau_1(u)}\big) \Big]\\
&= \E^1_{x,m}[\cP(\sigma_*, \tau_1(u))].
\end{aligned}
\end{equation*}
Integration with respect to $u$ yields the optimality of $\tau_\Gamma$ in \eqref{eq:ito8}:
\[
\E^1_{x,m}[\cP(\sigma_*, \tau_\Gamma)] = \int_0^1 \E^1_{x,m}[\cP(\sigma_*, \tau_1(u))] \ud u = \int_0^1 u^1_*(x,m) \ud u = u^1_*(x,m).
\]
This concludes the proof.
\end{proof}

\subsection{A candidate PBE}\label{sec:candidatePBE}
We are now in a position to construct a candidate PBE.
Consistent with our construction from the earlier sections, in our candidate equilibrium $(\sigma_*,\tau_*)$ the type $\theta=0$ of Player 1 does not act. The type $\theta=1$ of Player 1 does two things: it stops immediately with some probability (possibly zero, to be determined) and, subsequently, acts so as to make Player 2's belief follow the dynamics
\begin{equation}\label{eq:p*}
p_t=\pi_*(X_t,M_t),\quad t\in[0,T_\e),
\end{equation}
where $\pi_*=\pi_{\Gamma_*}$ is from \eqref{eq:pi}, $M_t = M_0 \vee \sup_{s \le t} X_s$ and $M_0\ge X_0$ must be chosen endogenously at time 0 so as to appropriately match Player 2's belief. In order to perform such matching we start by noticing that for $m\in(0,m^1_*]$ we have $\pi_*(x,m)=x$, because $\Gamma_*(m)=1$, and for $m\in(m^1_*,1]$ we have $\pi_*(x,m)<x$, because
$\Gamma_*(m)<1$. In particular, $\partial_m \pi_*(x,m)=0$ for $m\in(0,m^1_*)$ and $\partial_m \pi_*(x,m)<0$ for $m\in(m^1_*,1)$. Therefore, from \eqref{eq:p*} we deduce that, on the set where the game does not end immediately,
\[
\pi_*(X_t,1)\le p_t=\pi_*(X_t,M_t)\le\pi_*(X_t,X_t) \le X_t,\quad t\in[0,T_\e),
\]
where the second inequality holds because $M_t\ge X_t$ by definition of $M_t$.

At the very beginning of the game Player 1 has performed no action and therefore Player 2's belief about the distribution of $\theta$ corresponds to $p = x$. Since $X_0=x\in(0,1)$ can be chosen arbitrarily and the Brownian dynamics visits all points with positive probability, an admissible set $A=A^{X,\tau_*}$ for the process $(X_t,p_t)$, as per Definition \ref{def:admissible}, reads as $A = A_1 \cup A_2$, where
\begin{equation}
\begin{aligned}
A_1 &= \{ (x, p)\in (0,1)^2 \colon p=x\ \text{and}\ \pi_*(x, x)<  x \},\\
A_2 &= \{ (x, p) \in (0,1)^2 \colon \pi_*(x, 1) < p \le \pi_*(x, x) \}.
\end{aligned}
\end{equation}
Player 2's belief at the very beginning of the game is $(x,p)=(x,x)$ and it belongs to either $A_1$ or $A_2$ depending on $x$. When $x \le m^1_*$, we have $\pi_*(x,x) = x$ and $(x,x) \in A_2$. However, when $x > m^1_*$, we have $\pi_*(x,x) < x$ and $(x,x) \in A_1$. We will see that this distinction is crucial in determining the appropriate value of $M_0$.

In order to construct a PBE, suppose that the game started at time zero from initial state $(X_0,p_{0-})=(x,x)\in A$, where $p_{0-}$ is the belief {\em before} any action is taken by either player. Suppose the game has not ended yet at time $t\ge 0$. Then, assuming that it is possible for Player 1 to make Player 2's belief evolve according to \eqref{eq:p*}, in the belief-augmented state-space the process is in state
\[
(X_t,p_t)=\big(X_t,\pi_*(X_t,M_t)\big)\in A_2,
\]
where the final inclusion holds because $\pi_*(X_t,M_t)\le \pi_*(X_t,X_t)$, as seen before. Thus, restarting the game at time $t$, we only need to be concerned with initial states $(x,p)\in A_2$. However, if at time zero $(X_0,p_{0-})=(x,x)\in A_1$ we will later see that Player 1 may stop immediately with positive probability.

First we consider a game that started at some time in the past and is currently being played starting from the current state $(x,p)\in A_2$. With no loss of generality we reset time so that time zero is the current time. Then, we want to determine $M_0$ from the identity
\[
p=p_{0-}=\pi_*(X_0,M_0)=\pi_*(x,M_0).
\]
From \eqref{eq:pi} we can solve for $\Gamma_*(M_0)$ and get,
\[
\Gamma_*(M_0) = \frac{p_{0-}(1 - X_0)}{X_0(1 - p_{0-})}\eqqcolon\Psi(X_0, p_{0-}).
\]
Note that $\Psi(x,p) = \frac{p/(1-p)}{x/(1-x)}$ is the odds ratio and for $(x,p) \in A_2$ we have $\Psi(x,p) \in [\Gamma_*(1), \Gamma_*(x)]$.
The inverse $\Gamma^{-1}_*$ is well defined on $[\Gamma_*(1), 1)$, due to the strict monotonicity of $\Gamma_*$ on $(m^1_*,1]$ (see Lemma \ref{lem:Gamma}), and it is extended by continuity to $[\Gamma_*(1), 1]$ with $\Gamma^{-1}_*(1) = m^1_*$.
Then we set
\begin{equation} \label{eq:M0}
M_0=\Gamma_*^{-1}\big(\Psi(X_0, p_{0-})\big)\ge X_0.
\end{equation}
Notice that the inequality holds because $p_{0-}=\pi_*(X_0,M_0)$ and $p_{0-}\le \pi_*(X_0,X_0)$ imply necessarily $M_0\ge X_0$ by monotonicity of $\pi_*(x,\cdot)$.

Next we provide the candidate optimal randomised stopping time for Player 1. Set  
\begin{equation}\label{eq:xiA2}
\xi_t^{1,p} = 
\begin{cases}
1 - \frac{\Gamma_*(M_t)}{\Gamma_*(M_0)}, & t < \He,\\
1, & t \ge \He.
\end{cases}
\end{equation}
and note that $\xi^{1,p}$ is continuous on $[0, \He)$, non-decreasing and satisfies $\xi^{1,p}_0=0$; further, $\xi^{1,p}$ depends on $p$ via $M_0$. Then
our candidate is
\begin{equation}
\label{eq:taudef6}
\tau^{*,p} = \tau^* = \ind_{\{\theta=0\}} \He + \ind_{\{\theta=1\}} \tau^{*,p}_1,
\end{equation}
where
$\tau^{*,p}_1 = \inf\{t \ge 0\colon  \xi^{1,p}_t > U\}$
and $U \sim U(0,1)$ is independent of $\cF^1_\infty$. Notice that $\tau^{*,p}_1\le \He$.

From \eqref{eq:hatpX} and the formula for the odds ratio $\Psi$ we have that
\[ 
\hat{p}_t = \frac{X_t \Psi(x,p)}{X_t \Psi(x,p) + (1- X_t)}. 
\]
Then from \eqref{eqn:pt}, with $\xi^0_t\equiv 0$ and $\xi^1_t=\xi^{1,p}_t$, we find using $\Gamma_*(M_0) = \Psi(x,p)$
\[ 
p_t = \frac{\Gamma_*(M_t) X_t}{\Gamma_*(M_t) X_t + (1-X_t)} = \pi_*(X_t,M_t). 
\]
Hence \eqref{eq:p*} holds for our choice of $\tau^{*,p}$ as given in \eqref{eq:taudef6}. Since $(X_t,p_t)\in A_2$, there is a one-to-one correspondence between the dynamics $(X_t,p_t)$ and $(X_t,M_t)$ with
\[
M_t=\Gamma_*^{-1}\big(\Psi(X_t, p_t)\big),\quad t\in[0,\He).
\]
Thus, we have shown that if $(x,p)\in A_2$ we can determine endogenously $M_0$ and construct a family of generating processes $(\xi^{1,p}_t)$ which is {\em continuous} on $[0,\He)$ and such that the uninformed player's belief evolves with dynamics given by \eqref{eq:p*}. In what follows we write $\xi^p$ as a shorthand for $\xi^{1,p}$.

Next we consider the case when at the very beginning of the game $(X_0,p_{0-})=(x,x)\in A_1$. If we just try to repeat arguments as above we run into trouble, because $\Gamma_*(M_0)=\Psi(x,x)=1$ and $M_0=\Gamma_*^{-1}(1)=m^1_*$. However, $(x,x)\in A_1\implies x>m^1_*$ and therefore $M_0<X_0$, violating the requirement that $M_0\ge X_0$. Intuitively, the fact that $p_{0-}>\pi_*(X_0,X_0)$ means that Player 2's initial belief is {\em too high}. So, allowing Player 1 to stop immediately with some probability will lower the belief to $p_0=\pi_*(X_0,X_0)<p_{0-}$. From that moment onward Player 1 can guarantee that the belief evolves as in \eqref{eq:p*}.

More precisely, recalling that the initial point is $(x,x) \in A_1$, we select a generating process of the form
\begin{equation}\label{eqn:xi-with-jump}
\xi^x_{0-} = 0, \quad 
\xi^x_t = \begin{cases}
\Delta\xi^x_0 + (1-\Delta \xi^x_0) \left(1 - \frac{\Gamma_*(M_t)}{\Gamma_*(M_0)}\right), & t < \He,\\
1, & t \ge \He,
\end{cases}
\end{equation}
where we must determine both $\Delta\xi^x_0$ and $M_0$ endogenously. Setting $\tau^{*,x}$ as in \eqref{eq:taudef6} with $\xi^x$ as in \eqref{eqn:xi-with-jump} we can compute (cf.\ \eqref{eqn:pt})
\[
p_0=\hP_{x}^{x}\big(\theta=1|\cF^2_0,\tau^*>0\big)=\frac{p_{0-}(1-\Delta\xi^x_0)}{(1-p_{0-}) + p_{0-}(1-\Delta\xi^x_0)}=\frac{x(1-\Delta\xi^x_0)}{(1-x) + x(1-\Delta\xi^x_0)}.
\]
Setting
\begin{equation}\label{eqn:qq}
p_0=\pi_*(x,x)=\frac{x \Gamma_*(x)}{1-x + x \Gamma_*(x)}\eqqcolon q(x)=q
\end{equation}
we get
\begin{equation}
\label{eqn:Deltaxi0}
\Delta\xi^x_0=\Delta\xi^x_0(x) = 1 - \frac{q}{1-q} \frac{1-x}{x} = 1-\Gamma_*(x).
\end{equation}
As $(x, q) \in A_2$, we have $\Psi(x,q)=\Psi(x,\pi_*(x,x))=\Gamma_*(x)$. Then, the running maximum process $(M_t)$ is initialised at $M_0 = X_0 = x$ and the dynamics of $(\xi^x_t)$ in \eqref{eqn:xi-with-jump} is fully determined. We notice that, once again combining the formula for $\hat p_t$ from \eqref{eq:hatpX} with the one for $p_t$ from \eqref{eqn:pt}, with $p=x$, $\xi^0_t\equiv 0$ and $\xi^1_t$ as in \eqref{eqn:xi-with-jump}, we can verify that \eqref{eq:p*} holds. For completeness, whenever $(x,p) \notin A = A_1 \cup A_2$, we set $\xi^{1,p}_t = 0$ for all $t \le \inf\{s \ge 0\colon (X_s, \hat p_s) \in A\}$. These values are not on an equilibrium path, so their choice does not affect the game.

So far we have focussed on actions of the informed player. To complete the construction of our candidate PBE we now consider the stopping rule of Player 2. Given the pair $(X_t,p_t)$ with belief process $(p_t)$ constructed as above,
the ex-ante expected payoff reads as
\[
\hE_{x}^{p}[\mathcal{P}(\sigma, \tau^*)] = p \E_x^1[\mathcal{P}(\sigma, \tau^*)] + (1 - p) \E_x^0[\mathcal{P}(\sigma, \tau^*)],
\]
with $\tau^*$ as in \eqref{eq:taudef6}. We expect that Player 2 maximises such payoff by stopping at $\sigma^{*,p}$ given by 
\begin{equation}
\label{eq:sigmadef6}
\begin{aligned}
\sigma^{*,p} &= \inf\{t \ge 0\colon X_t \le d_*(\Gamma^{-1}_*(\Psi(X_t, p_t))) \}\wedge\He= \inf\{t \ge 0\colon X_t \le d_*(M_t) \}\wedge\He.
\end{aligned}
\end{equation}
Notice that, in particular $\hP^p_x(\sigma^{*,p}>0)=1$ for $(x,p)\in A_1$ because $X_0=M_0=x>m^1_*$ and $d_*(x)=c_*(x)<x$ for $x>m^1_*$, by construction of $c_*$. This means that Player 2 never stops at time zero with initial state $(X_0,p_{0-})\in A_1$. We also notice that $\sigma^{*,p}$ is effectively the same as $\sigma_*$ from Lemma \ref{lem:semim} but we use both notations to emphasise a small subtlety: $\sigma^{*,p}$ is defined in terms of the pair $(X_t,p_t)$ and the process $(M_t)$ is constructed as described in this section; instead, $\sigma_*$ is defined for the pair $(X_t,M_t)$ with arbitrary $M_0\ge X_0$. The construction of $\sigma^{*,p}$ makes it obvious that it is a path functional of $X$ dependent on $p$ via $p_t$, so it is of the form required in the definition of state-belief Markovian strategies.

We recall that we omit the index $p$ in $\tau^*$, $\xi$ and $\sigma^*$ whenever its value is identical to the subscript of the probability measure $\hP^p_x$.

\subsection{Verification of the PBE}
Throughout the section we refer to the pair $(\sigma^*,\tau^*)$ from \eqref{eq:taudef6} and \eqref{eq:sigmadef6}.
We start the analysis with $(x,p) \in A_2$ as it requires a number of preparatory remarks. The Radon-Nikodym derivative of
$\hP^p_x$ with respect to $\P_x$ is given by (cf. Lemma \ref{lem:Y-Pyp})
\begin{equation}\label{eqn:RN_xp}
\frac{\ud \hP^p_x}{\ud \P_{x}}
= \ind_{\{\theta=0\}} \frac{1-p}{1-x} + \ind_{\{\theta=1\}} \frac{p}{x}
= \frac{\ind_{\{\theta=0\}} + \Gamma_*(M_0) \ind_{\{\theta=1\}}}{1 - x + x \Gamma_*(M_0)},
\end{equation}
where the second equality follows from the definition of $M_0 = \Gamma^{-1}_*(\Psi(x,p))$ and direct calculations.
The restriction of the Radon-Nikodym derivative to $\cF^2_t$ defines a process
\begin{equation}\label{eqn:RNprocess}
\frac{\ud \hP^p_x}{\ud \P_{x}} \Big|_{\cF^2_t} =
\E_{x} \bigg[ \frac{\ind_{\{\theta=0\}} + \Gamma_*(M_0) \ind_{\{\theta=1\}}}{1 - x + x \Gamma_*(M_0)} \bigg| \cF^2_{t} \bigg] = \frac{1-X_{t\wedge\He} + X_{t\wedge\He} \Gamma_*(M_0)}{(1-x) + x \Gamma_*(M_0)}.
\end{equation}
We will use \eqref{eqn:RNprocess} to rewrite the expectation $\hE^p_x$ in terms of $\E_{x}$ so that we can apply results of preceding sections written in terms of
the pair $(X_t,M_t)$ with $M_t = M_0 \vee \sup_{s \le t} X_s$ constructed as in Section \ref{sec:candidatePBE}.

Consider the process for $t\le \He$
\begin{equation}\label{eqn:Lambda}
\Lambda_t = \hE^p_x \left[ \ind_{\{\tau^* \ge t\}} \left( \ind_{\{\theta=1\}} u^1_* + \ind_{\{\theta=0\}} u^0_* \right)(X_{t}, M_{t}) +
\ind_{\{\theta=1\}} \int_{[0,{t})}f(X_s) \ud\xi_s \bigg| \cF^2_t\right],
\end{equation}
which will be used in several instances below.
\begin{lemma}\label{lem:6.5}
Let $(x, p) \in A_2$ and recall $v_*$ as in \eqref{eq:v*def} for $0\le x\le m\le 1$. We have $\hP^p_x$\as on $\{t \le \He\}$
\begin{equation}\label{eqn:Lambdat}
\Lambda_t = \frac{(1-X_t) + X_t \Gamma_*(M_t)}{ (1-X_{t}) + X_{t} \Gamma_*(M_0) } v_*(X_t, M_t)
- \frac{X_t \int_{[0,t)} f(X_u) \Gamma'_*(M_u) \ud M_u } {(1-X_t) + X_t \Gamma_*(M_0) }.
\end{equation}
\end{lemma}
\begin{proof}
Using \eqref{eqn:RN_xp} and \eqref{eqn:RNprocess}, we have
\begin{equation}\label{eqn:u01}
\begin{aligned}
&\hE^p_x \big[ \ind_{\{\tau^* \ge t\}} \big( \ind_{\{\theta=1\}} u^1_* + \ind_{\{\theta=0\}} u^0_*)(X_{t}, M_{t}) \big|
\cF^2_{t} \big]\\
&=
\E_{x} \Big[ \frac{\ud \hP^p_x}{\ud \P_x} \ind_{\{\tau^* \ge t\}} \big(
\ind_{\{\theta=1\}} u^1_* + \ind_{\{\theta=0\}} u^0_*)(X_{t}, M_{t}) \Big| \cF^2_{t} \Big]
\Big/{\frac{\ud \hP^p_x}{\ud \P_x}\Big|_{\cF^2_t}}.
\end{aligned}
\end{equation}
Observing that
\begin{equation}\label{eq:rationRN}
\frac{\ud \hP^p_x}{\ud \P_x}\Big/\left(\frac{\ud \hP^p_x}{\ud \P_x}\Big|_{\cF^2_t}\right)=\frac{\ind_{\{\theta=0\}} + \Gamma_*(M_0) \ind_{\{\theta=1\}}}{1-X_t + \Gamma_*(M_0) X_t},
\end{equation}
the expression in \eqref{eqn:u01} simplifies to
\begin{equation}
\begin{aligned}
&\E_{x} \Big[ \frac{\ind_{\{\theta=0\}} + \Gamma_*(M_0) \ind_{\{\theta=1\}}}{{1 - X_t + X_t \Gamma_*(M_0)}} \ind_{\{\tau^* \ge t\}} \big(
\ind_{\{\theta=1\}} u^1_* + \ind_{\{\theta=0\}} u^0_*)(X_{t}, M_{t}) \Big| \cF^2_{t} \Big]\\
&=\frac1{1 - X_t + X_t \Gamma_*(M_0)} \left(\E_{x} \big[\Gamma_*(M_0) \ind_{\{\theta=1\}} \ind_{\{\tau^*_1 \ge t\}}u^1_*(X_{t}, M_{t})
\big| \cF^2_{t} \big]
+ \E_{x} \big[\ind_{\{\theta=0\}} u^0_*(X_{t}, M_{t}) \big| \cF^2_{t} \big] \right)\\
&=\frac1{1 - X_t + X_t \Gamma_*(M_0)} \left( \E_{x} \big[\Gamma_*(M_0) \ind_{\{\theta=1\}} (1-\xi_{t-}) u^1_*(X_{t}, M_{t}) \big|
\cF^2_{t} \big]
+ \E_{x} \big[\ind_{\{\theta=0\}} u^0_*(X_{t}, M_{t}) \big| \cF^2_{t} \big] \right)\\
&=\frac1{1 - X_t + X_t \Gamma_*(M_0)} \left( \Gamma_*(M_0) X_{t} \frac{\Gamma_*(M_{t})}{\Gamma_*(M_0)} u^1_*(X_{t}, M_{t})
+ (1-X_{t}) u^0_*(X_{t}, M_{t}) \right)\\
&=\frac{ X_{t} \Gamma_*(M_{t}) u^1_*(X_{t}, M_{t}) + (1-X_{t}) u^0_*(X_{t}, M_{t}) }{1 - X_t + X_t \Gamma_*(M_0)},
\end{aligned}
\end{equation}
where in the first equality we used that $\tau^* = \He$ on $\theta=0$, and in the second equality we integrated out the randomisation
device of Player 1. We justify the third equality separately on $\{t < \He\}$ and on $\{t = \He\}$. On the set $\{t < \He\}$, we
apply the definition of $\xi_t$ and use the continuity of $t \mapsto \Gamma_*(M_t)$. On the set $\{ t = \He \}$, we have
\[
1 - \xi_{\He-} = \frac{\Gamma_*(M_{\He-})}{\Gamma_*(M_0)} = \frac{\Gamma_*(M_{\He})}{\Gamma_*(M_0)},
\]
due to the continuity of $\Gamma_*$ on $[0,1]$. This shows that the above formula extends to $t \le \He$.
We have obtained
\begin{equation}\label{eqn:v_bayes}
\hE^p_x \big[ \ind_{\{\tau^* \ge t\}} \big( \ind_{\{\theta=1\}} u^1_* + \ind_{\{\theta=0\}} u^0_*)(X_{t}, M_{t}) \big| \cF^2_{t} \big]
=
\frac{ X_{t} \Gamma_*(M_{t}) u^1_*(X_{t}, M_{t}) + (1-X_{t}) u^0_*(X_{t}, M_{t}) }{ X_{t} \Gamma_*(M_0) + (1-X_{t}) },
\end{equation}
which in turn is equal to
\[
\frac{1-X_t + X_t \Gamma_*(M_t)}{ 1-X_{t} + X_{t} \Gamma_*(M_0) }v_*(X_t, M_t),
\]
thanks to \eqref{eq:v*def}.

Using the definition of $(\xi_t)$ from \eqref{eq:xiA2} and \eqref{eq:rationRN}, we obtain on $\{t < \He\}$
\begin{equation*}
\begin{aligned}
\hE^p_x \bigg[\ind_{\{\theta=1\}} \int_{[0,t)} f(X_u) \ud \xi_u \bigg| \cF^2_t\bigg] &=
\E_{x} \Big[ \frac{\ind_{\{\theta=0\}} + \Gamma_*(M_0) \ind_{\{\theta=1\}}}{1 - X_t + X_t \Gamma_*(M_0)} \int_{[0,t)}
f(X_u) \ind_{\{\theta=1\}} \ud \xi_u \Big| \cF^2_t\Big]\\
&= -
\frac{\E_{x} \Big[ \Gamma_*(M_0) \ind_{\{\theta=1\}} \int_{[0,t)} f(X_u) \Gamma'_*(M_u) / \Gamma_*(M_0) \ud M_u \Big|
\cF^2_t\Big] } {1-X_t + \Gamma_*(M_0) X_t}\\
&= - \frac{\E_{x}[\ind_{\{\theta=1\}} | \cF^2_t] \int_{[0,t)} f(X_u)  \Gamma'_*(M_u) \ud M_u} {1-X_t + \Gamma_*(M_0) X_t}\\
&= - \frac{X_t \int_{[0,t)} f(X_u){\Gamma'_*}(M_u) \ud M_u } {1-X_t + \Gamma'_*(M_0) X_t},
\end{aligned}
\end{equation*}
where in the third equality we used $\cF^2_t$-measurability of the integral and in the final one we used $\P_x(\theta=1|\cF^2_t)=X_t$ because $(X_t)$ is a Brownian motion under $\P_x$. Combining the above calculations yields \eqref{eqn:Lambdat}.
\end{proof}

\begin{remark}\label{rem:taugt}
It follows from the computations in the above lemma that on $\{t\le\He\}$
\begin{equation}
\begin{aligned}
\hE^p_x \big[ \ind_{\{\tau^* \ge t\}} \big| \cF^2_t \big] &=\E_x\Big[\frac{\ind_{\{\theta=0\}}+\Gamma_*(M_0)\ind_{\{\theta=1\}}}{1-X_t+X_t\Gamma_*(M_0)}\ind_{\{\tau^* \ge t\}}\Big|\cF^2_t\Big]\\
&=\E_x\Big[\frac{\ind_{\{\theta=0\}}\ind_{\{\He\ge t\}}+\Gamma_*(M_0)\ind_{\{\theta=1\}}\ind_{\{\tau^*_1\ge t\}}}{1-X_t+X_t\Gamma_*(M_0)}\Big|\cF^2_t\Big]\\
&= \frac{1-X_t + X_t \Gamma_*(M_t)}{1-X_{t} + X_{t} \Gamma_*(M_0)},
\end{aligned}
\end{equation}
which provides an interpretation of the multiplicative factor in front of $v_*$ as the probability perceived by Player 2 that strictly prior to time $t$ Player 1 has not ended the game yet.
\end{remark}

Set
\[
\hat v_*(x,m) = (1 - x + x \Gamma_*(m))v_*(x,m) = x \Gamma_*(m) u^1_*(x,m) + (1-x)u^0_*(x,m).
\]
In the next lemma we calculate conditional expectation of increments of the process $(\Lambda_t)$.
\begin{lemma}\label{lem:mart-Lambda}
For $(x,p) \in A_2$ and $t \ge s \ge 0$, we have on $\{s < \He\}$
\begin{equation}\label{eqn:Lambda_m}
\begin{aligned}
&\hE^p_x \big[ \Lambda_{t \wedge \He} \!-\! \Lambda_s  \big| \cF^2_s \big]\\
&=
\frac1{1\!-\!X_s\! +\! \Gamma_*(M_0) X_s} \E_{x} \Big[ \hat v_*(X_{t \wedge \He}, M_{t \wedge \He})\! -\! \hat v_*(X_{s}, M_{s})
-\!\! \int_{[s,t \wedge \He)}\!\!  f(X_u) X_u \Gamma'_*(M_u) \ud M_u \Big| \cF^2_s\Big].
\end{aligned}
\end{equation}
\end{lemma}
\begin{proof}
Fix $t \ge s \ge 0$ and, on the event $\{s<\He\}$, consider the difference
\begin{align*}
&\hE^p_x \big[ \Lambda_{t \wedge \He} - \Lambda_s  \big| \cF^2_s \big]\\
&=\hE^p_x \bigg[\frac{1}{1-X_{t \wedge \He} + X_{t \wedge \He} \Gamma_*(M_0)} \hat v_*(X_{t \wedge \He}, M_{t \wedge \He}) \bigg| \cF^2_s \bigg]
- \frac{1}{1-X_{s} + X_{s} \Gamma_*(M_0)} \hat v_*(X_s, M_s)\\
&\quad - \hE^p_x\bigg[\frac{X_{t\wedge\He}\int_{[0,t \wedge \He)} f(X_u) \Gamma'_*(M_u) \ud M_u } {1-X_{t \wedge \He} + X_{t \wedge \He}\Gamma_*(M_0) } \bigg| \cF^2_s \bigg]
+ \frac{X_s\int_{[0,s)} f(X_u) {\Gamma'}_*(M_u) \ud M_u } {1-X_s + X_s\Gamma_*(M_0) },
\end{align*}
where we used \eqref{eqn:Lambdat} and the definition of $\hat v$.
We change the measure to $\P_x$ using \eqref{eq:rationRN} and the fact that for any non-negative $\cF^2_{t\wedge\He}$-measurable random variable $\Xi$
\begin{equation}
\begin{aligned}
\hE^p_x\big[\Xi\big|\cF^2_s\big]&=\E_x\bigg[\frac{\ind_{\{\theta=0\}}+\Gamma_*(M_0)\ind_{\{\theta=1\}}}{1-X_s+X_s\Gamma_*(M_0)}\Xi\bigg|\cF^2_s\bigg]
=\frac{\E_x\Big[\E_x\big[\ind_{\{\theta=0\}}+\Gamma_*(M_0)\ind_{\{\theta=1\}}\big|\cF^2_{t\wedge\He}\big]\Xi\Big|\cF^2_s\Big]}{1-X_s+X_s\Gamma_*(M_0)}\\
&=\E_x\bigg[\frac{1-X_{t\wedge\He}+X_{t\wedge\He}\Gamma_*(M_0)}{1-X_s+X_s\Gamma_*(M_0)}\Xi\bigg|\cF^2_s\bigg].
\end{aligned}
\end{equation}
That yields
\begin{equation}\label{eq:t00}
\begin{aligned}
&\hE^p_x \big[ \Lambda_{t \wedge \He} - \Lambda_s  \big| \cF^2_s \big]\\
&=\E_{x} \bigg[\frac{1}{1-X_{s} + X_{s} \Gamma_*(M_0)} \hat v_*(X_{t \wedge \He}, M_{t \wedge \He}) \bigg| \cF^2_s \bigg]
- \frac{1}{1-X_{s} + X_{s} \Gamma_*(M_0)} \hat v_*(X_s, M_s)\\
&\quad - \E_{x}\bigg[\frac{X_{t \wedge \He} \int_{[0,t \wedge \He)} f(X_u) \Gamma'_*(M_u) \ud M_u } {1-X_s + X_s\Gamma_*(M_0)} \bigg| \cF^2_s \bigg]
+ \frac{X_s \int_{[0,s)} f(X_u) \Gamma'_*(M_u) \ud M_u } {1-X_s + X_s\Gamma_*(M_0)}.
\end{aligned}
\end{equation}
By the product rule, we have
\begin{equation*}
\begin{aligned}
X_{t \wedge \He} \int_{[0,t \wedge \He)} f(X_u) \Gamma'_*(M_u) \ud M_u 
&=X_s \int_{[0,s)} f(X_u) \Gamma'_*(M_u) \ud M_u+\int_{[s,t \wedge \He)}X_u f(X_u) \Gamma'_*(M_u) \ud M_u\\
&\quad+\int_{[s,t \wedge \He)}\Big(\int_{[0,r)} f(X_u) \Gamma'_*(M_u) \ud M_u\Big)\ud X_r,
\end{aligned}
\end{equation*}
where the final integral is a $\P_x$-martingale. Plugging the expression above into \eqref{eq:t00} 
we obtain \eqref{eqn:Lambda_m}.
\end{proof}

\begin{lemma}\label{lem:belief-state-Markovian}
The pair $\big((\sigma^{*,p})_{p \in (0,1)}, (\tau^{*,p})_{p \in (0,1)}\big)$ constructed as in \eqref{eq:taudef6} and \eqref{eq:sigmadef6} with $\xi^p$ as in \eqref{eq:xiA2} or \eqref{eqn:xi-with-jump} is belief-state-Markovian in the sense of Definition \ref{def:rst-markov}.
\end{lemma}
\begin{proof}
By the explicit construction, the randomised stopping times $(\sigma^*, \tau^*)$ are determined by belief-dependent generating processes. Condition \eqref{eqn:unfinished} is also satisfied: the generating process for $\tau^{0,p}$ is given by $\xi^{0,p}_t = \ind_{\{t \ge \He\}}$ and the generating process $(\xi^p_t)$ for $\tau^{1,p}$ satisfies $\xi^{p}_t < 1$ whenever $M_t < 1$, i.e., for all $t < \He$, because $\Gamma_*(m) > 0$ for $m \in (0,1)$ by Lemma \ref{lem:Gamma}. The identity \eqref{eqn:rst-markov-xi} easily holds for $\xi^{0,p}$. We will verify it in a more complex case of the generating process $(\xi^p_t)$ of $\tau^*_1$ defined in \eqref{eq:xiA2}. Explicit indication of the superscript $p$ will be important for this proof.

We start with an auxiliary identity that we will use several times in this proof. For $(x,p) \in A_2$, stopping times $\eta, \eta' \in \cT(\F^2)$, $\eta < \He$,  $\eta \le \eta'$ and an $\cF^2_{\eta'}$-measurable bounded random variable $\Xi$, we have
\begin{equation}\label{eqn:aux_ind}
\begin{aligned}
\hE^p_x \Big[ \Xi \ind_{\{\tau^{*,p} > \eta\}} \Big| \cF^2_\eta \Big]
&=
\E_{x} \Big[ \Xi \frac{\ind_{\{\theta=0\}} + \ind_{\{\theta = 1\}} \Gamma_*(M_0)}{1 - X_\eta + X_\eta \Gamma_*(M_0)} \big( \ind_{\{\theta=0\}} + \ind_{\{\theta=1\}} (1-\xi^p_\eta) \big) \Big| \cF^2_\eta \Big]\\
&=
\E_{x} \Big[ \Xi \frac{\ind_{\{\theta=0\}} + \ind_{\{\theta = 1\}} \Gamma_*(M_\eta)}{1 - X_\eta + X_\eta \Gamma_*(M_0)} \Big| \cF^2_\eta \Big]
=
\E_{x} \Big[ \Xi \frac{1 - X_{\eta'} + X_{\eta'} \Gamma_*(M_\eta)}{1 - X_\eta + X_\eta \Gamma_*(M_0)}  \Big| \cF^2_\eta \Big],
\end{aligned}
\end{equation}
where $M_0 = \Gamma^{-1}_*(\Psi(x,p))$ and we used the definition of $(\xi^p_t)$ in the penultimate equality and the tower property with $\cF^2_{\eta'}$ in the last equality. For $\Xi=1$ the expression further simplifies to 
\[
\frac{1-X_\eta+X_\eta\Gamma_*(M_\eta)}{1-X_\eta+X_\eta\Gamma_*(M_0)}.
\]
Using the latter we arrive at
\begin{equation}\label{eqn:aux_ind_full}
\hE^p_x [ \Xi | \cF^2_\eta, \tau^{*,p} > \eta ]
=
\frac{\hE^p_x \big[ \Xi \ind_{\{\tau^{*,p} > \eta\}} \big| \cF^2_\eta \big]}{\hE^p_x [ \ind_{\{\tau^{*,p} > \eta\}} | \cF^2_\eta ] }
=
\E_{x} \bigg[ \frac{1 - X_{\eta'} + X_{\eta'} \Gamma_*(M_\eta)}{1 - X_\eta + X_\eta \Gamma_*(M_\eta)} \Xi \bigg| \cF^2_\eta \bigg].
\end{equation}

We return to the main line of the proof. We will treat separately $(x,p)$ in $A_1$ and in $A_2$. First fix $(x,p) \in A_2$ and $\eta \in \cT(\F^2)$, $\eta < \He$. 
Thanks to \eqref{eqn:aux_ind_full} with $\eta' = (\eta+s)\wedge \He$, we have
\begin{equation}\label{eqn:eq for A2}
\begin{aligned}
&\hE^p_x \bigg[ \varphi\bigg(\frac{1-\xi^p_{\eta + s}}{1-\xi^p_\eta} \bigg) \ind_{\{\eta+s<\He\}}\bigg| \cF^2_\eta, \tau^{*,p} > \eta \bigg]\\
&=
\E_{x} \bigg[ \frac{1 - X_{\eta+s} + X_{\eta+s} \Gamma_*(M_\eta)}{1 - X_\eta + X_\eta \Gamma_*(M_\eta)} \varphi\bigg(\frac{\Gamma_*(M_{\eta + s})}{\Gamma_*(M_\eta)} \bigg) \ind_{\{\eta+s<\He\}} \bigg| \cF^2_\eta \bigg].
\end{aligned}
\end{equation}
Denoting $m = M_0 =\Gamma_*^{-1}(\Psi(x,p))$, recall that, $M_t = m \vee \sup_{s \le t} X_s$. The pair $(X_t, m \vee \sup_{s \le t} X_s)$ under $\P_x$ describes the same stochastic process as $(X_t, M_t)$ under $\P_{x,m}$, i.e., the two processes are equal in distribution on the path space. 
Hence
\begin{equation}
\begin{aligned}
&\E_{x} \bigg[ \frac{1 - X_{\eta+s} + X_{\eta+s} \Gamma_*(M_\eta)}{1 - X_\eta + X_\eta \Gamma_*(M_\eta)} \varphi\bigg(\frac{\Gamma_*(M_{\eta+s})}{\Gamma_*(M_\eta)} \bigg) \ind_{\{\eta+s<\He\}} \bigg| \cF^2_\eta \bigg]\\
&=\E_{x,m} \bigg[ \frac{1 - X_{\eta+s} + X_{\eta+s} \Gamma_*(M_\eta)}{1 - X_\eta + X_\eta \Gamma_*(M_\eta)} \varphi\bigg(\frac{\Gamma_*(M_{\eta + s})}{\Gamma_*(M_\eta)} \bigg) \ind_{\{\eta+s<\He\}} \bigg| \cF^2_\eta \bigg]=Q_s(X_\eta,M_\eta),
\end{aligned}
\end{equation}
where in the last equality we use the strong Markov property of the pair $(X_t, M_t)$ with respect to the filtration $\F^2$ and
\begin{equation}\label{eq:Qs}
Q_s(z,\ell)\coloneqq\E_{z,\ell}\bigg[\frac{1 - X_{s} + X_{s} \Gamma_*(\ell)}{1 - z + z \Gamma_*(\ell)} \varphi\bigg(\frac{\Gamma_*(M_{s})}{\Gamma_*(\ell)} \bigg) \ind_{\{s < \He\}} \bigg].
\end{equation}
Since $(x,p)\in A_2$, we have $(X_t,p_t)=(X_t,\pi_*(X_t,M_t))$ for $t\in[0,\He)$, and moreover $M_t=\Gamma_*^{-1}(\Psi(X_t,p_t))$. In particular, $(X_\eta,p_\eta)\in A_2$ and in \eqref{eq:Qs} we can consider $(z,\pi_*(z,\ell))\in A_2$. Then, for $s < \He$, the ratio $\Gamma_*(M_{s})/\Gamma_*(\ell)$ under the measure $\P_{z,\ell}$ is equal to $1-\xi^{q}_{s}$ defined in \eqref{eq:xiA2} with $q=\pi_*(z,\ell)$. In other words
\begin{equation}\label{eq:Qs2}
\begin{aligned}
Q_s(z,\ell)&=\E_{z,\ell}\bigg[\frac{1 - X_{s} + X_{s} \Gamma_*(\ell)}{1 - z + z \Gamma_*(\ell)} \varphi\big(1-\xi^q_{s}\big) \ind_{\{s < \He\}}\bigg]\\
&=\E_{z,\ell}\bigg[\frac{\ind_{\{\theta=0\}} + \ind_{\{\theta=1\}} \Gamma_*(\ell)}{1 - z + z \Gamma_*(\ell)} \varphi\big(1-\xi^q_{s}\big) \ind_{\{s < \He\}}\bigg]
=\hE_{z}^{q}\big[ \varphi\big(1-\xi^q_{s}\big)\ind_{\{s < \He\}}\big],
\end{aligned}
\end{equation}
where again we used the change of measure \eqref{eqn:RN_xp}.

In conclusion, we have shown for $t<\He$ 
\begin{equation}\label{eqn:eta_s_less_He}
\hE^p_x \bigg[ \varphi\bigg(\frac{1-\xi^p_{\eta + s}}{1-\xi^p_\eta} \bigg) \ind_{\{\eta + s < \He\}} \bigg| \cF^2_\eta, \tau^{*,p} > \eta \bigg]=Q_s(X_\eta,M_\eta)=\hE_{X_\eta}^{p_\eta}\big[ \varphi\big(1-\xi^{p_\eta}_{s}\big) \ind_{\{s < \He\}}\big],
\end{equation}
where in the last equality we used \eqref{eq:Qs2} with $p_\eta=\pi_*(X_\eta,M_\eta)$. On the event $\{\eta+s \ge \He\}$, we have $\xi^p_{\eta+s} = 1$, so
\begin{equation}\label{eqn:eta_s_He}
\begin{aligned}
&\hE^p_x \bigg[ \varphi\bigg(\frac{1-\xi^p_{\eta + s}}{1-\xi^p_\eta} \bigg) \ind_{\{\eta+s \ge \He\}} \bigg| \cF^2_\eta, \tau^{*,p} > \eta \bigg]\\
&=
\varphi(0)\, \hP^p_x (\eta + s \ge \He | \cF^2_\eta, \tau^{*,p} > \eta)
=
\varphi(0)\, \hP^{p_\eta}_{X_\eta} (s \ge \He) =
\hE_{X_\eta}^{p_\eta}\big[ \varphi\big(1-\xi^{p_\eta}_{s}\big)\ind_{\{s \ge \He\}}\big],
\end{aligned}
\end{equation}
where the second equality follows from \eqref{eqn:eta_s_less_He} with $\varphi \equiv 1$. Combining the above two identities proves \eqref{eqn:rst-markov-xi}.

We turn our attention to $(x,p) \in A_1$ (hence $p=x$). Recall that the process $(X_t, p_t)$ is immediately shifted to $A_2$. In particular, the generating process $(\xi^x_t)$ reads as in \eqref{eqn:xi-with-jump}, so $1-\xi^x_{s} = (1 - \Delta \xi^x_0)(1-\xi^q_s)$, $s \ge 0$, with $q$ given in \eqref{eqn:qq}. Noting that $\hP_{x}^x = \P_{x,x}$, we have
\begin{equation}
\begin{aligned}
&\hE_{x}^x \bigg[ \varphi\bigg(\frac{1-\xi^x_{\eta + s}}{1-\xi^x_\eta} \bigg) \ind_{\{\eta+s<\He\}} \bigg| \cF^2_\eta,\, \tau^{*,x} > \eta\bigg]\\
&=
\E_{x,x} \bigg[ \varphi\bigg(\frac{1-\xi^x_{\eta + s}}{1-\xi^x_\eta} \bigg) \ind_{\{\eta+s<\He\}} \bigg| \cF^2_\eta,\, \tau^{*,x} > \eta\bigg]
=
\frac{\E_{x,x} \Big[ \varphi\Big(\frac{1-\xi^x_{\eta + s}}{1-\xi^x_\eta}\Big) \ind_{\{\eta+s<\He\}} \ind_{\{\tau^{*,x} > \eta\}}  \Big| \cF^2_\eta\Big]}{\P_{x,x}(\tau^{*,x} > \eta|\cF^2_\eta)}\\
&=\frac{\E_{x,x}\Big[\varphi\Big(\frac{1-\xi^x_{\eta + s}}{1-\xi^x_\eta}\Big) \Big((1-X_{\eta + s}) + X_{\eta + s} (1-\Delta\xi^x_0)\Gamma_*(M_\eta)/\Gamma_*(x) \Big) \ind_{\{\eta+s<\He\}} \Big|\cF^2_\eta\Big]}{1-X_\eta +X_\eta (1-\Delta\xi^x_0)\Gamma_*(M_\eta)/\Gamma_*(x)}\\
&=\E_{x,x}\bigg[\varphi\Big(\frac{1-\xi^x_{\eta + s}}{1-\xi^x_\eta}\Big) \frac{(1-X_{\eta + s}) + X_{\eta + s} \Gamma_*(M_\eta)}{1-X_\eta +X_\eta \Gamma_*(M_\eta)} \ind_{\{\eta+s<\He\}} \Bigg|\cF^2_\eta\Bigg]\\
&=
\hE_{x}^{q} \bigg[ \varphi\bigg(\frac{1-\xi^x_{\eta + s}}{1-\xi^x_\eta} \bigg) \ind_{\{\eta+s<\He\}} \bigg| \cF^2_\eta, \tau^{*,q} > \eta \bigg]
=
\hE_{x}^{q} \bigg[ \varphi\bigg(\frac{1-\xi^q_{\eta + s}}{1-\xi^q_\eta} \bigg) \ind_{\{\eta+s<\He\}} \bigg| \cF^2_\eta, \tau^{*,q} > \eta \bigg]
\end{aligned}
\end{equation}
where \eqref{eqn:xi-with-jump} and \eqref{eqn:aux_ind} were used in the third equality, \eqref{eqn:Deltaxi0} in the fourth and \eqref{eqn:aux_ind_full} in the fifth one (as $\Gamma_*(x) = \Psi(x,q)$). The last equality follows from \eqref{eqn:xi-with-jump} and \eqref{eqn:Deltaxi0} which yield the identities 
\[
\frac{1-\xi^x_{\eta + s}}{1-\xi^x_\eta} = \frac{\Gamma_*(M_{\eta + s})}{\Gamma_*(M_\eta)} = \frac{1-\xi^q_{\eta + s}}{1-\xi^q_\eta}
\] 
on the event $\{\eta+s<\He\}$. It remains to refer now to the first part of the proof as $(x,q) \in A_2$. The identity \eqref{eqn:eta_s_He} is justified analogously as in the case $(x,p) \in A_2$.

We omit the proof for $\sigma^*$ because the arguments are analogous to those above with the simplification that the generating process of $\sigma^p$ jumps from $0$ to $1$ at a hitting time for the process $(X_t, M_t)$.
\end{proof}

\begin{theorem}\label{thm:mainscaled}
The pair $((\sigma^{*,p}), (\tau^{*,p}))$ with $\sigma^*$ as in \eqref{eq:sigmadef6} and $\tau^*$ as in \eqref{eq:taudef6}, generated by $(\xi^p_t)$ as in \eqref{eq:xiA2} for $(x,p)\in A_2$ and as in \eqref{eqn:xi-with-jump} for $(x,p)\in A_1$, is a perfect Bayesian equilibrium for the game with payoff $\cP$. That is, for any
$(x,p) \in A_1 \cup A_2$ and any pair $(\sigma,\tau)\in\cT_R(\F^2)\times\cT_R(\F^1)$
\[
\hE^p_x\big[\cP(\sigma,\tau^*)\big]\le \hE^p_x\big[\cP(\sigma^*,\tau^*)\big]\le \hE^p_x\big[\cP(\sigma^*,\tau)\big].
\]
Moreover, the value of the game is equal to
\begin{equation}\label{eqn:Vxp}
V(x,p) =
\begin{cases}
x f(x) (1-\Gamma_*(x)) + \big(1-x + x \Gamma_*(x)\big) v_*(x,x), & (x,p) \in A_1,\\
v_*\big(x,\Gamma^{-1}_*(\Psi(x,p))\big), & (x,p) \in A_2.
\end{cases}
\end{equation}
\end{theorem}
\begin{proof}
The proof is divided into two parts, addressing separately $(x,p)\in A_1$ and $(x,p)\in A_2$.

\textbf{Part 1.} Fix $(x,p) \in A_2$ and take $m = \Gamma^{-1}_*(\Psi(x,p))$.
First, we show the optimality of $\tau^*$. By the definition of $\cP(\sigma,\tau)$, for any $\tau\in\cT_R(\F^1)$, recalling the
decomposition $\tau=\tau_0\ind_{\{\theta=0\}}+\tau_1\ind_{\{\theta=1\}}$, we have
\begin{equation}\label{eqn:v*lower}
\begin{aligned}
\hE^p_x\big[\cP(\sigma^*,\tau)\big]&=\int_0^1\hE^p_x\big[\ind_{\{\theta=0\}}\cP(\sigma^*,\tau_0(u))\big]\ud
u+\int_0^1\hE^p_x\big[\ind_{\{\theta=1\}}\cP(\sigma^*,\tau_1(u))\big]\ud u\\
&=(1-p)\int_0^1\E^0_{x}\big[\cP(\sigma^*,\tau_0(u))\big]\ud u+p\int_0^1\E^1_{x}\big[\cP(\sigma^*,\tau_1(u))\big]\ud u\\
&=(1-p)\int_0^1\E^0_{x,m}\big[\cP(\sigma_*,\tau_0(u))\big]\ud u+p\int_0^1\E^1_{x,m}\big[\cP(\sigma_*,\tau_1(u))\big]\ud u\\
&\ge (1-p)u^0_*(x,m)+p u_*^1(x,m)\\
&= v_*(x, m),
\end{aligned}
\end{equation}
where in the first equality we explicitly integrated out the randomisation device for $\tau$, in the second equality we used the
definition of the measure $\hP^p_x$ and the third one is based on three observations: first, the process $(X_t, m \vee \sup_{s \le t} X_s)$ under $\P^i_x$ has the same law on the path space as $(X_t, M_t)$ under $\P^i_{x,m}$, $i=0,1$;
second, $\sigma^*$, $\tau_0(u)$ and $\tau_1(u)$ are $\F^2$-stopping times, hence functions of the paths of $X$, and their laws are fully determined by the law of the pair $(X_t, m \vee \sup_{s \le t} X_s)$ under $\P^i_x$ and of the pair $(X_t, M_t)$ under $\P^i_{x,m}$, $i=0,1$;
third, for $(x,p) \in A_2$, the stopping time $\sigma_*$ coincides with $\sigma^*$.
The inequality follows from Proposition \ref{prop:u0u1}. The last equality above is because $p = \pi_*(x,m)$. Taking $\tau=\tau^*$, the inequality above becomes equality, thanks to Proposition \ref{prop:u0u1}, and we obtain
\[
\hE^p_x\big[\cP(\sigma^*,\tau^*)\big] = (1-p)u^0_*(x,m)+pu_*^1(x,m)=v_*(x,m),
\]
which proves the optimality of $\tau^*$.

Next, we prove the optimality of $\sigma^*$.
Let $\sigma\in\cT_R(\F^2)$ with $\sigma \le \He$ be arbitrary; we will use $\sigma(u)$ for arbitrary $u\in(0,1)$.
Consider \eqref{eqn:Lambda_m} with $t = \sigma(u)$ and $s = 0$. We have
\begin{equation}\label{eq:ito10a}
\begin{aligned}
&\hE^p_x\big[\Lambda_{\sigma(u)} - \Lambda_0 \big]\\
&=
\frac1{1-x + \Gamma_*(m) x} \E_{x} \bigg[ \hat v_*(X_{\sigma(u)}, M_{\sigma(u)}) - \hat v_*(x, m)
- \int_{[0,\sigma(u))} f(X_s) X_s \Gamma'_*(M_s) \ud M_s\bigg]\\
&=
\frac1{1-x + \Gamma_*(m) x} \E_{x}\bigg[ \int_{0}^{\sigma(u)} \partial_x \hat v_*(X_s,M_s) \ud X_s
+\frac12\int_{0}^{\sigma(u)} \partial_{xx} \hat v_*(X_s,M_s) \ud s\\
&\quad+\int_0^{\sigma(u)}\partial_{m} \hat
v_*(X_s,M_s)\ud M_s -\int_{[0,\sigma(u))}f(X_s) X_s \Gamma'_*(M_s) \ud M_s \bigg],
\end{aligned}
\end{equation}
where the second equality is by It\^o's formula. The integral with respect to $\ud X_s$ has zero expectation because $\partial_x
\hat v$ is bounded on $[0,1]$ and $(X_t)$ is a Brownian motion under $\P_{x}$. The integral with respect to $\ud s$ is
non-positive because of \eqref{eq:vsh},
which translates into $\partial_{xx} \hat v \le 0$. For the remaining
integral with respect to $\ud M_s$ we have:
\[
\int_0^{\sigma(u)}\!\! \big( \partial_{m} \hat v_*(X_s,M_s)\! -\! f(X_s) X_s \Gamma'_*(M_s) \big) \ud M_s
= \int_0^{\sigma(u)}\!\! \big( \partial_{m} \hat v_*(M_s,M_s)\! -\! f(M_s) M_s \Gamma'_*(M_s) \big) \ud M_s
= 0,
\]
where the first equality uses that $(M_t)$ increases only when $X_s = M_s$ and the last equality holds by \eqref{eqn:m_tildev}, where we can treat $m < m^1_*$ and $m>m^1_*$ separately because the measure $\ud M_s$ is continuous.

In conclusion, we have
\[
\hE^p_x\big[\Lambda_{\sigma(u)} - \Lambda_0 \big] \le 0,
\]
and so
\begin{equation}\label{eqn:vLa}
v_*(x, m) = 
\Lambda_0 \ge \hE^p_x[\Lambda_{\sigma(u)}] .
\end{equation}

Next we show that $\hE^p_x[\Lambda_{\sigma(u)}] \ge \hE^p_x[\cP(\sigma(u), \tau^*)]$, from which it will follow that $v_*(x,m) \geq \hE^p_x [\cP(\sigma,\tau^*)]$. Recall that
\begin{equation}\label{eqn:PP}
\hE^p_x[\cP(\sigma(u), \tau^*)] = \hE^p_x \big[ \ind_{\{\tau^* < \sigma(u)\}} f(X_{\tau^*})
+ \ind_{\{\tau^* = \sigma(u)\}} h(X_{\sigma(u)})
+ \ind_{\{\tau^* > \sigma(u)\}} g(X_{\sigma(u)})\big]
\end{equation}
and (cf.\ \eqref{eqn:Lambda})
\begin{equation}\label{eqn:Elambda}
\hE^p_x[\Lambda_{\sigma(u)}] = \hE^p_x \bigg[ \ind_{\{\tau^* \ge {\sigma(u)}\}} \!\!\left( \ind_{\{\theta=1\}} u^1_* \!+\!
\ind_{\{\theta=0\}} u^0_* \right)(X_{{\sigma(u)}}, M_{{\sigma(u)}})\! +\! \ind_{\{\theta=1\}}\int_{[0,{\sigma(u)})}\!\!f(X_s) \ud\xi_s  \bigg].
\end{equation}
We begin by addressing the last term on the right-hand side of \eqref{eqn:Elambda} showing that it equals the first term on the right-hand side of \eqref{eqn:PP}.
Since $(\xi_t)$ from \eqref{eq:xiA2} generates the randomised stopping time $\tau^*_1$, using the well-known change-of-variable formula \cite[Ch.\ 0,
Prop.\ 4.9]{revuzyor} we get
\begin{equation}
\begin{aligned}
\ind_{\{\theta=1\}}\int_{[0,\sigma(u))}f(X_s) \ud\xi_s &=
\ind_{\{\theta=1\}}\int_0^1\ind_{\{\tau^*_1(r)<\sigma(u)\}}f(X_{\tau^*_1(r)})\ud r.
\end{aligned}
\end{equation}
Taking expectations yields
\begin{equation}
\begin{aligned}
&\hE^p_x\Big[\ind_{\{\theta=1\}} \int_{[0,\sigma(u))}f(X_s)\ud\xi_s\Big]\\
&=\hE^p_x\Big[\ind_{\{\theta=1\}} \int_0^1\ind_{\{\tau^*_1(r)<\sigma(u)\}}f(X_{\tau^*_1(r)})\ud r\Big]
=\hE^p_x \Big[\ind_{\{\tau^*< \sigma(u)\}}f(X_{\tau^*})\Big],
\end{aligned}
\end{equation}
where, for the last equality, we recall that $\ind_{\{\theta=1\}}\ind_{\{\tau^*_1<\sigma(u)\}}=\ind_{\{\theta=1\}}\ind_{\{\tau^*<\sigma(u)\}}$, because $\tau^* = \He \ge \sigma(u)$ on $\{\theta=0\}$.

For the second term under expectation in \eqref{eqn:PP} we notice that $\hP^p_x$-a.s.
\[
\{\sigma(u)=\tau^*\}\cap\{\theta=0\} =\{\sigma(u)=\H_0\}\cap\{\theta=0\}
\]
and also $\hP^p_x$-a.s.
\begin{align*}
\{\sigma(u)=\tau^*\}\cap\{\theta=1\}
&= \{\sigma(u)=\tau_1^*\} \cap\{\theta=1\}
= \{\sigma(u)=\H_1\}\cap\{\theta=1\},
\end{align*}
where the last equality follows from the path-by-path continuity of the process $(\xi_t)$ for $t \in [0, \H_1)$ (cf.\ \eqref{eq:xiA2}), which implies
\begin{equation} \label{eq:refforabstract}
\hP^p_x(\sigma(u)=\tau_1^* < \H_1)= \hP^p_x(\Delta \xi_{\sigma(u)} > 0 \text{ and } \sigma(u) < \H_1) = 0.
\end{equation}
Then, we have $\hP^p_x$-a.s.
\begin{equation}\label{eq:v=h}
\begin{aligned}
\ind_{\{\sigma(u)= \tau^*\}} h(X_{\sigma(u)})
&=\ind_{\{\sigma(u)= \tau^*\}\cap\{\theta=0\}}h(X_{\sigma(u)})+\ind_{\{\sigma(u)= \tau^*\}\cap\{\theta=1\}}h(X_{\sigma(u)})\\
&=\ind_{\{\sigma(u)= \tau^* = \He \}\cap\{\theta=0\}}g(0)+\ind_{\{\sigma(u)= \tau^*= \He \}\cap\{\theta=1\}}f(1)\\
&=\ind_{\{\sigma(u)= \tau^* = \He\}\cap\{\theta=0\}}u^0_*(0,M_{\sigma(u)})+\ind_{\{\sigma(u)= \tau^* =
\He\}\cap\{\theta=1\}}u^1_*(1,1)\\
&=\ind_{\{\sigma(u)= \tau^* = \He\}} \big( \ind_{\{\theta=0\}}
u^0_*(X_{\sigma(u)},M_{\sigma(u)})+\ind_{\{\theta=1\}}u^1_*(X_{\sigma(u)},M_{\sigma(u)})\big),
\end{aligned}
\end{equation}
where the second equality holds because $g(0)=h(0)$ and $f(1)=h(1)$ (Assumption \ref{ass:fghscaled}-(ii)), the third one because $u^0_*(z,\ell)=g(z)$ for $z\le d_*(\ell)$ and by the third condition in \eqref{eq:PDEu1}.

For the final term under expectation in \eqref{eqn:PP}, using \eqref{eqn:RN_xp} for the second equality below with $\Gamma_*(M_0)=\Gamma_*(m)=\Psi(x,p)$ yields
\begin{equation*}
\begin{aligned}
&\hE^p_x \big[\ind_{\{\tau^* > \sigma(u)\}} g(X_{\sigma(u)}) \big]\\
&= \hE^p_x \Big[ \ind_{\{\theta=1\}} \ind_{\{\tau^* > \sigma(u)\}} g(X_{\sigma(u)}) + \ind_{\{\theta=0\}} \ind_{\{\sigma(u) <
\He\}} g(X_{\sigma(u)})\Big]\\
&=\frac1{1 - x + x \Gamma_*(m)} \E_{x} \Big[\Gamma_*(m) \ind_{\{\theta=1\}} (1-\xi_{\sigma(u)}) g(X_{\sigma(u)})
+ \ind_{\{\theta=0\}} \ind_{\{\sigma(u) < \He\}} g(X_{\sigma(u)}) \Big]\\
&=\frac1{1 - x + x \Gamma_*(m)} \E_{x}\bigg[ \Gamma_*(m)  \ind_{\{\sigma(u) < \He\}} X_{\sigma(u)}
\frac{\Gamma_*(M_{\sigma(u)})}{\Gamma_*(m)} g(X_{\sigma(u)}) + \ind_{\{\sigma(u) < \He\}} (1-X_{\sigma(u)}) g(X_{\sigma(u)})
\bigg]\\
&= \frac{1}{1 - x + x \Gamma_*(m)} \E_{x} \bigg[ \ind_{\{\sigma(u) < \He\}} \bigg( X_{\sigma(u)}\Gamma_*(M_{\sigma(u)}) + 1 - X_{\sigma(u)}\bigg) g(X_{\sigma(u)}) \bigg]\\
& \leq\frac1{1 - x + x \Gamma_*(m)} \E_{x}\Big[ \ind_{\{\sigma(u) < \He\}} \Big(X_{\sigma(u)} \Gamma_*(M_{\sigma(u)}) + 1-X_{\sigma(u)}\Big)
v_*(X_{\sigma(u)}, M_{\sigma(u)}) \Big] \\
&=\frac1{1 - x + x \Gamma_*(m)} \E_{x}\Big[ \ind_{\{\sigma(u) < \He\}} \Big(X_{\sigma(u)} \Gamma_*(M_{\sigma(u)})u^1_*(X_{\sigma(u)},M_{\sigma(u)}) + (1-X_{\sigma(u)})u^0_*(X_{\sigma(u)},M_{\sigma(u)})\Big)\Big]\\
&=\frac1{1 - x + x \Gamma_*(m)} \E_{x}\Big[ \ind_{\{\sigma(u) < \He\}} \Big(\ind_{\{\theta=1\}} \Gamma_*(M_{\sigma(u)})u^1_*(X_{\sigma(u)},M_{\sigma(u)}) + \ind_{\{\theta=0\}}u^0_*(X_{\sigma(u)},M_{\sigma(u)})\Big)\Big],
\end{aligned}
\end{equation*}
where we used \eqref{eq:vgL} for the inequality, \eqref{eq:v*def} for the penultimate equality and tower property combined with $X_{\sigma(u)}=\P_x(\theta=1|\cF^2_{\sigma(u)})$, for the final equality. Then, using again
\[
\frac{\Gamma_*\big(M_{\sigma(u)}\big)}{\Gamma_*(m)}=1-\xi_{\sigma(u)}=\P_x\big(\tau^*_1>\sigma(u)\big|\cF^2_\infty\big)
\]
and the change of measure \eqref{eqn:RN_xp}, the final expression above turns into
\begin{equation*}
\begin{aligned}
&\hE^p_x \Big[ \ind_{\{\sigma(u) < \He\}} \ind_{\{\tau^* \ge \sigma(u)\}} \big(\ind_{\{\theta=1\}} u^1_*(X_{\sigma(u)},
M_{\sigma(u)}) + \ind_{\{\theta=0\}} u^0_*(X_{\sigma(u)}, M_{\sigma(u)})\big)\Big]\\
&= \hE^p_x \Big[ \ind_{\{\tau^* > \sigma(u)\}} \big(\ind_{\{\theta=1\}} u^1_*(X_{\sigma(u)}, M_{\sigma(u)}) +
\ind_{\{\theta=0\}} u^0_*(X_{\sigma(u)}, M_{\sigma(u)})\big)\Big],
\end{aligned}
\end{equation*}
where the equality holds because the process $(\xi_t)$ is continuous for $t < \He$ and we have $\hP^p_x$\as
\[
\{\sigma(u) < \He\} \cap \{\tau^* \ge \sigma(u)\} = \{\sigma(u) < \He\} \cap \{\tau^* > \sigma(u)\} = \{\tau^* > \sigma(u)\}.
\]
In conclusion, we have obtained $\hE^p_x[\cP(\sigma(u), \tau^*)] \le \hE^p_x[\Lambda_{\sigma(u)}] \le v_*(x, m)$ for each $u \in (0,1)$ with the last inequality by \eqref{eqn:vLa}, so $\sup_{\sigma \in \cT_R(\F^2)} \hE^p_x[\cP(\sigma, \tau^*)] \le v_*(x,m)$.

To prove the optimality of $\sigma^*$, notice that $\partial_{xx} \hat v_*(X_s, M_s) = 0$ for $s \le
\sigma^*$. Then, taking $\sigma(u)=\sigma^*$ in \eqref{eq:ito10a} we obtain the analogue of \eqref{eqn:vLa} but with equality. That is,
$v_*(x, m) = \Lambda_0 = \hE^p_x[\Lambda_{\sigma^*}]$. It remains to show the equality $\hE^p_x[\Lambda_{\sigma^*}] =
\hE^p_x[\cP(\sigma^*, \tau^*)]$. For this, we only need that the above arguments for the third term of \eqref{eqn:PP} hold with
equality. Indeed,
\begin{align*}
&\hE^p_x \big[\ind_{\{\tau^* > \sigma^*\}} g(X_{\sigma^*}) \big]\\
&= \hE^p_x \Big[ \ind_{\{\theta=1\}} \ind_{\{\tau^* > \sigma^*\}} g(X_{\sigma^*}) + \ind_{\{\theta=0\}} \ind_{\{\sigma^* <
\He\}} g(X_{\sigma^*})\Big]\\
&= \hE^p_x \Big[ \ind_{\{\theta=1\}} \ind_{\{\tau^* > \sigma^*\}} u^1_*(X_{\sigma^*}, M_{\sigma^*}) + \ind_{\{\theta=0\}}
\ind_{\{\sigma^* < \He\}} u^0_*(X_{\sigma^*}, M_{\sigma^*})\Big]\\
&= \hE^p_x \Big[ \ind_{\{\tau^* > \sigma^*\}} \big(\ind_{\{\theta=1\}} u^1_*(X_{\sigma^*}, M_{\sigma^*}) + \ind_{\{\theta=0\}}
u^0_*(X_{\sigma^*}, M_{\sigma^*})\big)\Big]
\end{align*}
where the second equality is by the third conditions in \eqref{eq:PDEu0} and in \eqref{eq:PDEu1}. This completes the proof of the optimality of $\sigma^*$ and justifies the second line of \eqref{eqn:Vxp}.
\medskip

\textbf{Part 2.} Fix $(x,p) \in A_1$ (hence $p=x$ and $x>m^1_*$).
First, we show the optimality of $\tau^*$. Recall that Player 2 models the belief process $(p_t)$ assuming that Player 1 follows
$\tau^*$, i.e., $p_{0-} = p =x>\pi_*(x,x)$ and $p_0 = \pi_*(x,x)$, i.e., $(x, p_0) \in A_2$. It is useful to recall the discussion following \eqref{eq:sigmadef6}. In particular, $\hP^p_x(\sigma^*>0)=1$ for $(x,p)\in A_1$ and $M_0=x$. Then,
\begin{equation}
\sigma^*=\inf\Big\{t\ge 0: X_t\le d_*\Big(\sup_{s\le t}X_s\Big)\Big\}\wedge\He.
\end{equation}
Furthermore, using Proposition \ref{prop:u0u1}, we have, for $i=0,1$,
\begin{equation}\label{eq:sigmaspec}
\E^i_x\big[\cP(\sigma^*,\tau_i)\big]\ge u^i_*(x,x),\quad \tau_i\in\cT_R(\F^2).
\end{equation}

Take any $\tau\in\cT_R(\F^1)$ and recall its
decomposition $\tau=\tau_0\ind_{\{\theta=0\}}+\tau_1\ind_{\{\theta=1\}}$ for $\tau_0, \tau_1 \in \cT_R(\F^2)$. For any $u \in
(0,1)$, we have
\begin{equation}\label{eq:verif00}
\begin{aligned}
\hE^x_x \big[\cP(\sigma^*,\tau(u))\big]&=\hE^x_x\big[\ind_{\{\theta=0\}}\cP(\sigma^*,\tau_0(u))+\ind_{\{\theta=1\}}\cP(\sigma^*,\tau_1(u))\big]\\
&=(1-x)\Big(\ind_{\{\tau_0(u) = 0\}} \E^0_x[\cP(\sigma^*,0)] + \ind_{\{\tau_0(u) > 0\}}
\E^0_{x}\big[\cP(\sigma^*,\tau_0(u))\big] \Big)\\
&\quad + x \Big(\ind_{\{\tau_1(u) = 0\}} \E^1_x[\cP(\sigma^*,0)] + \ind_{\{\tau_1(u) > 0\}}
\E^1_{x}\big[\cP(\sigma^*,\tau_1(u))\big] \Big)\\
&\ge (1-x) \min\big(\E^0_x[\cP(\sigma^*,0)], u_*^0(x,x)\big)+x \min\big(\E^1_x[\cP(\sigma^*,0)], u_*^1(x,x)\big),
\end{aligned}
\end{equation}
where the second equality is by the definition $\hP^p_x$ (cf. Lemma \ref{lem:Y-Pyp}) and $x=p$ and the inequality holds by \eqref{eq:sigmaspec}. Using that $\P^i_x(\sigma^*>0)=1$ for $i=0,1$ then
\[
\E^0_x[\cP(\sigma^*,0)] = \E^1_x[\cP(\sigma^*,0)] = f(x).
\]
Using also that $u^0_*(x,x) < f(x)$ and $u^1_*(x,x) =
f(x)$, we obtain
\begin{equation*}
\begin{aligned}
\hE^x_x\big[\cP(\sigma^*,\tau(u))\big]
&\ge (1-x) u^0_*(x,x) + x u^1_*(x,x)\\
&=(1-x) u^0_*(x,x) + x\Gamma_*(x)u^1_*(x,x) + x\big(1-\Gamma_*(x)\big) f(x) \\
&= \big(1-x + x \Gamma_*(x)\big) v_*(x,x) + x \big(1-\Gamma_*(x)\big) f(x) = V(x,p)
\end{aligned}
\end{equation*}
with $V$ defined in \eqref{eqn:Vxp}.
Taking $\tau = \tau^*$ the inequality in \eqref{eq:verif00} becomes equality because $\ind_{\{\tau^*_0=0\}}=0$ and
\begin{equation*}
\begin{aligned}
\hE^x_x\big[\cP(\sigma^*,\tau^*)\big]&=(1-x)
\E^0_{x}\big[\cP(\sigma^*,\He)\big]\\
&\quad + x \Big(\E^1_x[\cP(\sigma^*,0)]\int_0^1\ind_{\{\tau^*_1(r) = 0\}}\ud r  + \int_0^1\ind_{\{\tau^*_1(r) > 0\}}\E^1_{x}\big[\cP(\sigma^*,\tau^*_1(r))\big]\ud r\Big)\\
 &=(1-x)u^0_*(x,x)+x \Big(\big(1-\Gamma_*(x)\big)f(x)+ \int_0^1\ind_{\{\tau^*_1(r) > 0\}}u^1_*(x,x)\ud r \Big)\\
 &=(1-x)u^0_*(x,x)+x \big(1-\Gamma_*(x)\big)f(x)+ x\Gamma_*(x)u^1_*(x,x)=V(x,p),
\end{aligned}
\end{equation*}
where for the second equality we used
\[
\int_0^1\ind_{\{\tau^*_1(r) = 0\}}\ud r =\Delta\xi_0= 1-\Gamma_*(x)\quad\text{and}\quad \E^1_{x}\big[\cP(\sigma^*,\tau^*_1(r))\big]=u^1_*(x,x),
\]
(cf.\ Proposition \ref{prop:u0u1}) and for the third one we used
\[
\int_0^1\ind_{\{\tau^*_1(r) > 0\}}\ud r=1-\Delta\xi_0=\Gamma_*(x).
\]
This proves optimality of $\tau^*$.

To prove the optimality of $\sigma^*$, we need to show that for any stopping time $\sigma \in \cT_R(\F^2)$ and any $u \in (0,1)$ we have
$\hE^p_x[\cP(\sigma(u), \tau^*)] \le V(x,p)$ for $p=x$ and that the equality is attained for a pure stopping time $\sigma = \sigma^*$. Since $\sigma(u) \in \cT(\F^2)$ for any $u \in (0,1)$, it is enough to prove the inequality $\hE^p_x[\cP(\sigma, \tau^*)] \le V(x,p)$ for $p=x$ for any pure stopping time $\sigma \in \cT(\F^2)$. Fix $\sigma \in \cT(\F^2)$. We have the decomposition
\begin{align*}
\hE^p_x[\cP(\sigma, \tau^*)] = \ind_{\{\sigma = 0\}} \hE^p_x[\cP(0, \tau^*)] + \ind_{\{\sigma > 0\}} \hE^p_x[\cP(\sigma,\tau^*)].
\end{align*}
The first term is evaluated easily as (recall that $x=p$)
\begin{equation}
\begin{aligned}
\hE^p_x[\cP(0, \tau^*)]
&= x \big(\Delta \xi_0 h(x) + (1-\Delta\xi_0)g(x)\big) + (1-x)g(x)\\
&\le x \Delta \xi_0 f(x) + \big(1-x + x(1-\Delta\xi_0)\big) g(x) \le V(x,p),
\end{aligned}
\end{equation}
with the last inequality following from $v_*(x,x) \ge g(x)$ by \eqref{eq:vgL}.

Computations for the second term are more involved. Recalling that we are working on the set $\{\sigma > 0\}$, and $(\xi_t)$ is continuous on $(0, \H_1)$ we have
\begin{align*}
&\hE^p_x[\cP(\sigma, \tau^*)]\\
&= \hE^p_x \bigg[ \ind_{\{\theta=0\}} g(X_{\sigma}) + \ind_{\{\theta=1\}}
\bigg( \int_{[0, \sigma)} f(X_s) \ud\xi_s + \ind_{\{\sigma = \He\}} h(X_{\sigma}) \Delta \xi_{\sigma}+ g(X_{\sigma}) (1-\xi_{\sigma}) \bigg)\bigg]\\
&= \E_x \bigg[ \ind_{\{\theta=0\}} g(X_\sigma)\\
&\hspace{30pt}+ \ind_{\{\theta=1\}}\bigg( f(x) \Delta \xi_0 +
\int_{(0, \sigma)} f(X_s) \ud\xi_s + \ind_{\{\sigma = \He\}}h(X_\sigma) \Delta \xi_{\sigma} + g(X_{\sigma}) (1-\xi_{\sigma})\bigg)\bigg]\\
&= x f(x) \Delta \xi_0 + \E_x \bigg[ \ind_{\{\theta=0\}} g(X_\sigma)\bigg]\\
&\quad +\E_x \bigg[ \ind_{\{\theta=1\}}
(1-\Delta\xi_0) \bigg(\int_{[0, \sigma)} f(X_s) \ud\hat\xi_s + \ind_{\{\sigma = \He\}}h(X_\sigma) \Delta \hat\xi_{\sigma}
+ g(X_{\sigma}) (1-\hat\xi_{\sigma})\bigg)\bigg],
\end{align*}
where the second equality holds because $\hP_{x}^{x}=\P_x$ as noted above, and $(\hat \xi_t)$ is defined as
\[
\hat \xi_t = \frac{\xi_t - \Delta \xi_0}{1-\Delta \xi_0}, \quad \hat \xi_{0-} = 0.
\]
By \eqref{eqn:xi-with-jump}, we have
$\hat \xi_t= 1 - \Gamma_*(M_t)/\Gamma_*(M_0)$,
where $M_0 = x$ and $M_t = \sup_{s \le t} X_s$. Using $\Delta\xi_0=1-\Gamma_*(x)$ and recalling $q$ defined in \eqref{eqn:qq} we get
\begin{equation}
\begin{aligned}
\E_x \bigg[ \ind_{\{\theta=0\}} g(X_\sigma)\bigg]=\big(1-x+x\Gamma_*(x)\big)\E_x \bigg[ \ind_{\{\theta=0\}}\frac{1-q}{1-x} g(X_\sigma)\bigg]
\end{aligned}
\end{equation}
and
\begin{equation}
\begin{aligned}
&\E_x \bigg[ \ind_{\{\theta=1\}}
(1-\Delta\xi_0) \bigg(\int_{[0, \sigma)} f(X_s) \ud\hat\xi_s + \ind_{\{\sigma = \He\}}h(X_\sigma) \Delta \hat\xi_{\sigma}
+ g(X_{\sigma}) (1-\hat\xi_{\sigma})\bigg)\bigg]\\
&=\Gamma_*(x)\E_x \bigg[ \ind_{\{\theta=1\}}
\bigg(-\int_{[0, \sigma)} f(X_s) \frac{\Gamma'_*(M_s)}{\Gamma_*(x)}\ud M_s + \ind_{\{\sigma = \He\}}h(X_\sigma) \frac{\Gamma_*(M_\He)}{\Gamma_*(x)}
+ g(X_{\sigma})\frac{\Gamma_*(M_\sigma)}{\Gamma_*(x)}\bigg)\bigg]\\
&=\big(1-x+x\Gamma_*(x)\big)\E_x \bigg[ \ind_{\{\theta=1\}}\frac{q}{x}
\bigg(-\int_{[0, \sigma)} f(X_s) \frac{\Gamma'_*(M_s)}{\Gamma_*(x)}\ud M_s + \ind_{\{\sigma = \He\}}h(X_\sigma) \frac{\Gamma_*(M_\He)}{\Gamma_*(x)}\\
&\hspace{170pt}+ g(X_{\sigma})\frac{\Gamma_*(M_\sigma)}{\Gamma_*(x)}\bigg)\bigg].
\end{aligned}
\end{equation}
Adding up the two expressions derived above and letting $\hat \tau$ be the randomised stopping time generated by $(\hat \xi_t)$  we get
\begin{equation}
\begin{aligned}
&\E_x \bigg[ \ind_{\{\theta=0\}} g(X_\sigma)+\ind_{\{\theta=1\}}
(1-\Delta\xi_0) \bigg(\int_{[0, \sigma)} f(X_s) \ud\hat\xi_s + \ind_{\{\sigma = \He\}}h(X_\sigma) \Delta \hat\xi_{\sigma}
+ g(X_{\sigma}) (1-\hat\xi_{\sigma})\bigg)\bigg]\\
&=\big(1-x+x\Gamma_*(x)\big)\E_x\bigg[\frac{\ud \hP_{x}^{q}}{\ud \P_x}\Big(\ind_{\{\theta=0\}}\cP(\sigma,\He)+\ind_{\{\theta=1\}}\cP(\sigma,\hat \tau)\Big)\bigg]\\
&=\big(1-x+x\Gamma_*(x)\big)\hE_{x}^{q}\big[\ind_{\{\theta=0\}}\cP(\sigma,\He)+\ind_{\{\theta=1\}}\cP(\sigma,\hat \tau)\big]\\
&=\big(1-x+x\Gamma_*(x)\big)\hE_{x}^{q}\big[\cP(\sigma,\tau^*)\big],
\end{aligned}
\end{equation}
where in the last equality we used that the process $(\hat \xi_t)$ is the same (path-by-path) as the generating process of $\tau^*_1$ in \eqref{eq:xiA2} starting from initial condition $(x,q)\in A_2$ (recall $\Psi(x,q)=\Gamma_*(x)$).

Thus we have obtained
\begin{equation*}
\hE^p_x[\cP(\sigma, \tau^*)]= x (1-\Gamma_*(x)) f(x) + \big(1-x+x \Gamma_*(x)\big) \hE_{x}^{q} [ \cP(\sigma, \tau^*) ].
\end{equation*}
Because $(x,q) \in A_2$ with $\Psi(x,q)=\Gamma_*(x)$, the first part of the proof yields $\hE_{x}^{q} [ \cP(\sigma, \tau^*)] \le v_*(x,x)$ and $\hE_{x}^{q} [ \cP(\sigma^*, \tau^*)] = v_*(x,x)$. So, recalling that $x=p$,
\[
\hE^p_x[\cP(\sigma, \tau^*)] \le x (1-\Gamma_*(x))f(x) + \big(1-x+x \Gamma_*(x)\big) v_*(x,x) = V(x,p)=\hE^p_x[\cP(\sigma^*, \tau^*)].
\]
That shows \eqref{eqn:Vxp} for $(x,p)\in A_1$ and $\sigma^*$ is optimal.
\end{proof}

\begin{remark}\label{rem:Vxx}
We can strengthen the statement of Remark \ref{rem:3.6} to a general setting of our game, i.e., we can show that the additional information available to Player 1 lowers her value compared to the game with symmetric information. At the inception of the game, Player 2 has no special information about $X_{\He}$, so $(x,p) = (x,x)$ and the value is $V(x,x)$ for the function $V$ defined in \eqref{eqn:Vxp}. We will show below that $V(x,x) < f(x)$ for all $x \in (0,1)$. However, in the game without a foresight, whenever equilibrium involves stopping by the minimiser, the value of the game coincides with $f$ on her stopping set, hence it is strictly larger there than $V(x,x)$.

We divide the proof that $V(x,x) < f(x)$ into two parts. When $(x,x) \in A_1$, i.e., $x > m^1_*$, the value is given by
\begin{equation}\label{eqn:Vxx-A1}
V(x,x) = x f(x) + (1-x) \Du0_{c_*}(x).
\end{equation}
Indeed, we have
$v_*(x,x) = \pi_*(x,x) u^1_*(x,x) + (1 - \pi_*(x,x)) u^0_*(x,x)$.
Moreover, $x > m^1_*$ implies $u^1_*(x,x) = f(x)$ by \eqref{eq:PDEu1} and $u_*^0(x,x) = \Du0_{c_*}(x)$ by Proposition \ref{prop:u0}. Inserting $\pi_*=\pi_{\Gamma_*}$ from \eqref{eq:pi} into \eqref{eqn:Vxp} and simplifying yields \eqref{eqn:Vxx-A1} for $V(x,x)$. The continuity of $u^0$ and Corollary \ref{cor:upbu0} give $u^0_*(x,x) \le g(x)$, so
\[
V(x,x) \le x f(x) + (1-x) g(x) < f(x), \quad x \in (m^1_*,1),
\]
where we used that $g(x) < f(x)$ for $x \in (0,1]$, cf. Assumption \ref{ass:fghscaled}-(iii).

When $(x,x) \in A_2$, i.e., $x \le m^1_*$, we have
\[
V(x,x) = v_*(x,m^1_*) = x u^1_*(x,m^1_*) + (1 - x) u^0_*(x,m^1_*).
\]
Corollary \ref{cor:upbu0} and \eqref{eq:u0u1*} imply $u^0_*(x,m^1_*) \le g(x)$  and Lemma \ref{lem:apriori} yields $u^1_*(x,m^1_*) \le f(x)$. This gives the bound $V(x,x) \le xf(x) + (1-x) g(x) < f(x)$, where we again use $g(x) < f(x)$ for $x \in (0,1]$.
\end{remark}

\section{Examples}\label{sec:examples}

In examples below, function $h$ is any continuous function such that $g \le h \le f$ and $h(0) = g(0)$ and $h(1) = f(1)$.
\subsection{Linear payoffs}\label{subsec:linear payoffs}
In this example we assume that $f(x)=x$ and $g(x)=\alpha x$ for $\alpha \in (0,1)$. 
We first obtain $c_*=c_*(m)$ as a root of $\phi(\cdot,m)$, where $\phi$ is given in \eqref{eq:phidef}. Fix $m \in
(0,1]$. Then $c=c_*(m)$ solves
\[ m = \frac{ c}{ m - c } (m - \alpha c). \]
This equation has solutions $c = \lambda_\pm m$ where $\lambda_\pm$ are the roots of $Q(\lambda)=0$ and $Q(\lambda) =1 - 2 \lambda
+ \alpha \lambda^2$. The two roots are given by $\lambda_\pm = \frac{1}{\alpha} \pm \frac{1}{\alpha} \sqrt{1- \alpha}$. We want
the (unique) root in $[0,1]$ (since we must have $0 \leq c \leq m$) and hence
\[ c_*(m) = \lambda m; \hspace{20mm} \lambda = \frac{1}{\alpha} - \frac{1}{\alpha} \sqrt{1- \alpha}.  \]
Note that since $Q(\frac{1}{2}) = \frac{\alpha}{4}>0$ and $Q(1) = \alpha-1<0$ we must have $\lambda \in (\frac{1}{2},1)$.

The second task is to calculate $\Du0_{c_*}(m)=u^0_*(m,m)$ given by \eqref{eq:Ddef}. Using the expression for $c_*(m)$ we find
\[
\Du0'_{c_*}(m) = \frac{(1 - \lambda m)}{m(1-m)(1-\lambda)} (\Du0_{c_*}(m) - \alpha \lambda m).
\]
This is a linear equation which can be integrated (using the boundary condition $\Du0_{c_*}(1) = \alpha \lambda$) to obtain
\begin{eqnarray*}
\Du0_{c_*}(m) & = & \frac{\lambda \alpha}{1 - \lambda} \frac{m^{1/(1-\lambda)}}{(1-m)} \int_m^1 (1- \lambda y) y^{-1/(1-\lambda)} \ud y \\
& = & \frac{\alpha}{1-m} \left(m- m^{1/(1-\lambda)} \right)  - \frac{1}{1-m} \left( m^2 - m^{1/(1-\lambda)} \right),
\end{eqnarray*}
where we use $\alpha = \frac{2 \lambda - 1}{\lambda^2}$ to simplify this expression.

The third task is to solve for $\gamma_*=\gamma_{c_*}$, where $\gamma_{c_*}(m)$ is given in \eqref{eq:gamma*def}. After some algebra we find 
\[ \gamma_*(m)= 1 - \frac{1}{2 \lambda} m^\frac{2 \lambda-1}{1 - \lambda}.  \]
Note that $\frac{2 \lambda-1}{1 - \lambda}>0$ and thus that $\gamma_*(0)=1$, $\gamma_*(1)=\frac{2 \lambda - 1}{2 \lambda}$ and
$\gamma_*$ is decreasing on $(0,1]$. In particular, $\Gamma_* = \gamma_{c_*}$ and $m^1_* = 0$ so that $d_*=c_*$ on $(0,1]$.

Value functions for the players take the following explicit form: for $\lambda m \leq x \leq m \leq 1$:
\begin{eqnarray*}
u^0_*(x,m) & = & \alpha \lambda m   + \frac{(x - \lambda m)(\alpha(1-\lambda) - m (1-\alpha\lambda) + (1- \alpha)m^{\frac{\lambda}{1 - \lambda}})}{(1 - \lambda)(1-x)} , \\
u^1_*(x,m) & = & \alpha \lambda m   + \frac{m(1 - \alpha \lambda)(x - \lambda m)}{(1 - \lambda)x} , \\
v_*(x,m) & = & (1 - \pi_*(x,m)) u^0_*(x,m) + \pi_*(x,m) u^1_*(x,m)
\end{eqnarray*}
where $\pi_*(x,m)  =  \frac{x \gamma_*(m)}{ 1-x + x \gamma_*(m)}$.

Note that $\frac{\lambda}{1-\lambda}>1$ and so for $m \in (0,1)$ we have $m^{\frac{\lambda}{1-\lambda}} < m$ and the numerator in the second element of the expression for $u^0_*$ is bounded above by
\[ (x - \lambda m)(\alpha(1-\lambda) - m (1-\alpha\lambda) + (1- \alpha)m) = (x - \lambda m)\alpha(1-\lambda)(1-m) \]
It follows that $\Du0_{c_*}(m) = u^0_*(m,m) < \alpha \lambda m + \frac{m(1-m)\alpha(1-\lambda)^2}{(1-\lambda)(1-m)} = \alpha m = g(m)$, which can also be deduced from Corollary~\ref{cor:upbu0}.

Since $m^1_* = 0$, $(x, x) \in A_1$ for any $x \in (0,1)$. Hence, $V(x,x)$ is given by \eqref{eqn:Vxx-A1}.  As $\Du0_{c_*}(x) \le g(x)$, we get $V(x,x) \le x^2 + (1-x)\alpha x < x$ for $x \in (0,1)$. This substantiates the claim of Remark \ref{rem:3.6} that the more informed Player 1 is able to lower the value of the game compared to the game with symmetric and full information.

\textbf{The equilibrium.}
Since $m^1_*=0$, every initial state $(x,x)$, $x\in(0,1)$, lies in $A_1$, and
$d_*=c_*$. The equilibrium therefore has the following simple form. Type $\theta=0$ of Player 1 never stops. Type $\theta=1$ of Player 1 first stops at time $0$ with probability
\[
1-\Gamma_*(x)=\frac{1}{2\lambda}\,x^{\frac{2\lambda-1}{1-\lambda}} .
\]
Afterwards she stops gradually only when $X$ reaches a new maximum according to the randomised rule
\[
\inf\Big\{t\ge 0:\tfrac{1}{2\lambda}M^{\frac{2\lambda-1}{1-\lambda}}_t>U\Big\}\wedge\He\quad\text{with}\quad M_t=\sup_{s\le t}X_s.
\]

Player 2, the maximiser, stops at the first time the process falls a fixed proportion below its running maximum,
\[
\sigma_*=\inf\{t\ge0:\ X_t\le\lambda M_t\}\wedge\He.
\]

If $\theta=1$ then the game may last until $\He=\H_1$ with neither player stopping. Indeed,
by excursion theory (see Rogers~\cite{rogers1989excursions} for an introduction) and the fact that for $a<y<b$ we have $\P^1_y(\H_a < \H_b) = \frac{a(b-y)}{y(b-a)}$ and then
$\P^1_y(\H_{\lambda y} < \H_{y + \ud y} ) \sim \frac{\lambda}{1-\lambda} \frac{\ud y}{y}$, we have
\( \P^1_x(\sigma_*=\H_1)
= \exp\left( - \int_x^1 \frac{\lambda}{1-\lambda} \frac{\ud y}{y} \right) = x^{\frac{\lambda}{1-\lambda}}. \)
It follows that
\[
\P^1_x\big(\tau^*_1=\sigma_*=\H_1\big)=\Gamma_*(1)\,\P^1_x(\sigma_* = \H_1)
=\frac{2\lambda-1}{2\lambda}\,x^{\frac{\lambda}{1-\lambda}}>0 .
\]

\subsection{Quadratic payoffs}

This example extends the previous one (which can be recovered by setting $\eta=0$) and is constructed so that the stopping boundary of the maximiser is again linear. Unlike in the previous example, now we find that $m^1_* \in (0,1)$.

Fix
\[ \lambda \in \left(\frac{1}{2}, \frac{2}{3}\right),  \qquad  \eta \in \left( 0, \frac{2\lambda(2 \lambda-1)}{3(2-3 \lambda)} \right);  \]
and define
\[ \alpha = \frac{2 \lambda - 1}{\lambda^2},   \qquad  \kappa = \eta \frac{2 - 3 \lambda}{\lambda^3} .  \]
Then $0 < \alpha < \frac{3}{4}$ and $0<\kappa < \frac{2 \alpha}{3}$.
Finally define $f,g:[0,1] \to \R$ by
\[ f(x)=x + \frac{\eta}{2}x^2 \hspace{10mm} g(x)= \alpha x - \frac{\kappa}{2} x^2.  \]
Note that parameters are such that $g$ is increasing.

We look for a solution $c_* = c_*(m)$ to  $m(m-c)f'(m) = c (f(m)-g(c))$, or equivalently $m(m-c)(1 + \eta m) = c(m + \frac{\eta}{2}m^2 - \alpha c + \frac{\kappa}{2} c^2)$. That is, $c_*$ solves $Q(c) = 0$, where
\[ Q(c) = \kappa c^3 - 2 \alpha c^2 + c(4m + 3 \eta m^2) - 2m^2(1 + \eta m).  \]
It is easy to verify that $Q(\lambda m)=0$, so $c_*(m) = \lambda m$. Note that ${\mathcal L}^1 g(x) = \frac{1}{2} g''(x) + \frac{1}{x} g'(x) = - \frac{\kappa}{2} + \frac{1}{x}(\alpha - \kappa x) = \frac{\alpha}{x} - \frac{3}{2} \kappa$ so that ${\mathcal L}^1 g(x) > 0$ on $(0,1)$ since $\alpha > \frac{3}{2}\kappa$.
Conditions (iii)–(iv) of Assumption \ref{ass:fghscaled} hold with $x_g=1$, because $\cL^0g(x)=-\frac{\kappa}{2}-\frac{\alpha-\kappa x}{1-x}<0$ and $\cL^1f(x)=\frac{\eta}{2}+\frac{1+\eta x}{x}>0$. By Lemma~\ref{lem:jan}, condition (v) is satisfied and in particular, the root of \eqref{eq:phidef} is unique in $(0,m)$. The first inequality in (vi) holds because $c_*(1)g'(c_*(1))=\lambda(\alpha-\kappa\lambda)<\lambda\alpha=\frac{2\lambda-1}{\lambda}<1<f'(1)$. The second inequality in (vi) follows from Lemma \ref{lem:newsuff}, since $g$ is concave.
Hence all the conditions of Assumption \ref{ass:fghscaled} hold and $c_*(m) = \lambda m$.

Solving \eqref{eq:Ddef} we find that
\begin{align*}
\Du0_{c_*}(m) & = \frac{m^{\frac{1}{1-\lambda}}}{1-m} \left( - \alpha + \frac{\lambda^2}{2 \lambda - 1} \left( \alpha + \frac{\kappa}{2} \right) + \frac{\lambda^3 \kappa}{2(2-3 \lambda)} \right) \\
& \quad - \frac{m}{(1-m)} \left( - \alpha + \frac{\lambda^2}{2 \lambda - 1} \left(\alpha + \frac{\kappa}{2} \right) m + \frac{\lambda^3 \kappa}{2(2-3 \lambda)} m^2 \right).
\end{align*}
After some more algebra we find that
\begin{equation}
\gamma_{c_*}(m)  =  \frac{1 +  \frac{(2 - 3 \lambda)}{2 \lambda(2 \lambda -1)} \eta + \frac{3}{2}  \eta m -  \left( \frac{1}{2 \lambda} + \eta \frac{1}{2(2 \lambda - 1)}  \right) m^{\frac{2 \lambda - 1}{1 - \lambda}} }{1 +  \frac{3}{2} \eta m }.
\end{equation}
Then $\gamma_{c_*}(0) = 1 + \eta \frac{(2 - 3 \lambda)}{2 \lambda(2 \lambda -1)}>1$ and $\gamma_{c_*}(1) = \frac{1 - \frac{1}{2 \lambda} -  \eta \frac{(2 - 3 \lambda)}{2 \lambda} }{1 + \frac{3}{2}  \eta} \in (0,1)$. Finally,
\[ m^1_* = \left( \frac{\eta(2 - 3 \lambda)}{(2 \lambda-1) + \eta \lambda} \right)^{\frac{1-\lambda}{2 \lambda - 1}} \in (0,1). \]
Note that $d_*(m)= c_*(m)=\lambda m$ on $[m^1_*,1]$ and $d_*(m) = \lambda m^1_* \wedge m$ on $(0,m^1_*)$ and the optimal strategy for the maximiser is $\sigma_* = \inf \{ u\ge 0 \colon X_u \leq d_*(M_u) \} \wedge \He$.

\textbf{The equilibrium.} Start from $(X_0,p_{0-})=(x,x)$ and set
$M_t\coloneqq(x\vee m^1_*)\vee\sup_{s\le t}X_s$. For every $x\in(0,1)$ the
equilibrium reads as follows. Type $0$ of Player 1 never stops. Type $1$ has the generating process $\xi_t=1-\Gamma_*(M_t)$, $t<\He$, $\xi_\He = 1$, with
\[
\Gamma_*(m)=\ind_{\{m\le m^1_*\}}+\gamma_{c_*}(m)\ind_{\{m>m^1_*\}}.
\]
She stops at time $0$ with probability $1-\Gamma_*(x\vee m^1_*)$, and afterwards only at times of increase of $M$ according to the rule
\[
\inf\big\{t\ge 0:1-\Gamma_*(M_t)>U\big\}\wedge\He.
\] 
Player 2 stops at time $\sigma_*=\inf\{t\ge 0\colon X_t\le\lambda M_t\}\wedge\He$.

Three regimes emerge. If $x\le\lambda m^1_*$, Player~2 stops immediately and the payoff is $g(x)$. If $\lambda m^1_*<x\le m^1_*$, nobody acts and the belief is $p_t=X_t$ until either $X$ hits $\lambda m^1_*$, where Player~2 stops, or $X$ hits $m^1_*$. In the latter case play continues, with zero probability of stopping at $m^*_1$, because $\Gamma_*(m^1_*)=1$. If $x>m^1_*$, then $(x,x)\in A_1$: type~$1$ stops at once with probability $1-\Gamma_*(x)$, which lowers the belief to $\pi_*(x,x)<x$. After that, both players follow the strategies where the Player~2 stops at the boundary $\lambda M_t$ and type $1$ of Player 1 acts continuously only at reaching a new maximum.

\subsection{A non-linear stopping boundary}
In the previous examples the boundary $c_*$ is linear. To obtain an example with a curved boundary we fix
$\beta\in(0,1)$ and take
\[
g(x)=x,\qquad f(x)=\frac{1+\beta x+\beta^2x^2}{3\beta}.
\]
We verify Assumption~\ref{ass:fghscaled} explicitly. We have $f(x)-g(x)=\frac{(1-\beta x)^2}{3\beta}>0$ on $[0,1]$ and $f'(x)=\frac{1+2\beta x}{3}>0$.
Moreover $\cL^0g(x)=-\frac{1}{1-x}<0$ and $\cL^1f(x)=\frac{\beta}{3}+\frac{1+2\beta x}{3x}>0$,
and $g$ is concave with $x_g=1$. Since $\cL^1g(x)=\frac{1}{x}>0$, Lemma~\ref{lem:jan} gives
condition~(v), and we canverify $m^2f'(m)>c_*^2(m)g'(c_*(m))$ for all $m$, so that $m_0(c_*)=0$.
A direct calculation shows that
\[
c_*(m)=\beta m^2
\]
solves $c\,(f(m)-g(c))=m(m-c)f'(m)$; hence it is the unique root of $\phi(\cdot,m)$.
For~(vi), we have $c_*(1)g'(c_*(1))/f'(1)=3\beta/(1+2\beta)<1$, and the second inequality
follows from Lemma~\ref{lem:newsuff} because $g$ is linear.

The coefficient of \eqref{eq:Ddef} decomposes as
$\frac{1-\beta m^2}{(1-m)\,m(1-\beta m)}=\frac1m+\frac{1}{1-m}+\frac{\beta}{1-\beta m}$.
Hence the bounded solution is 
\[
\Du0_{c_*}(m)=\beta m^2+\frac{\beta m(1-m)(3-\beta-2\beta m)}{3(1-\beta m)},
\]
with $\Du0_{c_*}(1)=\beta=g(c_*(1))$. For $\beta m^2\le x\le m\le1$, \eqref{eq:u0*def} and
\eqref{eq:u1*def} read
\[
u^0_*(x,m)=\beta m^2+\frac{\beta(x-\beta m^2)(1-m)^2(3-\beta-2\beta m)}{3(1-x)(1-\beta m)^2},
\qquad
u^1_*(x,m)=\beta m^2+\frac{(x-\beta m^2)(1+2\beta m)}{3\beta x}.
\]
From \eqref{eq:gamma*def},
\[
\gamma_{c_*}(m)=\frac{\beta\big[3(1-\beta m^2)(1-\beta m)^2-\beta(1-m)^2(3-\beta-2\beta m)\big]}
{(1-\beta m)^3(1+3\beta m)} .
\]
The function $\gamma_{c_*}$ decreases from $\gamma_{c_*}(0+)=1-(1-\beta)^3$ to
$\gamma_{c_*}(1)=\frac{3\beta}{1+3\beta}$. Thus $m^1_*=0$, $\Gamma_*=\gamma_{c_*}$ on $(0,1]$ and
$d_*=c_*$.

\textbf{The equilibrium.} Starting from $(x,x)$, type $1$ of Player 1 stops at time $0$ with
probability $1-\gamma_{c_*}(x)$, and thereafter at times of increase of $M$ according to the rule $\inf\{t\ge 0: 1 - \gamma_{c_*}(M_t)>U\}\wedge \He$. Type $0$ never stops. Player~2 stops at
\[
\sigma_*=\inf\{t\ge0\colon X_t\le\beta M_t^2\}\wedge\He .
\]
In contrast with the previous examples, the initial atom in the stopping distribution of type $1$ does not vanish as $x\downarrow0$; it tends to $(1-\beta)^3$.

\appendix
\section{Proofs of lemmas from Section \ref{sec:problem}}\label{app0}

Here we provide proofs of Lemmas \ref{lem:jan}, \ref{lem:newsuff} and \ref{lem:c*}.

\begin{proof}[Proof of Lemma \ref{lem:jan}]
Let $c=c(m)$ be such that $\phi(c(m),m)=0$. By Remark \ref{rem:c*} we know that $c(m)$ exists and it lies in $(0,m)$.
If $m>x_g$ then by the assumption that $\phi(x_g,m)>0$ we conclude that $c$ may be taken to lie in $(0, x_g \wedge m)$.
On $(0,m)$ we have
\[ \partial_x \phi(x,m) = \frac{1}{(m-x)^2}\psi(x,m)\]
and
\[ \partial_{xx} \phi(x,m) = \frac{1}{(m-x)^3}(2 \psi(x, m) + (m-x) \partial_x \psi(x,m)),\]
where $\psi(x,m) = mf(m)-mg(x)- x(m-x)g'(x)$. Assume that there are at least two roots of $\phi(\cdot, m)$ on $(0, x_g \wedge m)$. Given that $\phi(0, m) < 0$ and $\lim_{x \uparrow x_g \wedge m} \phi(x,m) > 0$, there must be $x_{max},x_{min} \in (0, x_g \wedge m)$ at which $\phi(\cdot, m)$ attains a local maximum and a local minimum, respectively, and $\phi(x_{max}, m) \ge 0$, $\phi(x_{min},m)\le 0$. Suppose $\cL^1_X g >0$ on $(0,x_g]$. Then, $\partial_{x} \phi(x_{min}, m) = 0$ and $\partial_{xx} \phi(x_{min}, m) \ge 0$ implying $\psi(x_{min},m) = 0$ and $(m-x_{min}) \partial_x \psi(x_{min},m) \ge 0$. However,
\begin{equation}\label{eqn:partial_x_psi}
\begin{aligned}
\partial_x \psi\big(x_{min},m\big) &= -  2x\big(m-x_{min}\big) \left[ \frac{1}{2}g''\big(x_{min}\big) + \frac{1}{x_{min}} g'\big(x_{min}\big) \right]\\
 &= - 2x_{min}\big(m-x_{min}\big) \cL^1_X g\big(x_{min}\big) < 0,
\end{aligned}
\end{equation}
hence a contradiction. Analogously, suppose $\cL^1_X g <0$ on $(0,x_g]$. Then, $\partial_{x} \phi(x_{max}, m) = 0$ and $\partial_{xx} \phi(x_{max}, m) \le 0$ implying $\psi(x_{max},m) = 0$ and $(m-x_{max}) \partial_x \psi(x_{max},m) \le 0$. However,
\begin{equation}
\partial_x \psi\big(x_{max},m\big) = - 2x\big(m-x_{max}\big) \cL^1_X g\big(x_{max}\big) < 0,
\end{equation}
hence a contradiction. This shows that the root is unique on $(0, x_g \wedge m)$.
\end{proof}

\begin{proof}[Proof of Lemma \ref{lem:newsuff}]
In this proof we write $c_*$ as shorthand for $c_*(1)$. By definition of $c_*$ and the first inequality in \eqref{eq:newcondition} we get
\[
c_* \frac{f(1)-g(c_*)}{1-c_*}=f'(1)>c_* g'(c_*).
\]
Then $f(1)> g(c_*) + (1- c_*)g'(c_*)$.
Notice that \eqref{eq:newsufficient} is equivalent to
\[
g(c_*) + (1-c_*)g'(c_*) \geq g(x) + (1-x) \frac{g(x)-g(c_*)}{x-c_*},\quad \quad x\in(c_*,1),
\]
and therefore
\[
f(1)> g(x) + (1-x) \frac{g(x)-g(c_*)}{x-c_*}, \quad  \quad x\in(c_*,1).
\]
The latter can be rewritten as $f(1)( (1- c_*) - (1-x)) > g(x)(1-c_*) - g(c_*)(1-x)$ or equivalently
\[
\frac{f(1)-g(x)}{1-x} > \frac{f(1)-g(c_*)}{1-c_*} = \frac{f'(1)}{c_*},\quad x\in(c_*,1),
\]
where the equality uses again the definition of $c_*$.
Then, using \eqref{eq:newsufficient} again,
\[ \frac{f(1)-g(x)}{1-x} > \frac{f'(1)}{c_*} \geq \frac{f'(1)}{c_*} \frac{1}{g'(c_*)} \frac{g(x)-g(c_*)}{x-c_*} \]
which is sufficient for the second inequality in \eqref{eq:newcondition}.
\end{proof}

\begin{proof}[Proof of Lemma \ref{lem:c*}]
We have that $c_*(1)<1$ from Remark \ref{rem:c*}.
Take any sequence $(m_j)_{j\in\N}\subset(0,1)$ converging to $m\in(0,1]$. Then, $\limsup_{j\to\infty}c_*(m_j)\in(0,m \wedge x_g]$ and
$\liminf_{j\to\infty}c_*(m_j)\in(0,m \wedge x_g]$. Moreover, by the joint continuity of $(x,m)\mapsto {\oldh}(x,m)$ and
$\oldh(c_*(m_j),m_j)=0$ we have
\[
\oldh\big(\limsup_{j\to\infty}c_*(m_j),m\big)=\oldh\big(\liminf_{j\to\infty}c_*(m_j),m\big)=0.
\]
Note that since $\oldh(0, m) = -m f'(m) \ne 0$ we must have $\liminf_{j\to\infty}c_*(m_j) > 0$.
By uniqueness of the root in $(0,m \wedge x_g)$ we deduce $c_*(m)=\limsup_{j\to\infty}c_*(m_j)=\liminf_{j\to\infty}c_*(m_j)$. Then, $c_*$ is continuous on $(0,1]$.
As $0< c_*(m) < m$ for $m \in (0,1)$ it follows that $c_*(0)\coloneqq\lim_{m\to 0}c_*(m) = 0$.

Using that
\[
\partial_x\phi(c_*(m),m) = \frac{\big(m^2 f'(m) - c_*(m)^2 g'(c_*(m))\big)}{c_*(m)(m-c_*(m))},
\]
on the set $m^2 f'(m)> c_*(m)^2 g'(c_*(m))$ we have $\partial_x\phi(c_*(m),m)>0$. Then $c_*$ is differentiable by the implicit function theorem with
\begin{equation}
\label{eq:c'}
 c_*'(m) = - \frac{\partial_m \phi(c_*(m),m)}{\partial_x\phi(c_*(m),m)} = \frac{2m c_*(m) (m-c_*(m))}{m^2 f'(m) - c_*(m)^2 g'(c_*(m))} \cL^1 f(m) > 0
\end{equation}
where the final inequality holds by Assumption~\ref{ass:fghscaled}-(iv), which gives us that $\cL^1 f>0$. This proves the final claim in the lemma.
\end{proof}

\section{The case of linear diffusions}\label{app}

The aim of this appendix is to outline the generalisation of our framework to the case when $X$ is replaced by a 1-dimensional, time-homogeneous diffusion (so-called linear diffusions). The basic idea is standard: a linear diffusion can be obtained as a suitable scaling of a time-changed Brownian motion. Because this is a textbook approach we only outline the procedure and provide formulae that allow to immediately map the general case to our Brownian setting.

In this section the tuple $(\Omega,\cF,\G,(Y_t)_{t\ge 0},(W_t)_{t\ge 0},(\P_y)_{y\in\R})$ defines a one-dimensional time-ho\-mo\-ge\-neous
diffusion with dynamics described by
\begin{equation}\label{eq:Y}
Y_t =Y_0 +\int_0^t a (Y_s) \ud W_s +\int_0^t b(Y_s) \ud s, \quad t\ge 0,
\end{equation}
where $W$ is a standard Brownian motion under all measures $\P_y$, the coefficients $a$ and $b$ are assumed continuous on $\R$ and $a(y)>0$ on a bounded interval $[\underline y,\overline
y]$. As usual $\P_y(Y_0=y)=1$ and $\G=(\cG_t)_{t\ge 0}$ is a filtration with respect to which the pair $(Y,W)$ is adapted.

We assume
$y\in(\underline y,\overline y)$ and $Y$ runs until it hits either $\underline{y}$ or $\overline{y}$. Letting
$\HY_z =\inf\{t\ge 0:
Y_t=z\}$ for $z\in[\underline y,\overline y]$, we set $\HY_{\e}:=\HY_{\underline y}\wedge \HY_{\overline y}$.
The filtration $\G^Y\coloneqq(\cG^Y_t)_{t\ge 0}$ is generated by the process $Y$ and $\theta^Y :=  Y_{\HYe}$. We set $\cG^2_t=\cG^Y_t$ and $\cG^1_t=\cG^Y_t\vee\sigma(\theta^Y)$ and denote $\G^k=(\cG^k_t)_{t\ge 0}$, $k=1,2$. The less informed player (maximiser, Player 2) picks a randomised stopping time $\sigma^Y \in\cT_R(\G^2)$, whereas Player 1 (minimiser) picks $\tau^Y\in\cT_R(\G^1)$. At
time $\tau^Y\wedge \sigma^Y$ the game ends and the following payoff is paid by Player 1 to Player 2:
\begin{equation}
\cP^Y(\sigma^Y,\tau^Y)\coloneqq
F\big(Y_{\tau^Y}\big)\ind_{\{\tau^Y<\sigma^Y\}}+G\big(Y_{\sigma^Y}\big)\ind_{\{\sigma^Y<\tau^Y\}}+H\big(Y_{\tau^Y}\big)
\ind_{\{\tau^Y=\sigma^Y\}}.
\end{equation}
Let $(\underline y,\overline y)\ni y\mapsto s(y)$ be the scale function associated to $Y$. The scale function is defined up to an affine transformation and therefore we normalise it so as to guarantee $s(\underline y)=0$ and $s(\overline y)=1$. In analogy with \eqref{eq:tau} $\tau^Y\in \cT_R(\G^1)$ is decomposed as
\[
\tau^Y=\ind_{\{\theta^Y=\underline y\}}\tau^Y_0+\ind_{\{\theta^Y=\overline y\}}\tau^Y_1,
\]
where $\tau^Y_0,\tau^Y_1\in\cT_R(\G^2)$. Then for any pair $(\sigma^Y,\tau^Y)\in\cT_R(\G^2)\times\cT_R(\G^1)$ we have
\begin{equation*}
\E_y\big[\cP^Y(\sigma^Y,\tau^Y)\big]=\big(1-s(y)\big)\E^0_y\big[\cP^Y\big(\sigma^Y,\tau^Y_0\big)\big]+s(y)\E^1_y\big[\cP^Y\big(\sigma^Y,\tau^Y_1\big)\big],
\end{equation*}
where we used the changes of measure prescribed by
\begin{equation}
\frac{\ud \P^0_y}{\ud\P_y}\Big|_{\cF^2_t}:=\frac{1-s(Y_t)}{1-s(y)}\quad\text{and}\quad \frac{\ud
\P^1_y}{\ud\P_y}\Big|_{\cF^2_t}:=\frac{s(Y_t)}{s(y)},\quad t\le \HYe,\quad \P_y-a.s.
\end{equation}
The measure $\P^0_y$ is the one under which the process $Y$ leaves the interval $(\underline y,\overline y)$
crossing $\underline y$ and $\P^1_y$ the one under which $Y$ leaves the interval crossing $\overline y$.
By the same arguments using Girsanov's theorem as in the Brownian formulation,
the process
\[
W^0_t=W_t+\int_0^t\frac{a(Y_u)s'(Y_u)}{1-s(Y_u)}\ud u,
\]
is a $\P^0_y$-Brownian motion and the dynamics of $Y$ under $\P^0_y$ reads as
\[
\ud Y_t= a(Y_t)\ud W^0_t + \Big(b(Y_t)-\frac{a^2(Y_t)s'(Y_t)}{1-s(Y_t)}\Big)\ud t.
\]
Similarly, the process
\[
W^1_t=W_t-\int_0^t\frac{a(Y_u)s'(Y_u)}{s(Y_u)}\ud u,
\]
is a $\P^1_y$-Brownian motion and the dynamics of $Y$ under $\P^1_y$ reads as
\[
\ud Y_t=a(Y_t)\ud W^1_t+\Big(b(Y_t)+\frac{a^2(Y_t)s'(Y_t)}{s(Y_t)}\Big)\ud t.
\]

Setting $Z=s(Y)$ and $\alpha(z)=(a\cdot s')\circ s^{-1}(z)$, by It\^o's formula we get
\begin{equation*}
\begin{aligned}
\ud Z_t=s'(Y_t)\ud Y_t+\tfrac{1}{2} s''(Y_t)a^2(Y_t)\ud t= \alpha(Z_t)\ud W_t,
\end{aligned}
\end{equation*}
where $W$ is the Brownian motion under $\P_y$. It is now clear that the functions $F$, $G$, $H$, are related to $f$, $g$, $h$ from \eqref{eq:P2} via
\[
f=F\circ s^{-1},\,\,g=G\circ s^{-1},\,\, h=H\circ s^{-1},
\]
where $s^{-1}$ is the inverse of the scale function. Similarly, the upper and lower value of the game are denoted by
\begin{equation}
\begin{aligned}
&\underline V(y)\coloneqq\sup_{\sigma^Y\in\cT_R(\G^2)}\inf_{\tau^Y\in\cT_R(\G^1)}\E_y\big[\cP^Y(\sigma^Y,\tau^Y)\big]\\
&\text{and}\quad
\overline V(y)\coloneqq\inf_{\tau^Y\in\cT_R(\G^1)}\sup_{\sigma^Y\in\cT_R(\G^2)}\E_y\big[\cP^Y(\sigma^Y,\tau^Y)\big],
\end{aligned}
\end{equation}
and the value reads as $V(y)$ when $\underline V=\overline V$. These functions are related to the ones from the Brownian case via
\[
\underline v=\underline V\circ s^{-1},\,\,\overline v=\overline V\circ s^{-1},\,\, v=V\circ s^{-1}.
\]

It only remains to run the process $Z$ according to a different clock in order to complete the analogy with the problem we solved in the Brownian framework. In terms of the $\P^0_y$-Brownian motion $W^0$ and of the $\P^1_y$-Brownian motion $W^1$ the dynamics read as
\begin{equation*}
\begin{aligned}
\ud Z_t=\alpha(Z_t)\ud W^0_t-\frac{\alpha^2(Z_t)}{1-Z_t}\ud t\quad\text{and}\quad \ud Z_t= \alpha(Z_t)\ud W^1_t+\frac{\alpha^2(Z_t)}{Z_t}\ud t.
\end{aligned}
\end{equation*}
Since $a>0$ on $[\underline y,\overline y]$, we can perform a time-change using the process $A(t)=\int_0^t \alpha^2(Z_s)\ud s$. Letting $T(t)=(A^{-1})(t)$ be the inverse of $A$, defining $X_t=Z_{T(t)}$ and setting $x=s(y)$, we obtain
\begin{equation}\label{eq:SDEX}
\begin{aligned}
X_t=x+B_t,\quad X_t&=x+B^0_{t}-\int_0^t\frac{1}{1-X_s}\ud s,\quad X_t=x+B^1_{t}+\int_0^t\frac{1}{X_s}\ud s,
\end{aligned}
\end{equation}
where $B$, $B^0$ and $B^1$ are $\P_y$, $\P^0_y$ and $\P^1_y$ Brownian motions, respectively, for the time changed filtration $\F=(\cG^Y_{T(t)})_{t\ge 0}$. For $\F^2=\F$ and $\F^1=(\cF^1_t)_{t\ge 0}$ with $\cF^1_t=\cF^2_t\vee\sigma(\theta^X)$, where $\theta^X=X_{\He}$ is defined as in \eqref{eq:thetaX}, we recover precisely the setup from Section \ref{sec:problem} upon observing that optimising over $\cT_R(\G^1)$ and $\cT_R(\G^2)$ is equivalent to optimising over $\cT_R(\F^1)$ and $\cT_R(\F^2)$ because
\[
\cT_R(\G^1)\times\cT_R(\G^2)\ni(\tau,\sigma)\mapsto \big(A(\tau),A(\sigma)\big)\in \cT_R(\F^1)\times\cT_R(\F^2),
\]
is a bijection.
The assumptions for functions $F$, $G$ and $H$ are easily derived from Assumption \ref{ass:fghscaled}. We will use the infinitesimal generators $\cL_Y$, $\cL^0_Y$ and $\cL^1_Y$
of the process $Y$ under the measures $\P_y$, $\P^0_y$ and $\P^1_y$, respectively.
\begin{assumption}
Functions $F$, $G$ and $H$ are such that:
\begin{itemize}
\item[(i)] on $[\underline y,\overline y]$ we have $F\ge H\ge G$ with $F,G$ continuous and $H$ measurable;
\item[(ii)] $H(\underline y)=G(\underline y)$, $H(\overline y)=F(\overline y)$;
\item[(iii)] on $(\underline y,\overline y]$ functions $F$ and $G$ are twice continuously differentiable with $F>G$, $F'>0$, $G'>0$;
\item[(iv)] $\cL^0_Y G(y)< 0$ and $\cL^1_Y F(y)> 0$ for $y\in (\underline y,\overline y)$ and $\cL_Y G(y)\le 0$ on $(\underline{y},y_G]$ for some $y_G\in(\underline y,\overline y]$.
\item[(v)] For $\overline y > \nu > \underline y$, the map
\begin{equation}
[\underline y, \nu) \ni y\mapsto \oldtlh(y,\nu)\coloneqq
s(y)\frac{F(\nu)-G(y)}{s(\nu)-s(y)}-s(\nu)\frac{F'(\nu)}{s'(\nu)},
\end{equation}
admits a root $y_*(\nu)\in (\underline y,y_G\wedge \nu)$ and it is unique in $(\underline y,y_G\wedge \nu]$;
\item[(vi)] it holds that
\begin{equation}
\label{eq:newconditionY}
1>s\big(y_*(\overline y)\big)\frac{G'\big(y_*(\overline y)\big)}{s'\big(y_*(\overline
y)\big)}\frac{s'(\overline y)}{F'(\overline y)} \geq \frac{G(y)-G\big(y_*(\overline y)\big)}{s(y)-s\big(y_*(\overline y)\big)}\frac{s(\overline y)-s(y)}{F(\overline y) - G(y)},\quad  \forall y \in \big(y_*(\overline y),\overline y\big).
\end{equation}
\end{itemize}
\end{assumption}

For the equivalence with Assumption \ref{ass:fghscaled} it must be $x_g=s(y_G)$, $c_*(m)=s\circ y_*\circ s^{-1}(m)$. The proof of such equivalence is a matter of straightforward but tedious algebra and we omit it. Thanks to the established equivalence between the Brownian framework and the general diffusive one, we can transfer the statement from Theorem \ref{thm:mainscaled} to the general diffusive framework. For that, it is convenient to introduce some new notations.

We define functions
\[
W(y,p)\coloneqq V\big(s(y),p\big),\ w_*(y,m)\coloneqq v_*\big(s(y),s(m)\big),\ \vartheta(y,p)\coloneqq\Psi\big(s(y),p\big),\ I_*(y)\coloneqq \Gamma_*\big(s(y)\big),
\]
and the rescaled belief map for $y\le \ell$
\[
\beta_*(y,\ell)\coloneqq\frac{s(y)I_*(\ell)}{1-s(y)+s(y)I_*(\ell)}.
\]

Recall the sets $A_1$ and $A_2$ from Section \ref{sec:candidatePBE} and the function $\Gamma_*$ from Theorem \ref{thm:mainscaled}. Let $M^Y_t=M^Y_0\vee \sup_{s\le t} Y_s$ with
\[
M^Y_0= s^{-1}\circ\Gamma^{-1}_*\Big(\vartheta(y,p)\Big),\quad\text{if}\ \big(s(y),p\big)\in A_2,
\]
and
\[
M^Y_0 = s^{-1}\circ\Gamma^{-1}_*\Big(\vartheta\big(y,\beta_*(y,y)\big)\Big) =y,\quad\text{if}\ \big(s(y),p\big)\in A_1.
\]
In equilibrium the belief process shall read as
\[
p_t=\pi_*\big(s(Y_t),s(M^Y_t)\big)=\beta_*(Y_t,M^Y_t),\quad t\in[0,\HYe).
\]
The generating process of the $\theta^Y=1$ type of Player 1 reads, for $\big(s(y),p\big)\in A_2$
\[
\Xi_t=\Xi^p_t=1-I_*(M^Y_t)/I_*(M^Y_0),\quad t\in[0,\He^Y),\quad \text{and}\ \Xi_{\HYe}=1.
\]
Similarly, for $\big(s(y),p\big)\in A_1$ it reads as
\[
\Xi_t=\Xi^p_t=\Delta\Xi_0+\big(1-\Delta\Xi_0\big)\Big(1-I_*(M^Y_t)/I_*(M^Y_0)\Big),\quad t\in[0,\He^Y),\quad \text{and}\ \Xi_{\HYe}=1,
\]
with $\Delta\Xi_0=1-I_*(y)$.
Then, we set
\begin{equation*}
\begin{aligned}
\sigma^Y_* =\sigma^{Y,p}&\coloneqq\inf\big\{t\ge 0: Y_t\le y_*\big(M^Y_t \vee y^1_{I_*}\big)\big\}\wedge\HYe\quad\text{and}\quad
\tau^Y_*=\tau^{Y,p}\coloneqq \ind_{\{\theta^Y= \underline{y} \}}\HYe+\ind_{\{\theta^Y= \overline{y} \}}\tau^{Y,p}_1,
\end{aligned}
\end{equation*}
with $y^1_{I_*} = \sup \{m \in [\underline y,\overline y) \colon I_*(m) = 1\} = s^{-1}(m^1_*)$ and
$\tau^{Y,p}_1\coloneqq\inf\big\{t\ge 0\colon \Xi^p_t>U\big\}$.
Finally, define $\hP_{y}^{p} \coloneqq p \P^1_y + (1-p) \P^0_y$.
The notation above allows us to restate Theorem \ref{thm:mainscaled} in terms of the process $Y$ and its running maximum. The proof is omitted as it is just a rescaling of the formulae in Theorem \ref{thm:mainscaled}.
\begin{theorem}\label{thm:mainunscaled}
The pair $(\sigma^Y_*, \tau^Y_*)$ is a perfect Bayesian equilibrium for the game with payoff $\cP^Y$. That is, for any
$(s(y),p) \in A_1 \cup A_2$ and any pair $(\sigma,\tau)\in\cT_R(\G^2)\times\cT_R(\G^1)$
\[
\hE_{y}^{p}\big[\cP^Y(\sigma,\tau^Y_*)\big]\le \hE_{y}^{p}\big[\cP^Y(\sigma^Y_*,\tau^Y_*)\big]\le \hE_{y}^{p}\big[\cP^Y(\sigma^Y_*,\tau)\big].
\]
Moreover, the value of the game is equal to
\begin{equation}\label{eqn:Vxpunscaled}
W(y,p) =
\begin{cases}
s(y) F(y) \big[1-I_*(y)\big] + \big[1-s(y) + s(y) I_*(y)\big] w_*\big(y,s(y)\big), & \big(s(y),p\big) \in A_1,\\
w_*\big(y, s(I^{-1}_*(\vartheta(y,p)))\big), & \big(s(y),p\big) \in A_2.
\end{cases}
\end{equation}
\end{theorem}

\bibliographystyle{abbrv}
\bibliography{biblio}

\end{document}